%% file: DAEv4.tex
\documentclass[hidelinks,onefignum,onetabnum]{siamart250211}

\input{ex_shared}
\newtheorem{assumption}{Assumption}
\usepackage{mathtools}
\usepackage{microtype}
\usepackage{booktabs}
\usepackage{multirow}
\usepackage{graphicx}
\usepackage{epstopdf}
\usepackage{subfigure}
\usepackage{array}
\usepackage{diagbox}
\usepackage{verbatim}
\ifpdf
\hypersetup{
  pdftitle={Deep asymptotic expansion method},
  pdfauthor={Qiao Zhu, Dmitrii Chaikovskii, Bangti Jin, and Ye Zhang}
}
\fi

\begin{document}

\maketitle

\begin{abstract}
Physics-informed neural network (PINN) has shown great potential in solving partial differential equations. However, it faces challenges when dealing with problems involving steep gradients. The solutions to singularly perturbed time-dependent reaction-advection-diffusion equations exhibit internal moving transition layers with sharp gradients, and thus  the standard PINN becomes ineffective. In this work, we propose a deep asymptotic expansion (DAE) method, which is inspired by asymptotic analysis and leverages deep learning to approximate the smooth part of the expansion. We first derive the governing equations for transition layers, which are then solved using PINN. Numerical experiments show that the DAE outperforms the standard PINN, gPINN, and PINN with adaptive sampling. We also show its robustness with respect to training point distributions, network architectures, and random seeds.
\end{abstract}



\section{Introduction}
Let $\overline{D}_{\boldsymbol{x}}=[-a,a]\times \mathbb{R}^{d-1} \subset\mathbb{R}^d$ with $\boldsymbol{x}=(x_1,...,x_d)$
and $x_1\in [-a,a]$, and define $
\boldsymbol{x}^*= (x_2,\dots,x_d)$  if $ d \geq 2$, and $\boldsymbol{x}^*=
\emptyset$ if $d =1$. Consider the following reaction-advection-diffusion 
equation:
{\small
\begin{align} \label{mainproblem}
\begin{cases}
\displaystyle \mu \Delta u -\partial_t u=A(u,\boldsymbol{x}) ( \boldsymbol{1}_d \cdot \nabla u)+ f(\boldsymbol{x}), \quad \mbox{in } D_{\boldsymbol{x}}\times (0, T],\\
\displaystyle u(-a,\boldsymbol{x}^*, t)=L(\boldsymbol{x}^*), ~  u(a,\boldsymbol{x}^*, t)=R(\boldsymbol{x}^*), ~ u(\boldsymbol{x}+\boldsymbol{P},t) = u(\boldsymbol{x},t), ~ t\in \bar{\mathcal{T}}=[0, T],\\
u(\boldsymbol{x},0)={u}_0(\boldsymbol{x}, \mu), ~ \boldsymbol{x}\in \overline{D}_{\boldsymbol{x}},
\end{cases}
\end{align}
}where $0<\mu \ll 1$ is a small parameter, $A(u,\boldsymbol{x})$ and $f(\boldsymbol{x})$ are smooth functions. The notation $\Delta$ and $\nabla$ denote the Laplacian and gradient operators, respectively, and $\boldsymbol{1}_d\in \mathbb{R}^d$ is a  row vector of ones.
 $\boldsymbol{P} = (0, P_2, \ldots, P_d)$ is the period length in each direction, and the periodicity is imposed on all axes except $x_1$.

The model \eqref{mainproblem} is singularly perturbed due to the presence of a small parameter $\mu>0$ in front of the highest-order derivative. Its solution exhibits steep gradients within narrow regions, known as internal transition layers. This model has been widely applied in mathematical modeling, e.g., density distributions of liquids and gasses under spatial heterogeneity \cite{olchev2009application,sogachev2006modification}, nonlinear acoustic phenomena \cite{nefedov2019existence,rudenko2017inhomogeneous}, temperature in the near-surface layer
of the ocean \cite{levashova2017heat}, carrier wave functions in Si/Ge heterostructures \cite{orlov2015use}, and propagation of an autowave front in a medium with barriers \cite{levashova2019spatio,sidorova2019autowave}.

Conventional numerical methods for solving the model \eqref{mainproblem}, such as finite difference, finite element, and finite volume methods, struggle to resolve sharp transitions in solutions as $\mu \to 0$. In contrast, asymptotic methods decompose the solution into different scales and typically become more accurate as $\mu$ decreases. {Nonetheless the subproblems in asymptotic methods generally do not have closed-form solutions, especially on complex geometries or with nonlinear coefficients. In practice, one still applies standard discretization schemes to resolve the subproblems, which is mesh dependent and can be costly. This practical limitation motivates our work.}

In recent years, deep learning methods have emerged as a versatile framework for solving PDEs \cite{han2018solving,sirignano2018dgm,zang2020weak}. Among these, physics-informed neural network (PINN) \cite{raissi2019physics} integrates physical laws directly into the learning process, using automatic differentiation to provide a mesh-free approach that avoids discretization errors. PINN has been widely applied to forward and inverse PDE problems \cite{CenJinLiZhou:2025,guo2022monte,JinLiQuanZhou:2024,tang2023pinns,wang2021understanding}. However, standard PINN often fails to resolve steep gradients present in thin boundary or interior layers of singularly perturbed problems (SPPs), due to spectral bias in learning multiscale features \cite{wang2021eigenvector}. Several enhanced PINN-based methods have been proposed for SPPs, which can be broadly classified into three categories.

\begin{itemize}

\item {Improving architecture and embedding priors. This line of research integrates physical insights into the model structure. One idea is architectural decomposition, which employs separate subnetworks for different scales~\cite{zhang2024multi}. Inspired by classical perturbation theory and asymptotic expansions, Arzani et al.~\cite{arzani2023theory} designed parallel PINNs to represent different approximation orders for boundary-layer problems in the inner and outer regions. The other popular idea is trial-space enrichment, which augments the solution space to better capture sharp features or singularities \cite{ainsworth2022galerkin,ainsworth2025extended,gie2024semi,HuJinZhou:2024,wang2024chien}.}

\item {Optimizing the training strategy. This approach modifies the training process to improve convergence, including parameter-continuation schedules that train the model by gradually decreasing the governing parameter~\cite{cao2023physics,cao2024multistep}, curriculum-based reweighting of the loss function to handle difficult regions~\cite{wang2023less}, and gradient-boosting ensembles for sequential error correction~\cite{fang2024ensemble}.}

\item {Reformulating the learning objective. This approach redefines the loss function to handle steep gradients more effectively, including augmenting the loss with gradients of the PDE residual~\cite{yu2022gradient} and predicting the derivative fields of the solution, followed by recovering the solution via integration~\cite{cong2025solving}.}

\end{itemize}

{Despite impressive advances in PINN for SPPs, the majority of existing studies have focused on steady-state problems with stationary boundary or internal layers. Time-dependent SPPs with moving transition layers exhibit strongly coupled spatiotemporal multiscale behavior, where the interfaces evolve with time and interact with the outer solution. These characteristics induce significant stiffness and nonstationarity and thus pose significant challenges to conventional numerical solvers and learning-based approaches. The idea of integrating asymptotic analysis with deep learning remains largely underexplored due to significant theoretical and computational challenges. In this work, we address this gap by developing a unified asymptotics-informed learning framework for time-dependent SPPs with a general nonlinear advection coefficient \( A(u,\boldsymbol{x}) \). It substantially  extends relevant theoretical developments for the linear case $A(u,\boldsymbol{x})=-u$ \cite{chaikovskii2023asymptotic,chaikovskii2022convergence,chaikovskii2023solving}. First, we conduct an asymptotic analysis of model~\eqref{mainproblem} and derive the representation of the zero-order asymptotic solution. Then we determine the transition layer position within a PINN framework by minimizing the residual of the layer equation, which avoids the mesh dependence. We refer to the proposed solution strategy as the deep asymptotic expansion (DAE) method. The hybrid approach preserves the strengths of asymptotic methods for singularly perturbed problems. To further enhance the accuracy, we couple the PINN for the transition-layer equation with residual-based adaptive refinement (RAR)~\cite{lu2021deepxde}, which adaptively enriches the collocation set at points with the largest residuals. The main contributions of this work are threefold:
\begin{itemize}
  \item We establish a unified asymptotic theory for a class of time-dependent reaction-advection-diffusion equations, extending existing results to include a general nonlinear advection coefficient \( A(u,\boldsymbol{x}) \).
  \item By combining asymptotic analysis with machine learning-based techniques, we developed the DAE, an accurate and efficient algorithm for solving singularly perturbed time-dependent reaction-advection-diffusion equations.
  \item 
  Compared with existing PINN-based methods, the DAE achieves much higher accuracy at a reduced training cost.
\end{itemize}
}

The rest of the paper is organized as follows. In Section \ref{sec2}, we conduct the asymptotic analysis of the model \eqref{mainproblem}. In Section \ref{sec3}, we propose the deep asymptotic expansion method. In Section \ref{sec4} we present numerical examples to illustrate the DAE. Finally, we conclude the work with further discussions.

\section{Asymptotic analysis}\label{sec2}

We denote the position of the transition layer at each time $t$ by $ h(\boldsymbol{x}^*,t,\mu)$,  which divides the region $\overline{D}_{\boldsymbol{x}}$ into two parts: $\overline{D}_{\boldsymbol{x}}^{(-)} =\{\boldsymbol{x}: x_1\in [-a,h(\boldsymbol{x}^*, t,\mu) ]$, $\boldsymbol{x}^* \in \mathbb{R}^{d-1} \}$ and $   \overline{D}_{\boldsymbol{x}}^{(+)}=\{\boldsymbol{x}: x_1\in [h(\boldsymbol{x}^*, t,\mu),a]$, $\boldsymbol{x}^* \in   \mathbb{R}^{d-1}\}$. First, we list the standard assumptions.

\begin{assumption}
The following conditions hold.\label{assumption}
\begin{itemize}
    \item[(a)] The equation $A(u, \boldsymbol{x}) ( \mathbf{1}_d\cdot \nabla u)+f(\boldsymbol{x})=0$ has a solution  $\varphi^{(-)}(\boldsymbol{x})<0$ for $x\in \overline{D}_{\boldsymbol{x}}^{(-)}$ (respectively $\varphi^{(+)}(\boldsymbol{x})>0$ for $x\in \overline{D}_{\boldsymbol{x}}^{(+)}$) under the condition $u(-a,\boldsymbol{x}^*)=L(\boldsymbol{x}^*)$ (respectively $u(a,\boldsymbol{x}^*)=R(\boldsymbol{x}^*)$), and $\varphi^{(\mp)}(\boldsymbol{x})\in C^1(\overline{D}_{\boldsymbol{x}}^{(\mp)})$. Moreover,   $\varphi^{(+)}(\boldsymbol{x} ) - \varphi^{(-)}(\boldsymbol{x}) >  2\mu^\rho$ for some fixed $\rho\in(0,1)$.

\item[(b)] $A(\varphi^{(-)}(\boldsymbol{x}), \boldsymbol{x}) > 0$ for $\boldsymbol{x} \in \overline{D}_{\boldsymbol{x}}^{(-)}$
and 
$A(\varphi^{(+)}(\boldsymbol{x}), \boldsymbol{x}) < 0$ for $\boldsymbol{x} \in \overline{D}_{\boldsymbol{x}}^{(+)}$.

\item[(c)]
Let $h_0(\boldsymbol{x}^*,t)\in C^1(\mathbb{R}^{d-1} \times \bar{\mathcal{T}})$ be the zero-order approximation of $h(\boldsymbol{x}^*,t,\mu)$, and $h_0\in \Omega=(-a,a)$. Then,
$
\max_{\boldsymbol{x}^* \in \mathbb{R}^{d-1},\; t \in \bar{\mathcal{T}}}
\sum_{i=2}^d  \partial_{x_i} h_0 (\boldsymbol{x}^*, t)  < 1
$.
\item[(d)] 
${u}_0(\boldsymbol{x}, \mu)=U_{n-1}(\boldsymbol{x},0)+\mathcal{O}(\mu^n)$, with $n\ge 1$ being the expansion order and $U_{n-1}$ constructed below. 
\item[(e)]
For every $\boldsymbol{x}\in {\overline{D}_{\boldsymbol{x}}}$ and
$
s \in(\varphi^{(-)}(h_0,\boldsymbol{x}^*),\;\varphi^{(+)}(h_0,\boldsymbol{x}^*))$:
\[
\int_{\varphi^{(-)}(h_0,\boldsymbol{x}^*)}^{s}
\bigg(
    A(u,\boldsymbol{x}) (1- \sum_{i=2}^d \partial_{x_i}{h_0})
    - \partial_{t}{h_0}(\boldsymbol{x}^*,t)
\bigg)du>0.
\]
\end{itemize}
\end{assumption}

Assumptions \ref{assumption} ensures that there is already a formed front at $t=0$, i.e., $u_0(\boldsymbol{x}, \mu)$ has an internal transition layer around some point $h_{00}\in [-a, a]$.

\begin{remark}
When $d=1$, $h(\cdot, t,\mu)$ is a scalar function $h(t,\mu)$. Assumption \ref{assumption} still holds trivially with $\boldsymbol{x}^*=\emptyset$, so that $\varphi^{(-)}(\cdot)$ and $\varphi^{(+)}(\cdot)$ are 1D functions.
\end{remark}

\subsection{Local coordinates}
To construct the asymptotic solution near the surface $h$, we employ local coordinates $\boldsymbol{r} = (r_1, \dots, r_d)$ to facilitate the analysis. Let $\boldsymbol{r}^*= (r_2,\dots,r_d)$ if $d \geq 2$ and $ \boldsymbol{r}^*=   \emptyset$ for $d =1$.
The normal equation to the surface $h(\boldsymbol{r}^*, t,\mu)$ is given by
\begin{equation*}
\frac{x_1-h(\boldsymbol{r}^*,t,\mu)}{-1}=\frac{x_2-r_2}{\partial_{r_2} h (\boldsymbol{r}^*, t,\mu)}=\dots=\frac{x_d-r_d}{\partial_{r_d}h (\boldsymbol{r}^*, t,\mu)}. \end{equation*}
The normal distance $r_1$ from $(x_1, \dots, x_d)$ to the surface $h(\boldsymbol{r}^*, t,\mu)$ is given by:
\begin{equation*}
r_1^2 =(x_1-h(\boldsymbol{r}^*,t,\mu))^2+ (x_2-r_2)^2 + \dots + (x_d-r_d)^2  = (x_1-h(\boldsymbol{r}^*,t,\mu))^2 \Big(1+\sum_{i=2}^d (\partial_{r_i} h)^2\Big).
\end{equation*}

In the inner transition layer (i.e., in the vicinity of the surface $ h (\boldsymbol{r}^*, t,\mu)$), we use the extended variable
$ \xi=\frac{r_1}{\mu}$ and the local coordinates $ (\boldsymbol{r},t) $ using the relations
\begin{align} \label{localrelations}
x_1 = h(\boldsymbol{r}^*, t,\mu) + \frac{r_1}{\sqrt{1 + \sum_{i=2}^d (\partial_{r_i} h)^2}},\ \dots,\ x_j = r_j - \frac{r_1 h_{r_j}}{\sqrt{1 + \sum_{i=2}^d(\partial_{r_i} h)^2}},
\end{align}
 for $j = 2, \dots, d$.
We assume that $ r_1> 0 $ in the domain $  \overline{D}_{\boldsymbol{x}}^{(+)}$, $r_1 <0 $ in the domain $ \overline{D}_{\boldsymbol{x}}^{(-)} $, and note that if $ x_1 = h (\boldsymbol{r}^*, t,\mu)$, then $ r_1 = 0 $, $r_2=x_2, \dots, r_d=x_d$. The derivatives of the function $ h (\boldsymbol{r}^*, t,\mu) $ in \eqref{localrelations} are also taken for $ r_2=x_2, \dots, r_d=x_d$.

To obtain the relationship between the local coordinate system $(\boldsymbol{r},t)$ and $(\boldsymbol{x},t)$, we define the vectors ${\bf X} = (\boldsymbol{x}, t)^T$ and ${\bf R} = (\boldsymbol{r}, t)^T$, with $r_i = r_i(\boldsymbol{x}, t)$.
By taking the total differential of ${\bf R}$, we obtain $\mathrm{d} {\bf R} = \nabla_{\bf X} {\bf R} \ \mathrm{d} {\bf X}$. This determines the differential changes in the local coordinate system $(\boldsymbol{r},t)$ from the changes in $(\boldsymbol{x},t)$. To determine the reverse relationship, i.e., the differential changes in $(\boldsymbol{x},t)$ from the changes in $(\boldsymbol{r},t)$, we use the relation $\mathrm{d} {\bf X} = \nabla_{\bf R} {\bf X} \ \mathrm{d} {\bf R}$, i.e., $\nabla_{\bf X} {\bf R} = \left(\nabla_{\bf R} {\bf X} \right)^{-1}$. This provides a convenient way to transform gradients and derivatives between the two coordinate systems:
\begin{align} \label{ComponentsOfNabla}
\left(\begin{smallmatrix}
\displaystyle \frac{\partial r_1 }{\partial x_1} & \displaystyle \dots &\displaystyle \frac{\partial r_1 }{\partial x_d} &\displaystyle \frac{\partial r_1 }{\partial t}\\\\
\displaystyle \dots & \displaystyle \dots & \displaystyle \dots  &\displaystyle \dots \\\\
\displaystyle \frac{\partial r_d}{\partial x_1} &\displaystyle \dots &\displaystyle \frac{\partial r_d}{\partial x_d}&\displaystyle \frac{\partial r_d }{\partial t} \\\\
\displaystyle 0 & \displaystyle \dots  &\displaystyle 0&\displaystyle 1
\end{smallmatrix} \right)= \left(\begin{smallmatrix}
\displaystyle\frac{\partial x_1 }{\partial r_1} &\displaystyle \dots  &\displaystyle \frac{\partial x_1 }{\partial r_d} &\displaystyle \frac{\partial x_1 }{\partial t}\\\\
\displaystyle \dots  &\displaystyle \dots  &\displaystyle \dots  &\displaystyle \dots \\\\
\displaystyle\frac{\partial x_d}{\partial r_1} & \displaystyle \dots  &\displaystyle \frac{\partial x_d}{\partial r_d} &\displaystyle \frac{\partial x_d}{\partial t} \\\\
\displaystyle 0 & \displaystyle \dots  &\displaystyle 0&\displaystyle 1
\end{smallmatrix} \right)^{-1}.
\end{align}

To analyze the behavior of the solution near the transition layer, we rewrite the differential operators in \eqref{mainproblem} using the local coordinates $(\boldsymbol{r},t)$ and the extended variable $\xi$:
\begin{align}\label{nabla}
\nabla = \big\lbrace \partial_{ x_1},\, \dots,\, \partial_{x_d} \big\rbrace
= \bigg\lbrace \sum_{j=1}^d \partial_{x_1} r_j \partial_{ r_j},\, \dots,\, \sum_{j=1}^d\partial_{x_d} r_j \partial_{ r_j} \bigg\rbrace.
\end{align}
The differential operator $\Delta$ in $\mathbb{R}^d$ takes the form:
\begin{align} \label{kdimensionallaplasian}
\Delta = \sum_{i=1}^d \frac{\partial^2}{\partial x_i^2}
= \sum_{i=1}^d \sum_{j=1}^d \frac{\partial r_j}{\partial x_i} \frac{\partial}{\partial r_j} \left( \frac{\partial}{\partial x_i} \right),
\end{align}
with $ {\partial }/{\partial x_i}$ given in \eqref{nabla}.
Next we express $\partial / \partial t$ in terms of $r_1,r_2,...,r_d$ and $t$:
\begin{align} \label{kdimensionaloperatordt}
\frac{\partial}{\partial t}= \frac{\partial}{\partial t} + \frac{\partial r_1}{\partial t}\frac{\partial }{\partial r_1}+\frac{\partial r_2}{\partial t} \frac{\partial }{\partial r_2} +...+\frac{\partial r_d}{\partial t} \frac{\partial }{\partial r_d},
\end{align}
with $\partial r_i/ \partial t$ given in \eqref{ComponentsOfNabla}. This separates the contributions of each variable to the overall behavior of the solution and is useful for developing an asymptotic approximation:

\begin{remark}
For $d=1$, from \eqref{localrelations}, we derive in local coordinates $x_1=h(t,\mu)+r_1$. Equations \eqref{nabla}-\eqref{kdimensionaloperatordt} take the following form
\begin{align*}\label{operatorequat}
\nabla =\frac{\partial}{\partial r_1}= \frac{1}{\mu} \frac{\partial}{\partial \xi}, \quad
\Delta =\frac{\partial^2}{\partial r_1^2}=\frac{1}{\mu^2}\frac{\partial^2}{\partial \xi^2}, \quad
\frac{\partial}{\partial t} =\frac{\partial}{\partial t}-\frac{1}{\mu} \frac{\partial h(t)}{\partial t}\frac{\partial}{\partial \xi}.
\end{align*}
\end{remark}

\subsection{Asymptotic representation of solution}

The solution to problem~\eqref{mainproblem} is approximated by an $n$-th order asymptotic expansion with respect to $\mu$, denoted by $U_n(\boldsymbol{x}, t, \mu)$:
\begin{equation*} \label{asymptoticsolution}
U_n(\boldsymbol{x},t, \mu)=
\begin{cases}
U_n^{(-)} (\boldsymbol{x},t, \mu) , \quad (\boldsymbol{x}, t) \in \overline{D}_{\boldsymbol{x}}^{(-)} \times \bar{\mathcal{T}} ,\\
U_n^{(+)}(\boldsymbol{x}, t, \mu) , \quad (\boldsymbol{x}, t )\in \overline{D}_{\boldsymbol{x}}^{(+)} \times \bar{\mathcal{T}},
\end{cases}
\end{equation*}
with
\begin{equation}
\label{asymptoticapproximation2terms}
U_n^{(\mp)}=\bar{u}^{(\mp)}(\boldsymbol{x}, \mu)+Q^{(\mp)}(\xi, h (\boldsymbol{r}^*, t,\mu), \boldsymbol{r}^*, t, \mu),
\end{equation}
where $\bar{u}^{(\mp)}$ and $Q^{(\mp)}$ represent the outer layer and inner transition layer, respectively, expanded in powers of $\mu$ up to order $n$:
\begin{align}
\bar{u}^{(\mp)}(\boldsymbol{x}, \mu) &= \sum_{i=0}^n \mu^i \bar{u}_i^{(\mp)}(\boldsymbol{x}), \label{expansionregularfunctions}\\
Q^{(\mp)}(\xi, h, \boldsymbol{r}^*, t, \mu) &= \sum_{i=0}^n \mu^i Q_i^{(\mp)}(\xi, h, \boldsymbol{r}^*, t).\label{expansiontransitionfunctions}
\end{align}

Let $x_1 = h(\boldsymbol{x}^*, t, \mu)$ denote the internal surface at which the solution equals the average of the outer expansions:
\begin{equation} \label{halfsum}
 \phi(h,\boldsymbol{x}^*,  \mu):=\tfrac{1}{2}(\bar{u}^{(-)}(h,\boldsymbol{x}^*,\mu)+\bar{u}^{(+)}(h,\boldsymbol{x}^*,\mu)).
\end{equation}
The functions $U_n^{(\mp)}(\boldsymbol{x}, t, \mu)$ and their derivatives along the normal to the surface $x_1=h$ are matched on the surface $h(\boldsymbol{x}^*, t,\mu )$:
\begin{equation}
\begin{aligned}
U_n^{(-)} (h,\boldsymbol{x}^*,  t, \mu) &= U_n^{(+)}(h,\boldsymbol{x}^*, t, \mu) = \phi(h,\boldsymbol{x}^*, \mu), \\
\partial_{\mathbf n} U_n^{(-)}(h,\boldsymbol{x}^*, t, \mu)&=\partial_{\mathbf n} U_n^{(+)}(h,\boldsymbol{x}^*, t, \mu).
\end{aligned}
\label{sewingcond2}
\end{equation}

The surface $h(\boldsymbol{x}^*, t, \mu)$ is expanded in $\mu$ up to order $n$:
\begin{equation} \label{curveexpansion}
 h(\boldsymbol{x}^*, t, \mu)=h_{0}(\boldsymbol{x}^*, t)+\mu h_{1}(\boldsymbol{x}^*,t)+\ldots+\mu^{n}h_{n}(\boldsymbol{x}^*,t).
\end{equation}

\subsection{Outer functions}

To find the outer functions of the zeroth order, we substitute the expansions \eqref{expansionregularfunctions} into the stationary equation:
\begin{equation*}
\mu \Delta \bar{u} = A(\bar{u},\boldsymbol{x}) \sum_{i=1}^d\partial _{x_i}\bar{u}  + f(\boldsymbol{x}).
\end{equation*}
To obtain the governing equation for the first order of $\bar{u}_{0}^{(\mp)}(\boldsymbol{x})$, we expand the terms into powers series $\mu$ and equate the coefficients at $\mu^{0}$:
\begin{equation}\label{zeroorderregularequation1}
\begin{aligned}
&\begin{cases}
\displaystyle -A(\bar{u}_{0}^{(-)},\boldsymbol{x}) \sum_{i=1}^d \partial_{x_i} \bar{u}_{0}^{(-)}=f(\boldsymbol{x}),\\
\bar{u}_{0}^{(-)}(-a, \boldsymbol{x}^*)=L(\boldsymbol{x}^*), \quad  \bar{u}_{0}^{(-)}(\boldsymbol{x}+\boldsymbol{P}) = \bar{u}_{0}^{(-)}(\boldsymbol{x})
\end{cases}   \\
&\begin{cases}
\displaystyle -A(\bar{u}_{0}^{(+)},\boldsymbol{x})\sum_{i=1}^d \partial_{x_i}\bar{u}_{0}^{(+)} =f( \boldsymbol{x}),\\
\bar{u}_{0}^{(+)}( a, \boldsymbol{x}^*)=R(\boldsymbol{x}^*), \quad \bar{u}_{0}^{(+)}(\boldsymbol{x}+\boldsymbol{P}) = \bar{u}_{0}^{(+)}(\boldsymbol{x})
\end{cases}
\end{aligned}
\end{equation}
Under Assumption \ref{assumption}(a), the solution is given by
\begin{align*}
\label{u0regu}
\bar{u}_{0}(\boldsymbol{x})= \begin{cases}
\bar{u}_{0}^{(-)}(\boldsymbol{x})=\varphi^{(-)}(\boldsymbol{x}), \quad \boldsymbol{x}\in \overline{D}_{\boldsymbol{x}}^{(-)},\\
\bar{u}_{0}^{(+)}(\boldsymbol{x})=\varphi^{(+)}(\boldsymbol{x}), \quad \boldsymbol{x}\in  \overline{D}_{\boldsymbol{x}}^{(+)}.
\end{cases}
\end{align*}

\subsection{Inner functions of zero order}
The equation that determines the functions of the inner transition layer in zero order $Q_{0}^{(\mp)}(\xi, h_0, \boldsymbol{r}^*, t)$ is derived by substituting the series \eqref{expansionregularfunctions}, \eqref{expansiontransitionfunctions}, and \eqref{curveexpansion} into    \eqref{mainproblem} and \eqref{sewingcond2}. This process involves  expanding all terms in power series of $\mu$, equating the coefficients at $\mu^{-1}$ for \eqref{mainproblem}, and at $\mu^{0}$ in equation \eqref{sewingcond2}. Further, the outer part is subtracted, and the independent variables are changed to $\xi$, $\boldsymbol{r}^*$, $t$ using \eqref{nabla}, \eqref{kdimensionallaplasian}, and \eqref{kdimensionaloperatordt}. By the decay of the transition functions at infinity, we get
\begin{align} \label{transitionalfunczeroord_kdim}
\begin{dcases}
\displaystyle \frac{\partial^{2}Q_{0}^{(\mp)}}{\partial\xi^{2}}
+\displaystyle \frac{\partial_t{h_0}+A\left( \varphi^{(\mp)}+Q_{0}^{(\mp)}, h_0, \boldsymbol{r}^* \right) \left( \sum_{i=2}^d \partial_{r_i} h_0-1\right)}{\sqrt{1+\sum_{i=2}^d(\partial_{r_i} h_0)^{2}}}\frac{\partial Q_{0}^{(\mp)}}{\partial\xi}=0;\\
\varphi^{(\mp)}(h_0,\boldsymbol{r}^*)+Q_{0}^{(\mp)}(0, h_0,\boldsymbol{r}^*, t)=\phi_0(h_0,\boldsymbol{r}^*),\quad Q_{0}^{(\mp)}(\mp\infty, h_0,\boldsymbol{r}^*, t)=0.
\end{dcases}
\end{align}

Note that $Q_{0}^{(-)}$ is defined for  $\xi\leq 0$, whereas $Q_{0}^{(+)}$ for $\xi\geq 0$. Assumption \ref{assumption}(e) ensures both unique solvability of \eqref{transitionalfunczeroord_kdim} and strict monotonicity of $Q_{0}^{(\mp)}$. To find $h_0(\boldsymbol{r}^*, t)$, let
\begin{equation*} \label{auxilaryfunction_kdim}
\tilde{u}(\xi, h_0(\boldsymbol{r}^*, t))=
\begin{cases}
\varphi^{(-)} (h_0,\boldsymbol{r}^*)+Q_{0}^{(-)} (\xi, h_0,\boldsymbol{r}^*, t) , \quad \xi\leq 0,\\
\varphi^{(+)}(h_0,\boldsymbol{r}^*)+Q_{0}^{(+)}(\xi, h_0,\boldsymbol{r}^*, t) , \quad \xi\geq 0.
\end{cases}
\end{equation*}
Then, \eqref{transitionalfunczeroord_kdim} can be expressed as:
\begin{equation} \label{replacement_kdim}
 \frac{\partial^{2}\tilde{u}}{\partial\xi^{2}}+ \frac{\partial_t{h_0}+A(\tilde{u},h_0, \boldsymbol{r}^*)\left(\sum_{i=2}^d\partial_{r_i} h_0-1\right)}{\sqrt{1+\sum_{i=2}^d(\partial_{r_i} h_0)^{2}}}\frac{\partial\tilde{u}}{\partial\xi}=0.
\end{equation}
Let $\frac{\partial \tilde{u} }{\partial\xi} = g(\tilde{u})$ and $\frac{\partial^2 \tilde{u} }{\partial\xi^{2}} = \frac{\partial g(\tilde{u}) }{\partial \tilde{u} } g(\tilde{u})$. Equation \eqref{replacement_kdim} is transformed into
{\small
\begin{align} \label{derivativetildeu_kdim}
\frac{\partial\tilde{u}}{\partial\xi}= \begin{cases}
\displaystyle \Phi^{(-)} (\xi, h_0(\boldsymbol{r}^*, t))=\int_{\varphi^{(-)}}^{\tilde{u}} \displaystyle \frac{-\partial_{t}{h_0}+A( {u},h_0, \boldsymbol{r}^*)\left(1-\sum_{i=2}^d\partial_{r_i} h_0 \right)}{\sqrt{1+\sum_{i=2}^d(\partial_{r_i} h_0)^{2}}} du ,~ \xi\leq 0,\\
\displaystyle \Phi^{(+)}(\xi, h_0(\boldsymbol{r}^*, t))= \int_{\varphi^{(+)}}^{\tilde{u}} \displaystyle \frac{-\partial_t{h_0}+A( {u},h_0, \boldsymbol{r}^*)\left(1-\sum_{i=2}^d\partial_{r_i} h_0 \right)}{\sqrt{1+\sum_{i=2}^d(\partial_{r_i} h_0)^{2}}} du ,~ \xi\geq 0.
\end{cases}
\end{align}
}
By using expansions \eqref{asymptoticapproximation2terms}--\eqref{sewingcond2}, we equate the $\mu^{-1}$ terms in \eqref{sewingcond2} and get
\begin{align}\label{sewindcondexpanded0_kdim}
\Phi^{(-)} (0, h_0(\boldsymbol{x}^*, t))-\Phi^{(+)} (0, h_0(\boldsymbol{x}^*, t))=0.
\end{align}
From \eqref{derivativetildeu_kdim} and \eqref{sewindcondexpanded0_kdim}, we deduce
\begin{equation} \label{integral_eq_kdim}
\int_{\varphi^{(-)}(h_0,\boldsymbol{x}^*)}^{\varphi^{(+)}(h_0,\boldsymbol{x}^*) } \frac{-\partial_t{h_0}-A({u},h_0, \boldsymbol{x}^*) (\sum_{i=2}^d\partial_{x_i} h_0-1)}{\sqrt{1+\sum_{i=2}^d(\partial_{x_i} h_0)^{2}}} du = 0.
\end{equation}
This is an integral equation for determining $h_0(\boldsymbol{x}^*,t)$ with the initial condition $ h_0(\boldsymbol{x}^*,0)=h_0^*$.
Under Assumptions \ref{assumption}(b)--(c), \eqref{integral_eq_kdim} has a solution inside the region $\Omega$ for
any $(\boldsymbol{x},t) \in \mathbb{R}^d\times \bar{\mathcal{T}}$. The functions $Q_{0}^{(\mp)}(\xi, h_0, \boldsymbol{x}^*, t)$
can be obtained by solving \eqref{derivativetildeu_kdim}.

Similar to the analysis in \cite{chaikovskii2022convergence}, one can show that the transition layer functions $Q_{0}^{(\mp)}(\displaystyle\xi, h_0, \boldsymbol{x}^*, t)$ decay exponentially as $\xi \rightarrow \mp \infty$ and satisfy the following estimates:
\begin{align}
\label{equatrk22}
C_{\min} e^{ \xi\kappa_{\min}} &\leq |Q_{0}^{(-)}(\xi, h_0, \boldsymbol{x}^*, t)| \leq  C_{\max} e^{ \xi \kappa_{\max}},
\quad \xi \leq 0, \quad t \in \bar{\mathcal{T}}, \\
\label{equatrk23}
C_{\min} e^{-\xi\kappa_{\min} } &\leq |Q_{0}^{(+)}(\xi, h_0, \boldsymbol{x}^*, t)| \leq  C_{\max} e^{- \xi \kappa_{\max}},
\quad \xi \geq 0, \quad t \in \bar{\mathcal{T}}.
\end{align}
where the constants $C_{\min}, C_{\max}$ and $\kappa_{\min} , \kappa_{\max}$ are independent of $\xi$, $\boldsymbol{x}$ and $ t$.

\begin{remark}\label{remark1}
For $A({u},h_0, \boldsymbol{x}^*)= -u$, equation \eqref{integral_eq_kdim} can be rewritten as:
    \begin{equation*} \label{h0rkequation}
\partial_{t}{h_0}=\frac{1}{2} \left(\sum_{i=2}^d \partial_{x_i} h_0-1\right) \left(\varphi^{(+)}(h_0,\boldsymbol{x}^* )+\varphi^{(-)}(h_0,\boldsymbol{x}^*) \right),
\end{equation*}
with an additional initial condition ${h_0}(\boldsymbol{x}^*, 0)=h_{0}^{*} \in \bar{\Omega}$ and periodic conditions along the axes. Moreover, the functions $Q_{0}^{(\mp)}(\xi, h_0, \boldsymbol{x}^*, t)$ take the form:
\begin{align*} \label{Q0rkequation}
\displaystyle Q_{0}^{(\mp)}(\xi, h_0, \boldsymbol{x}^*, t)=\frac{2P^{(\mp)}(h_0,\boldsymbol{x}^{*})}{\exp \left(-\xi \frac{P^{(\mp)}(h_0,\boldsymbol{x}^{*}) \left(1-\sum_{i=2}^d\partial_{x_i} h_0 \right)}{\sqrt{1+\sum_{i=2}^d(\partial_{x_i} h_0)^2}} \right)+1},
\end{align*}
where $
P^{(\pm)}(h_0, \boldsymbol{x}^{*}) = \frac{1}{2} \left( \varphi^{(\mp)}(h_0, \boldsymbol{x}^{*}) - \varphi^{(\pm)}(h_0, \boldsymbol{x}^{*}) \right).
$
\end{remark}

Similar to \cite{chaikovskii2022convergence}, from the estimates \eqref{equatrk22} and \eqref{equatrk23}, the width of the inner transition layer $ \Delta h$ can be estimated as $\Delta h \sim \mu |\ln \mu|$.

Now we can present the main result on the approximation properties of the asymptotic expansion. We only consider the zeroth-order approximation, as higher-order approximations are much more expensive yet yield only minor improvements.

\begin{theorem} \label{MainThm}
Suppose that  $f(\boldsymbol{x})$ and ${u}_0(\boldsymbol{x} , \mu)$ belong to $  C^{3}\bigl(\overline{D}_{\boldsymbol x}\bigr)$, that $ L(\boldsymbol{x}^*)$ and $ R(\boldsymbol{x}^*) $ belong to $  C^{3}\bigl(\mathbb R^{d-1}\bigr)$, that all four functions are $\boldsymbol{x}^*$-periodic, and fix $0 < \mu \ll 1$. Then, under Assumption \ref{assumption}, problem \eqref{mainproblem} has a unique smooth solution with an inner transition layer. The zeroth-order asymptotic solution $U_{0}(\boldsymbol{x},t,\mu)$ is given by 
\begin{align} \label{eq001}
U_{0}(\boldsymbol{x},t)=\begin{cases}
\varphi^{(-)}(\boldsymbol{x})+Q_{0}^{(-)}(\xi_0, h_0, \boldsymbol{x}^*, t) , \ (\boldsymbol{x},t)\in \bar{\Omega}^{(-)} \times\mathbb{R}^{d-1} \times \bar{\mathcal{T}},\\
\varphi^{(+)}(\boldsymbol{x})+Q_{0}^{(+)}(\xi_0, h_0, \boldsymbol{x}^*, t) , \ (\boldsymbol{x},t)\in \bar{\Omega}^{(+)} \times\mathbb{R}^{d-1} \times \bar{\mathcal{T}},
\end{cases}
\end{align}
with $\xi_0=\left(\boldsymbol{x}-h_0 \right) \sqrt{1+\sum_{i=2}^{d} (\partial_{x_i} h_0)^{2}}/ \mu$, $\bar{\Omega}^{(-)}=[-a,h_{0}(\boldsymbol{x}^*,t)], \bar{\Omega}^{(+)}=[h_{0}(\boldsymbol{x}^*,t),a]$. Moreover, the following estimates hold:
\begin{align*}
\forall (\boldsymbol{x},t)\in \bar{\Omega} \times\mathbb{R}^{d-1} \times  \bar{\mathcal{T}}:~ &|u(\boldsymbol{x},t)-U_{0}(\boldsymbol{x},t)|=\mathcal{O}(\mu), \\
\forall t \in \bar{\mathcal{T}}:~& |h(\boldsymbol{x}^*,t)-h_{0}(\boldsymbol{x}^*,t)|=\mathcal{O}(\mu |\ln \mu|) . 
\end{align*}
Furthermore, outside the narrow region $\left(h_{0}\left(\boldsymbol{x}^*,t\right)- \Delta h /2, h_{0}(\boldsymbol{x}^*,t)+\Delta h /2\right)$ with $\Delta h\sim\mu|\ln\mu|$, there exists a constant $C$ independent of $\boldsymbol{x}, t$ and $\mu$ such that
\begin{align}
|u(\boldsymbol{x},t)-\varphi^{(-)}(\boldsymbol{x})| &\leq C \mu,  &&(\boldsymbol{x},t) \in [-a, h_{0}(\boldsymbol{x}^*,t)- \Delta h /2] \ \times\mathbb{R}^{d-1} \times  \bar{\mathcal{T}},  \\
|u(\boldsymbol{x},t)-\varphi^{(+)}(\boldsymbol{x})| &\leq C \mu,  &&(\boldsymbol{x},t) \in [h_{0}(\boldsymbol{x}^*,t)+\Delta h /2, a] \times\mathbb{R}^{d-1} \times  \bar{\mathcal{T}}. 
\end{align}
\end{theorem}

The proof of Theorem \ref{MainThm} is given in the supplementary material.
\section{Deep learning for the problem}\label{sec3}
In this section, we develop the main algorithm. First, we describe the PINN  for solving problem (\ref{mainproblem}). Then, we propose the DAE method and present an adaptive sampling technique to enhance the fitting accuracy of DAE. Finally, we present an approximation theory of the DAE. Throughout this section, we abbreviate the trained neural network solution $\hat{h}_{0}(\boldsymbol{x}^*,t;\tilde{\theta}^*)$ to $\hat{h}_0$.

\subsection{PINN}\label{sec31}
In the standard PINN, we represent the unknown solution $u(\boldsymbol{x},t)$ of problem \eqref{mainproblem} using a deep neural network (DNN) $\hat{u}_\theta(\boldsymbol{x},t)\equiv \hat{u}(\boldsymbol{x},t;\theta)$, where $\theta$ denotes all trainable weights $W$ and biases $b$ in the neural network.

By substituting $\hat{u}_\theta(\boldsymbol{x},t)$ into equation \eqref{mainproblem}, we define the PDE residual:
\begin{align*}
    \displaystyle \mathcal{G}_\theta(\boldsymbol{x},t):= \mu \Delta \hat{u}_\theta(\boldsymbol{x}, t) - \partial_t \hat{u}_\theta(\boldsymbol{x}, t) - A(\hat{u}_\theta(\boldsymbol{x}, t), \boldsymbol{x}) ( \boldsymbol{1}_k  \cdot \nabla \hat{u}_\theta(\boldsymbol{x}, t)) - f(\boldsymbol{x}).
\end{align*}
All partial derivatives $\partial_t$, $\Delta$ and $\nabla$ can be computed using automatic differentiation (AD). The empirical loss for problem (\ref{mainproblem}) is given by
\begin{align*}
L(\theta;N) = L_{\text{PDE}}(\theta;N_f) +  L_{\text{BC}}(\theta;N_b) + L_{\text{IC}}(\theta;N_i),
\end{align*}
with 
\begin{align*}
    &L_{\text{PDE}}(\theta;N_f) = \frac{1}{N_{f}} \sum_{j=1}^{N_{f}} \big| \mathcal{G}_\theta(\boldsymbol{x}_f^j,t_f^j) \big|^2,~ L_{\text{IC}}(\theta;N_i) = \frac{1}{N_{i}} \sum_{j=1}^{N_{i}} \big| \hat{u}_\theta(\boldsymbol{x}_i^j, 0) - {u}_0(\boldsymbol{x}_i^j, \mu) \big|^2,\\
&L_{\text{BC}}(\theta;N_b) = \frac{1}{N_{b}} \sum_{j=1}^{N_{b}} \Big( \big| \hat{u}_\theta(-a, (\boldsymbol{x}^*)_{b_1}^j, t_b^j) - L((\boldsymbol{x}^*)_{b_1}^j) \big|^2  \\
& \qquad\qquad\qquad + \big| \hat{u}_\theta(a, (\boldsymbol{x}^*)_{b_1}^j, t_b^j) - R((\boldsymbol{x}^*)_{b_1}^j) \big|^2+| \hat{u}_\theta(\boldsymbol{x}_{b_2}^j + \boldsymbol{P}, t_b^j) - \hat{u}_\theta(\boldsymbol{x}_{b_2}^j, t_b^j) |^2\Big).
\end{align*}
The first term \(L_{\text{PDE}}\) measures the residuals of the PDE evaluated at the interior points \(\{\boldsymbol{x}_f^j, t_f^j\}_{j=1}^{N_f}\) randomly sampled in the domain \(D_{\boldsymbol{x}} \times (0,T]\). The second term \(L_{\text{BC}}\) enforces the Dirchlet boundary condition at \(\{(\boldsymbol{x}^*)_{b_1}^j, t_b^j\}_{j=1}^{N_b}\) and the periodic condition at \(\{ \boldsymbol{x}_{b_2}^j, t_b^j\}_{j=1}^{N_b}\). Similarly, the third term \(L_{\text{IC}}\) enforces the initial condition at  \(\{\boldsymbol{x}_i^j\}_{j=1}^{N_i}\). {The baseline PINN provides a general solution strategy. However, for SPPs, the approach struggles to resolve moving transition layers with sharp gradients, and the residual errors concentration near the layer, leading to degraded convergence.}

\subsection{Deep asymptotic expansion (DAE)}
The analysis in Section~\ref{sec2} provides the inner and outer solutions. Nonetheless,  equation \eqref{integral_eq_kdim} governing the transition-layer position often requires mesh-based discretization. 
Thus there is a trade-off  between PINN and asymptotic analysis: the former is meshless but struggles to resolve the transition layer accurately, whereas the latter exposes the structure yet its implementation relies on meshes. This observation motivates the present work.

\begin{figure}[htbp]
    \centering
    \includegraphics[width=\linewidth]{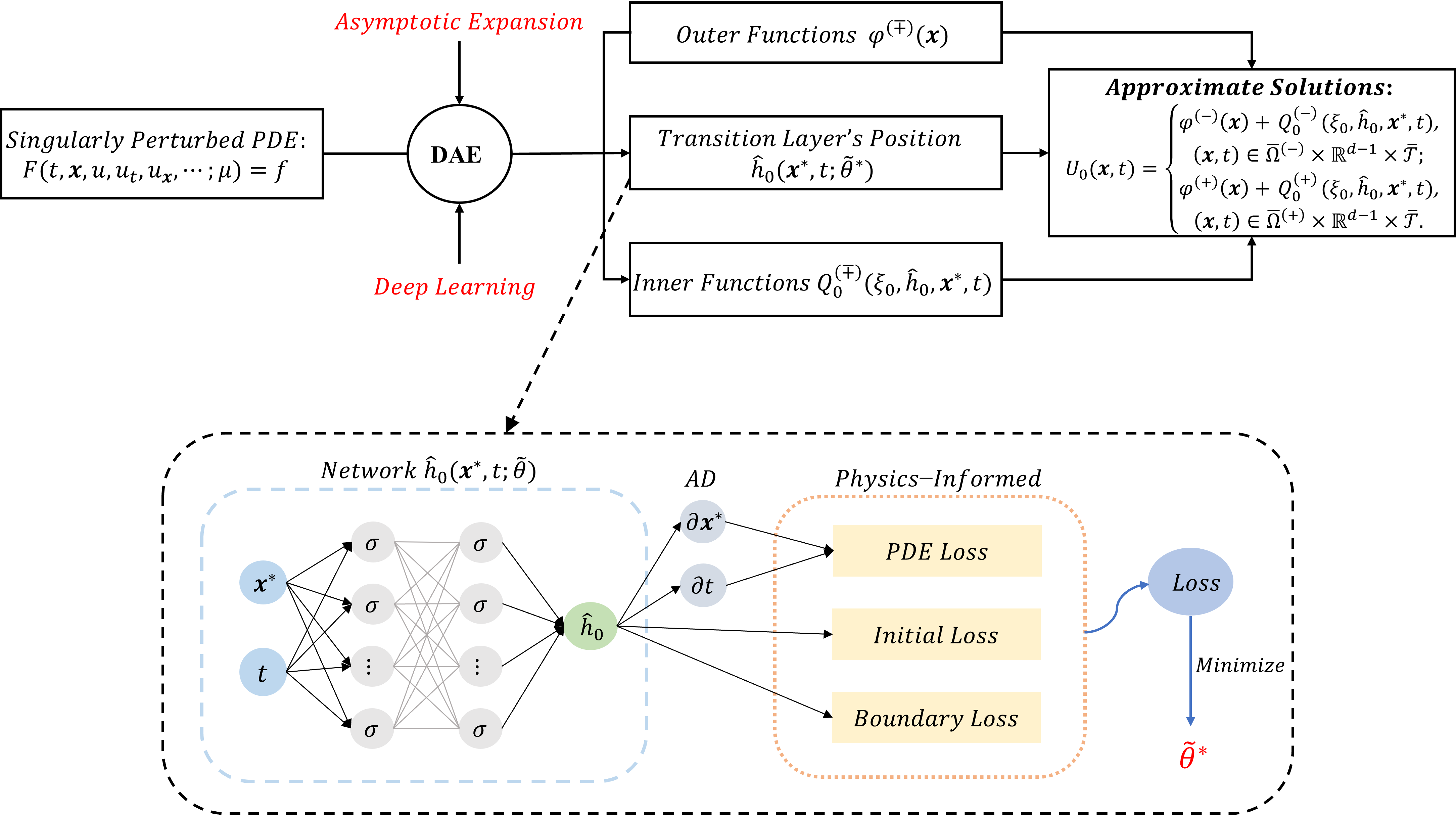}
    \caption{DAE algorithm framework.}
    \label{fig:DAE algorithm framework}
\end{figure}

To address these challenges, we propose the DAE framework in Fig.~\ref{fig:DAE algorithm framework}. The core idea of DAE is to combine the physical insights of asymptotic analysis with the universal approximation power of neural networks. It yields a mesh-free solution strategy. Instead of applying a neural network to the full problem \eqref{mainproblem} or discretizing the asymptotic reduction, DAE uses the network as a differentiable function approximator for the transition-layer position governed by Eq.~\eqref{integral_eq_kdim}. It not only avoids mesh dependence but also, by targeting a lower-dimensional and simpler subproblem, greatly reduces computational burden and training difficulty. Furthermore, since the network learns only the more tractable subproblem, the training is easier and it exhibits enhanced generalization across different parameters or initial conditions.

Specifically, we approximate the transition layer position $h_{0}(\boldsymbol{x}^*,t)$ by a neural network $\hat{h}_{0}(\boldsymbol{x}^*,t;\tilde{\theta})$, which is constructed as:
\begin{equation}\label{h_net}
\hat{h}_{0}(\boldsymbol{x}^*,t;\tilde{\theta}) = h_0^* + t \cdot {h}(\boldsymbol{x}^*,t;\tilde{\theta}),
\end{equation}
where $h(\cdot)$ is a neural network with parameters $\tilde{\theta}$. This choice enforces the initial condition $h_0(\boldsymbol{x}^*,0)=h_0^*$ exactly, which helps improve the training stability and overall solution accuracy.

The network parameters $\tilde{\theta}$ are optimized by minimizing a composite loss function $L_{\hat{h}_{0}}(\tilde{\theta})$ that consists of a PDE residual term and a periodicity constraint. First, we define the PDE residual at each collocation point $(\boldsymbol{x}^*_j, t_j)$, denoted as $\mathcal{R}_j$. It quantifies the error in satisfying the governing equation~\eqref{integral_eq_kdim} and is computed by discretizing the integral using an $M$-point Gauss-Legendre quadrature scheme. With the shorthand $\hat{h}_{0,j} = \hat{h}_{0}(\boldsymbol{x}^*_j,t_j;\tilde{\theta})$, the residual $\mathcal{R}_j$ is given by:
\begin{equation}
\mathcal{R}_j \equiv \sum_{m=1}^{M} w'_{m,j} \frac{-\partial_t\hat{h}_{0,j}-A({u_{m,j}},\hat{h}_{0,j}, \boldsymbol{x}^*_j) \left(\sum_{i=2}^d\partial_{x_i} \hat{h}_{0,j}-1\right)}{\sqrt{1+\sum_{i=2}^d(\partial_{x_i} \hat{h}_{0,j})^{2}}},
\end{equation}
where the mapped Gaussian nodes $u_{m,j}$ and weights $w'_{m,j}$ are dynamically computed for each point $j$ over the integration interval $[\varphi^{(-)}(\hat{h}_{0,j}, \boldsymbol{x}^*_j), \varphi^{(+)}(\hat{h}_{0,j}, \boldsymbol{x}^*_j)]$.

The empirical loss function consists of the mean squared residuals over $N_h$ points and the mean squared periodicity error over $N_{hp}$ points:
\begin{equation}
\label{eq:h0_loss}
L_{\hat{h}_{0}}(\tilde{\theta}) = \frac{1}{N_{h}} \sum_{j=1}^{N_{h}} \left| \mathcal{R}_j \right|^2 + \frac{1}{N_{hp}} \sum_{j=1}^{N_{hp}} \Big| \hat{h}_{0}(\boldsymbol{x}^*_j+\boldsymbol{P},t_j;\tilde{\theta}) - \hat{h}_{0,j} \Big|^2.
\end{equation}
Note that embedding the initial condition directly into the network architecture is effective for the DAE method, whereas for the standard PINN applied to the full system~\eqref{mainproblem} (cf. Section~\ref{sec31}), this approach tends to produce excessively large initial losses, which makes the training difficult. The complete procedure of the DAE is summarized in Algorithm~\ref{algorithmDAE}.

The success of the DAE in delivering high-fidelity numerical results depends crucially on the accurate resolution of the neural network approximation for the position of transition layer, $\hat{h}_{0}(\boldsymbol{x}^*,t;\tilde{\theta})$. The associated errors can severely compromise the accuracy of the overall solution $U_{0}(\boldsymbol{x},t;\tilde{\theta}^*)$, especially for small $\mu$. While there are various strategies to improve the network's predictive accuracy, e.g., modifying its architecture or tuning loss function weights, we focus on a particularly effective approach: optimizing the training data distribution via adaptive sampling.

We employ the RAR strategy proposed by Lu et al. \cite{lu2021deepxde}. The core idea is to use the magnitude of the point-wise PDE residual $|\mathcal{R}_j|$ as an error indicator, based on the observation that the largest residuals often emerge where the network's approximation is least accurate. This is particularly pronounced in regions where the interface geometry exhibits high curvature or undergoes rapid temporal changes. The RAR algorithm iteratively identifies points with the highest residuals and adding them to the training set. This process creates an adaptive sampling scheme that dynamically focuses the model's capacity on resolving these challenging regions of the solution domain, facilitating the learning of a more accurate and well-defined solution. The procedure for the adaptive refinement is detailed in Algorithm~\ref{algorithmRAR}.

\begin{remark}
The DAE is suitable when Eq.~\eqref{zeroorderregularequation1}--\eqref{transitionalfunczeroord_kdim} admits closed-form solutions so that $\varphi^{(\mp)}$ and $Q_{0}^{(\mp)}$ are explicit. This is not limited to $A(u,\boldsymbol{x})=-u$; any $A(u,\boldsymbol{x})$ that allows Eq.~\eqref{zeroorderregularequation1}--\eqref{transitionalfunczeroord_kdim} to have an explicit solution is also admissible. When explicit formulas are unavailable, $\varphi^{(\mp)}$ and $Q_{0}^{(\mp)}$ can be approximated by neural networks.
\end{remark}

\begin{algorithm}[htbp]
\caption{Deep Asymptotic Expansion (DAE)}
\label{algorithmDAE}
\begin{algorithmic}[1]
\STATE{\textbf{Input:} Initial condition $h_0^*$; training set sizes $N_h, N_{hp}$; number of iterations $K$.}
\STATE{\textbf{Output:} The zero-order asymptotic approximate solution $U_{0}(\boldsymbol{x},t;\tilde{\theta}^*)$.}

\STATE{Compute $\varphi^{(\mp)}(\boldsymbol{x})$ from Eq.~\eqref{zeroorderregularequation1}.}
\STATE{Construct the network $\hat{h}_{0}(\boldsymbol{x}^*,t;\tilde{\theta})$ as in Eq.~\eqref{h_net} and initialize it with Xavier weights and zero biases.}
\STATE{Generate training sets $\mathcal{N}_1 = \{\boldsymbol{x}_j^*, t_j\}_{j=1}^{N_{h}}$ (for the PDE residual) and $\mathcal{N}_2 = \{\boldsymbol{x}_j^*, t_j\}_{j=1}^{N_{hp}}$ (for the periodicity).}

\FOR{$k=1$ \textbf{to} $K$}
    \STATE{Compute the loss $L_{\hat{h}_{0}}(\tilde{\theta})$ on $\mathcal{N}_1$ and $\mathcal{N}_2$ according to (\ref{eq:h0_loss}).}
    \STATE{Update parameters $\tilde{\theta}$ using an Adam optimizer.}
\ENDFOR
\STATE{Obtain the trained parameters $\tilde{\theta}^* \leftarrow \tilde{\theta}$.}

\STATE{Compute $Q_{0}^{(\mp)}(\xi_0,\hat{h}_{0},\boldsymbol{x}^*,t)$ from Eq.~\eqref{transitionalfunczeroord_kdim} with $h_{0}=\hat{h}_{0}(\boldsymbol{x}^*,t;\tilde{\theta}^*)$.}
\STATE{Assemble $U_{0}(\boldsymbol{x},t;\tilde{\theta}^*)$ using (\ref{eq001}) with $\varphi^{(\mp)}(\boldsymbol{x})$ and $Q_{0}^{(\mp)}(\xi_0, \hat{h}_{0}, \boldsymbol{x}^*, t)$.}

\STATE{\textbf{Return:} $U_{0}(\boldsymbol{x},t;\tilde{\theta}^*)$.}
\end{algorithmic}
\end{algorithm}

\begin{algorithm}[htbp]
\caption{Residual-based Adaptive Refinement (RAR) for DAE}
\label{algorithmRAR}
\begin{algorithmic}[1]
\STATE \textbf{Input:} 
    Initial network $\hat{h}_{0}(\cdot;\tilde{\theta})$ (pre-trained on $\mathcal{T}_{\text{train}}$); 
    Candidate point set $\mathcal{S}$; 
    Refinement batch size $m$; 
    Residual tolerance $\mathcal{E}_0$.

\STATE \textbf{Output:} The refined network $\hat{h}_{0}(\cdot;\tilde{\theta}^*)$.

\STATE Compute residuals for all points in $\mathcal{S}$ and calculate the mean:
       $\mathcal{E}_r \leftarrow \frac{1}{|\mathcal{S}|} \sum_{(\boldsymbol{x}^*_j, t_j) \in \mathcal{S}} |\mathcal{R}_j|$

\WHILE{$\mathcal{E}_r \ge \mathcal{E}_0$}
    \STATE Select subset $\mathcal{S}_{\text{add}} \subset \mathcal{S}$ of $m$ points with the largest residual magnitudes $|\mathcal{R}_j|$.
    \STATE Augment training set: $\mathcal{T}_{\text{train}} \leftarrow \mathcal{T}_{\text{train}} \cup \mathcal{S}_{\text{add}}$.
    \STATE Retrain the network on the updated $\mathcal{T}_{\text{train}}$.
    \STATE Re-compute the mean residual $\mathcal{E}_r$ on $\mathcal{S}$ using the updated network.
\ENDWHILE
\end{algorithmic}
\end{algorithm}

\subsection{Approximation theory for DAE}
\label{sec:DAE_theory}
We develop an approximation theory for the DAE framework to quantify the gap between the approximate solution $U_0(\boldsymbol{x}, t;\tilde{\theta}^*)$ and the exact solution $u(\boldsymbol{x}, t)$ of \eqref{mainproblem}.

It is classical that sufficiently large neural networks can uniformly approximate functions and their derivatives on compact sets \cite{hornik1991approximation}.  However, practical networks are size-limited. Let $\mathcal{F}$ denote the set of functions representable by the chosen neural network architecture. Since the solution $h_0$ is unlikely to lie in $\mathcal{F}$, consider the best approximation in 
$\mathcal{F}$ defined by
$h_{\mathcal{F}} = \arg \min_{\widetilde{h}_0 \in \mathcal{F}} \| \widetilde{h}_0 - h_0 \|_{C^1(B\times \bar{\mathcal{T}})}$, where $B \subset \mathbb{R}^{d-1}$ is a bounded closed set.
Since the training uses only the dataset $\mathcal{N}=\{\mathcal{N}_1,\mathcal{N}_2\}$, we define $h_{\mathcal{N}}$ as the network that globally minimizes the empirical loss $L_{\hat{h}_{0}}(\tilde{\theta})$ on $\mathcal{N}$. We may assume that $h_0$, $h_{\mathcal{F}}$, and $h_{\mathcal{N}}$ are well-defined and unique.
Finding the exact minimizer $h_{\mathcal{N}}$ is generally intractable \cite{blum1988training}, the optimizer typically returns an approximate solution $\hat{h}_0$. Therefore, the total error can be decomposed into \cite{bottou2007tradeoffs}:
\begin{align}\label{total_error_h}
\| h_0 - \hat{h}_0 \|&_{C^1(B\times \bar{\mathcal{T}})}\\
\leq& \underbrace{\| h_0 - h_{\mathcal{F}} \|_{C^1(B\times \bar{\mathcal{T}})}}_{\text{Approximation Error } \varepsilon_{\mathrm{app}}}
+
\underbrace{\| h_{\mathcal{F}} - h_{\mathcal{N}} \|_{C^1(B\times \bar{\mathcal{T}})}}_{\text{Generalization Error } \varepsilon_{\mathrm{gen}}}+\underbrace{\| h_{\mathcal{N}} - \hat{h}_0 \|_{C^1(B\times \bar{\mathcal{T}})}}_{\text{Optimization Error } \varepsilon_{\mathrm{opt}}}\notag.
\end{align}
{The approximation error $\varepsilon_{\mathrm{app}}$ arise from the limited expressivity of the network class $\mathcal{F}$ and decreases as the network size grows.
The generalization error $\varepsilon_{\mathrm{gen}}$ reflects the discrepancy between the best-in-class solution and the one obtained from finite training data; it can be reduced by increasing the number of training samples.
Finally, the optimization error $\varepsilon_{\mathrm{opt}}$ is caused by the difficulty of finding the global minimizer of the non-convex loss; in practice, it depends on the optimizer, training iterations, and hyperparameter choices.}

From the error decomposition \eqref{total_error_h}, we can derive an error bound on the DAE. We denote function $Q_{0}^{(\mp)}(\xi_0, h_0, \boldsymbol{x}^*, t)$ by $\widetilde{Q}_{0}^{(\mp)}(h_0, \partial_{x_2}h_0, \ldots, \partial_{x_d}h_0, \boldsymbol{x}^*, t)$, since $\xi_0$ depends on $h_0$ and its spatial derivatives.
\begin{theorem}\label{DAE_error_thm}
Let the assumptions of Theorem~\ref{MainThm} hold, and  there exists an upper bound ${M} > 0$ such that the total error satisfies 
$\varepsilon_{\mathrm{app}} + \varepsilon_{\mathrm{gen}} + \varepsilon_{\mathrm{opt}} \leq {M}$. If $\widetilde{Q}_{0}^{(\mp)}(h_0, \partial_{x_2}h_0, \ldots, \partial_{x_d}h_0, \boldsymbol{x}^*, t)$ is continuously differentiable with respect to each component, then there exists $\hat{h}_{0}(\boldsymbol{x}^*, t; \tilde{\theta}^*) \in \mathfrak{S}$ such that, for all $(\boldsymbol{x}, t) \in \bar{\Omega} \times B \times \bar{\mathcal{T}}$, the following error estimate holds:\begin{equation} \label{DAE_error}
    | u(\boldsymbol{x}, t) - U_0(\boldsymbol{x}, t;\tilde{\theta}^*) | 
    \leq C_1 \mu + C_2(\varepsilon_{\mathrm{app}}+\varepsilon_{\mathrm{gen}}+\varepsilon_{\mathrm{opt}}),
\end{equation}
where $C_1$ is a constant, and $C_2$ is a constant depending on $M$ and the domain $B \times \bar{\mathcal{T}}$.
\end{theorem}
\begin{proof}
By the triangle inequality,
\begin{align*}
| u(\boldsymbol{x}, t) - U_0(\boldsymbol{x}, t; \tilde{\theta}^*) | 
&\leq | u(\boldsymbol{x}, t) - U_0(\boldsymbol{x}, t) | 
+ | U_0(\boldsymbol{x}, t) - U_0(\boldsymbol{x}, t; \tilde{\theta}^*) |,
\end{align*}
where the term $| u(\boldsymbol{x}, t) - U_0(\boldsymbol{x}, t) |$ is $\mathcal{O}(\mu)$ by Theorem \ref{MainThm}. Since $\widetilde{Q}_{0}^{(\mp)}$ is continuously differentiable with respect to each component, by the mean value theorem, we have
\begin{equation}\label{Q1_est}
\begin{aligned}
&|Q_{0}^{(\mp)}(\xi_0, h_0, \boldsymbol{x}^*, t)-Q_{0}^{(\mp)}(\xi_0, \hat{h}_0, \boldsymbol{x}^*, t)| \\
& =|\widetilde{Q}_{0}^{(\mp)}(h_0, \partial_{x_2}h_0,\cdots, \partial_{x_d}h_0,\boldsymbol{x}^*, t) - \widetilde{Q}_{0}^{(\mp)}(\hat{h}_0, \partial_{x_2}\hat{h}_0,\cdots, \partial_{x_d}\hat{h}_0,\boldsymbol{x}^*, t)| \\
& \le \sum_{i=1}^d \max_{a_1,\ldots,a_d \leq  2\| h_0 \|_{C^1(B\times \bar{\mathcal{T}})}+ M} |\partial_{a_i} \widetilde{Q}_{0}^{(\mp)}(a_1, a_2,\ldots, a_d,\boldsymbol{x}^*, t)| \cdot \| {h}_0-\hat{h}_0  \|_{C^1(B\times \bar{\mathcal{T}})}\\
& := C_2\| h_0 - \hat{h}_0 \|_{C^1(B\times \bar{\mathcal{T}})}.
\end{aligned}
\end{equation}Then the estimate of $| U_0(\boldsymbol{x}, t) - U_0(\boldsymbol{x}, t; \tilde{\theta}^*)|$ follows by combining (\ref{eq001}), (\ref{Q1_est}) and (\ref{total_error_h}).
The detailed derivation of (\ref{Q1_est}) is given in the supplementary material.
\end{proof}

\begin{remark}
    For $A({u},h_0, \boldsymbol{x}^*)= -u$, the assumption on $Q_{0}^{(\mp)}$ is satisfied. The function $\widetilde{Q}_{0}^{(\mp)}(h_0, \partial_{x_2}h_0,\ldots, \partial_{x_d}h_0,\boldsymbol{x}^*, t)$ is continuously differentiable with respect to $h_0$, the spatial derivatives $\partial_{x_2}h_0,\ldots, \partial_{x_d}h_0$, $\boldsymbol{x}^*$, and $t$; see Remark \ref{remark1}. 
\end{remark}
\section{Experimental results}
\label{sec4}
In this section, we employ the DAE, PINN and gPINN to solve problem \eqref{mainproblem} to assess the performance of DAE. Throughout, the hyperparameters for each experiment are listed in Table \ref{tab:exp_parameters}, unless otherwise noted. Since the DAE uses deep learning techniques after reducing the problem dimensionality, it employs fewer training points and iterations than PINN. All networks use tanh activation function. All computations were performed on an NVIDIA L40S GPU, within the PyTorch framework. In all the experiments, we set the seed to 1234 for reproducibility, except for the sensitivity experiments on random seeds. The data and associated codes are available at GitHub, \url{https://github.com/z1998w/DAE}. 

{We employ the problems and the asymptotic solutions from existing works \cite{chaikovskii2022convergence,chaikovskii2023asymptotic,chaikovskii2023solving} for the evaluation. This choice is motivated by the computational infeasibility of conventional mesh-based methods (e.g., finite volume method). Indeed,  even for the 1D problem with a small $\mu_1$, resolving the sharp layer demands a mesh so fine that the runtime is prohibitive.} The accuracy is evaluated on a test set using the final loss \( e_{{loss}} \), relative \( L^2 \) error \( e_2 \), and maximum absolute error \( e_{\infty} \) defined by
\begin{align*}
e_2 =
\frac{\sqrt{\sum_{i=1}^{N_{\text{test}}} |\hat{u}(\boldsymbol{x}_i, t_i; \boldsymbol{\theta}^*) - u(\boldsymbol{x}_i, t_i)|^2}}
     {\sqrt{\sum_{i=1}^{N_{\text{test}}} |u(\boldsymbol{x}_i, t_i)|^2}}, ~ e_{\infty} =
\max_{i=1, \dots, N_{\text{test}}} |\hat{u}(\boldsymbol{x}_i, t_i; \boldsymbol{\theta}^*) - u(\boldsymbol{x}_i, t_i)|,
\end{align*}
where $u(\boldsymbol{x}_i,t_i)$ represents the reference solution mentioned above. The set $\{(\boldsymbol{x}_i, t_i)\}_{i=1}^{N_{\text{test}}}$ is generated via Latin Hypercube Sampling \cite{mckay2000comparison}, which provides an efficient means of estimating error norms, particularly for higher-dimensional problems. The number of test points,  $N_{\text{test}}$, was set to 5000, 10000, and 13000 for the 1D, 2D, and 3D problems, respectively. All competing methods were evaluated on this identical set of points, establishing an unbiased basis for performance comparison.

\begin{table}[htbp]
\centering
\fontsize{6}{7}\selectfont
\caption{Hyperparameters selected in Sections \ref{1D}-\ref{3D}, unless otherwise noted.}
\label{tab:exp_parameters}
\begin{tabular}{cc|ccccccc}
\toprule
\textbf{Problem} & \textbf{Method} & \textbf{Depth} & \textbf{Width} & \textbf{[$N_f$/$N_h$, $N_b$/$N_{hp}$, $N_i$]} & \textbf{Iterations} & \textbf{Optimizer} \\
\midrule
\multirow{2}{*}{1D--\ref{1D}}
  & PINN/gPINN & 4 & 10 & [2000, 2000, 2000] & 20000 & Adam \\
  & DAE  & 4 & 10 & [1000, --, --] & 8000 & Adam \\
\midrule
\multirow{2}{*}{2D--\ref{2D}}
  & PINN/gPINN & 5 & 10 & [3000, 1000, 3000] & 30000 & Adam  \\
  & DAE  & 5 & 10 & [2000, 1000, --] & 15000 & Adam\\
\midrule
\multirow{2}{*}{3D--\ref{3D}}
  & PINN/gPINN & 6 & 10 & [6000, 6000, 6000] & 40000 & Adam \\
  & DAE  & 6 & 10 & [4000, 2000, --]  & 20000 & Adam \\
\bottomrule
\end{tabular}
\end{table}

\subsection{1D problem}\label{1D}

Consider problem \eqref{mainproblem} on the unit interval $\Omega_1=(0,1)$ and final time $T>0$:
\begin{align}\label{forwardexample1}
\begin{cases}
\mu_1 \frac{\partial^2 u}{\partial x^2} - \frac{\partial u}{\partial t} = -u \frac{\partial u}{\partial x} + f(x), & \quad (x,t) \in \Omega_1 \times [0,T], \\
u(0,t) = L, \, u(1,t) = R, & \quad t \in [0,T], \\
u(x,0) = {u}_0 (x,\mu_1), &\quad x\in\Omega_1,
\end{cases}
\end{align}
with $ L = -10 $, $ R = 5 $, $T=0.3$, a monotonically increasing source $ f(x)=x - x^2 + x^3 $, and ${u}_0 (x,\mu_1) = \frac{R - L}{2} \tanh(\frac{x - 0.1}{\mu_1}) + \frac{R + L}{2}$.

To ensure the existence of an asymptotic solution, we verify Assumption \ref{assumption}. By solving \eqref{zeroorderregularequation1}, we obtain the leading regular terms
$\varphi^{(-)} (x) = - \frac{\sqrt{600 + 6x^2 - 4x^3 + 3x^4}}{\sqrt{6}}$ and $
\varphi^{(+)} (x) = \frac{\sqrt{145 + 6x^2 - 4x^3 + 3x^4}}{\sqrt{6}}$. The leading term $ h_0(t) $ in \eqref{integral_eq_kdim} satisfies
\begin{align*}
\label{eq44Ex1}
\begin{cases}
\frac{dh_0(t)}{dt} = - \frac{1}{2} \left( \varphi^{(-)}(h_0(t)) + \varphi^{(+)}(h_0(t)) \right), \\
h_0(0) = 0.1.
\end{cases}
\end{align*}
Thus, if Assumption \ref{assumption} holds, equation~\eqref{forwardexample1} has the following solution:
{\small
\begin{align*} 
U_{0}(x,t)=\begin{cases}
\displaystyle \varphi^{(-)} (x) +\frac{ \left(\varphi^{(+)}(h_{0}(t))-\varphi^{(-)}(h_{0}(t)) \right)}{  \exp \left( \left(x-h_0(t) \right) \left( \frac{\varphi^{(-)}(h_{0}(t))-\varphi^{(+)}(h_{0}(t)) }{2\mu_1} \right) \right)+ 1}  ,\,\, x \in [0,h_0(t)],\\
\displaystyle \varphi^{(+)} (x)+\frac{ -\left(\varphi^{(+)}(h_{0}(t))-\varphi^{(-)}(h_{0}(t))\right)}{  \exp \left(\left(x-h_0(t) \right)  \left( \frac{\varphi^{(+)}(h_{0}(t))-\varphi^{(-)}(h_{0}(t)) }{2\mu_1}\right) \right)+ 1} ,\,\, x \in [h_0(t),1].
\end{cases}
\end{align*}
}Assumptions \ref{assumption}(a), (b) and (d) are satisfied by the choice of $A(u),L,R $ and ${u}_0$. Hence, we only need to check Assumptions \ref{assumption}(c) and (e).

\begin{figure}[htbp]
\centering
\subfigure[approximation of $h_0$ \label{fig:1Dh0}]{
  \includegraphics[width=0.35\linewidth]{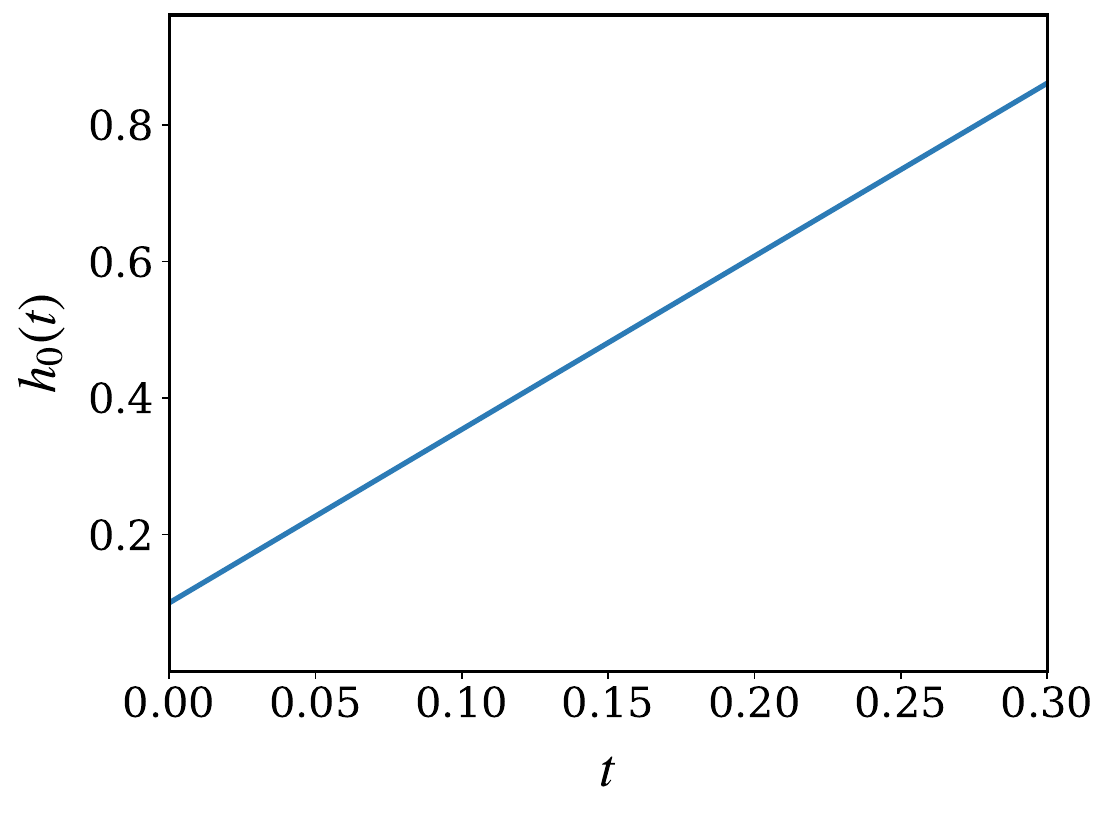}
}\hspace{-1ex}
\subfigure[approximation of ${h}_1$ \label{fig:1Dx1}]{
  \includegraphics[width=0.35\linewidth]{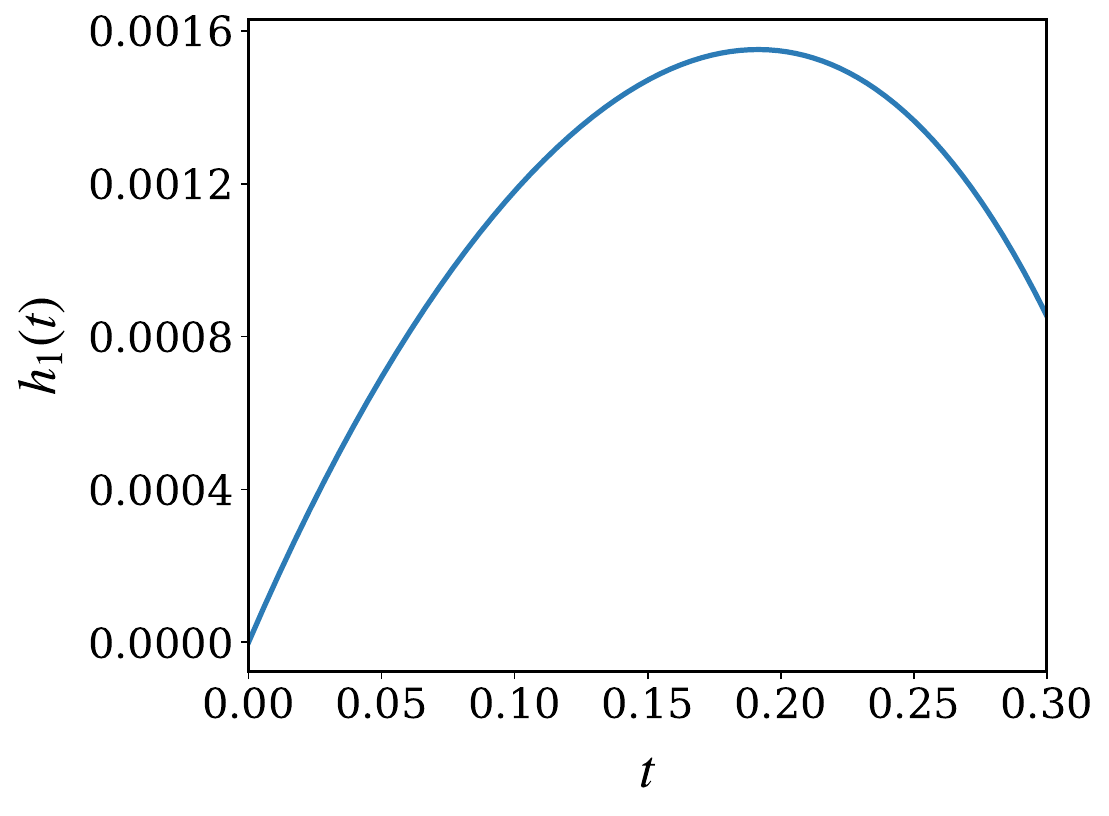}
}
\caption{1D forward problem}
\end{figure}

Fig.~\ref{fig:1Dh0} presents the approximate solution of ${h}_0(t)$ by the DAE, which numerically confirms that $0 < {h}_0(t) < 1$ for all $t \in \bar{\mathcal{T}}$ thereby validating Assumption~\ref{assumption}(c).  For the function $A=-u$, using the explicit formula for $({h_0}(t))' $ from  \eqref{integral_eq_kdim}, Assumption~\ref{assumption}(e) simplifies to:
$
\int_{\varphi^{(-)}(h_0)}^{s}
(-u-({h_0}(t))')du = \frac{1}{2}( s-\varphi^{(-)})( \varphi^{(+)}-s)>0.
$
For any
$s\in(\varphi^{(-)},\varphi^{(+)})$, both \((s-\varphi^{(-)})\) and \((\varphi^{(+)}-s)\) are positive. Hence, Assumption~\ref{assumption}(e) holds.

\begin{figure}[htbp]
\centering
\subfigure[]{
\includegraphics[width=0.25\linewidth]{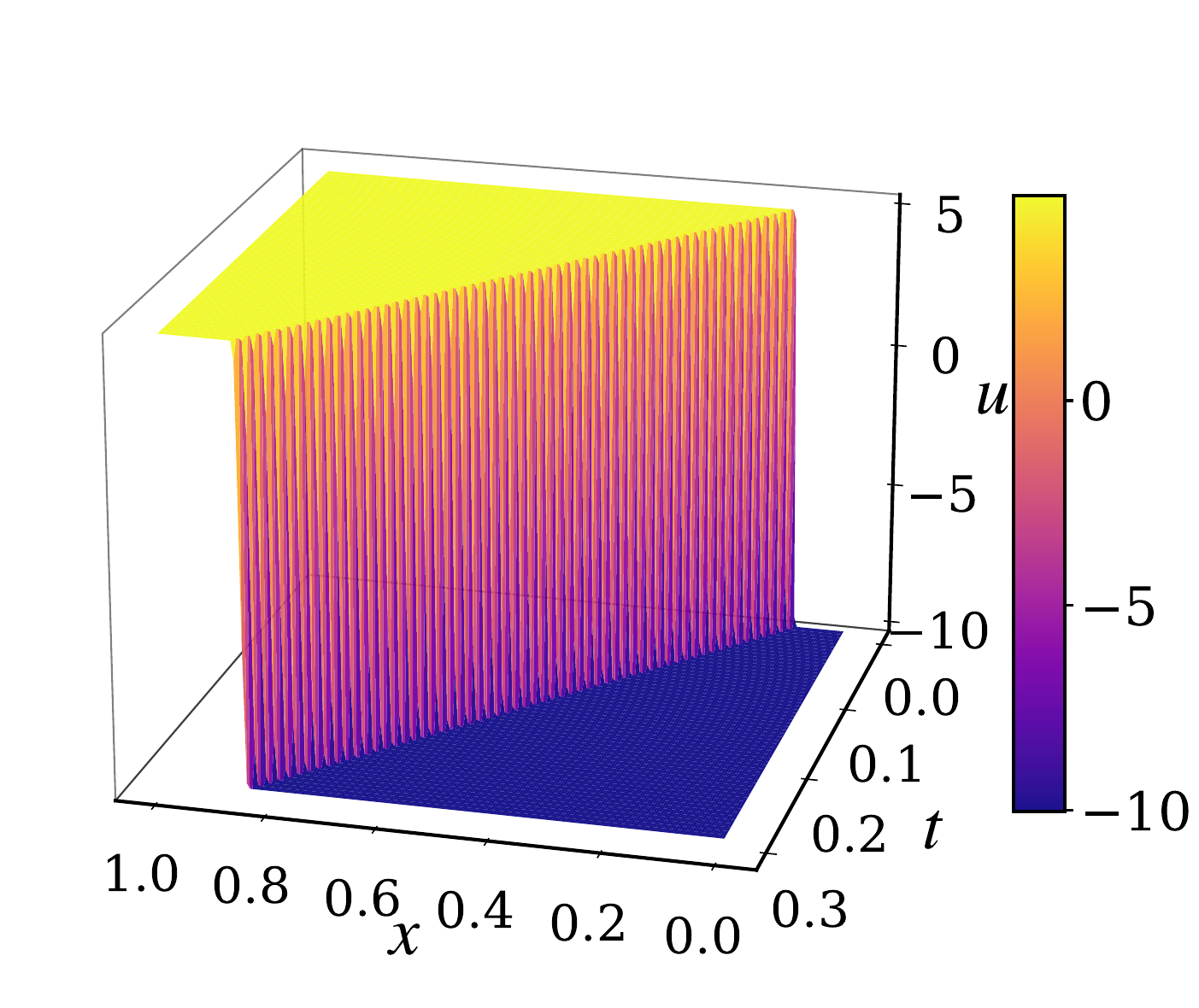} \label{fig:1d_true_solution_mu2} }
\subfigure[]{
\includegraphics[width=0.25\linewidth]{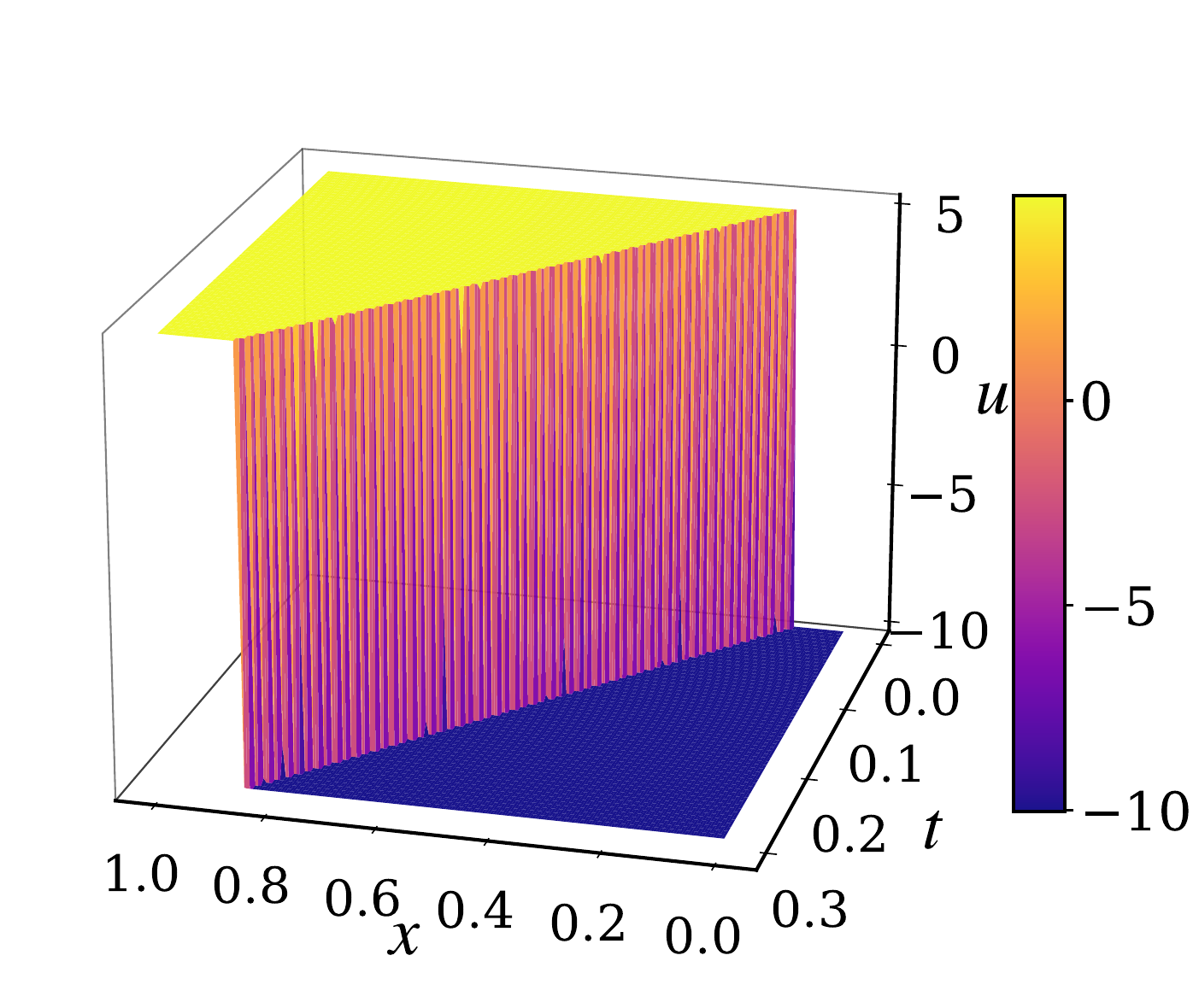} \label{fig:1d_true_solution_mu3} }
\subfigure[]{
\includegraphics[width=0.25\linewidth]{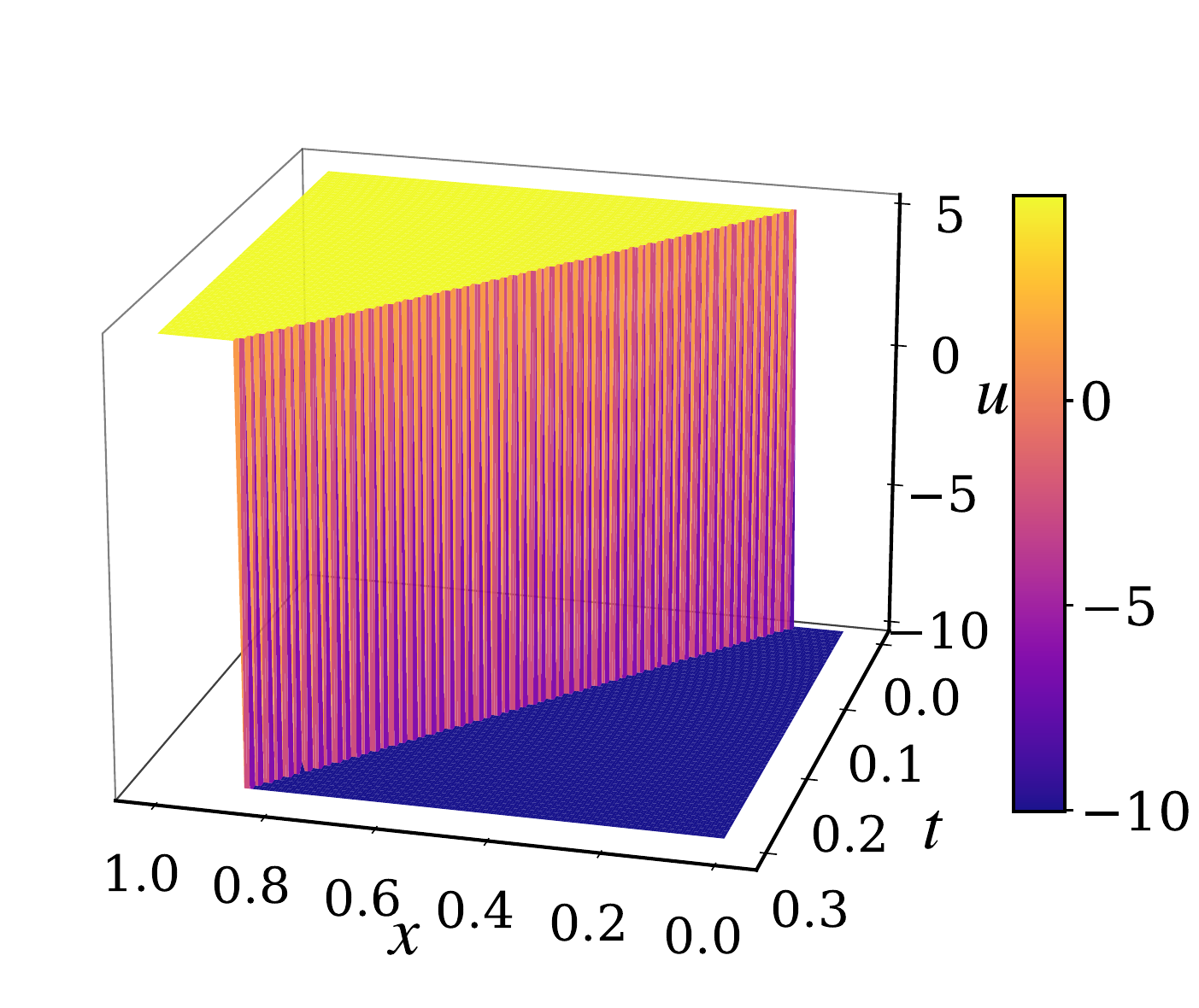} \label{fig:1d_true_solution_mu4} }
\vspace{-1ex}
\subfigure[]{
\includegraphics[width=0.25\linewidth]{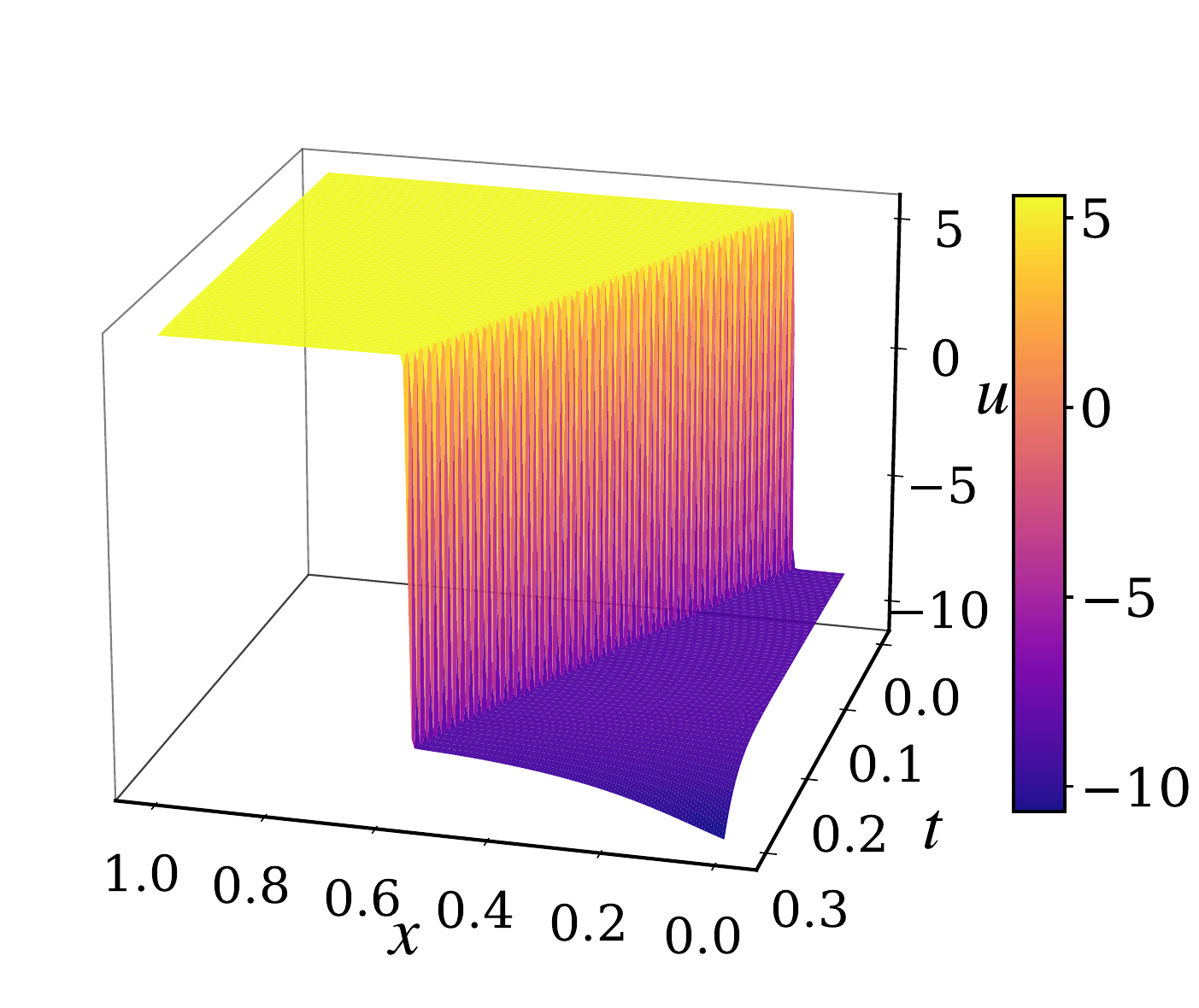} \label{fig:1d_PINN_solution_mu2} }
\subfigure[]{
\includegraphics[width=0.25\linewidth]{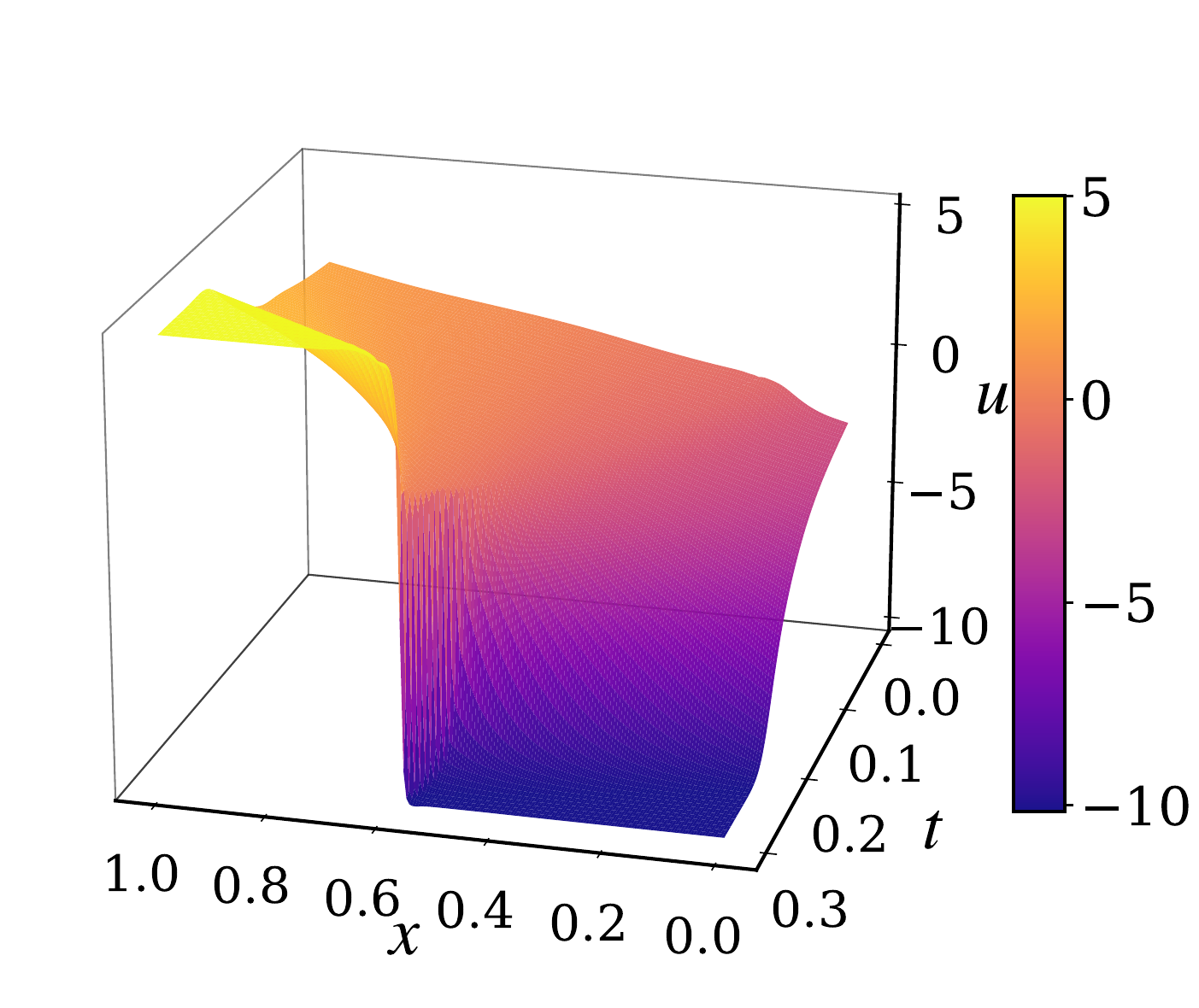} \label{fig:1d_PINN_solution_mu3} }
\subfigure[]{
\includegraphics[width=0.25\linewidth]{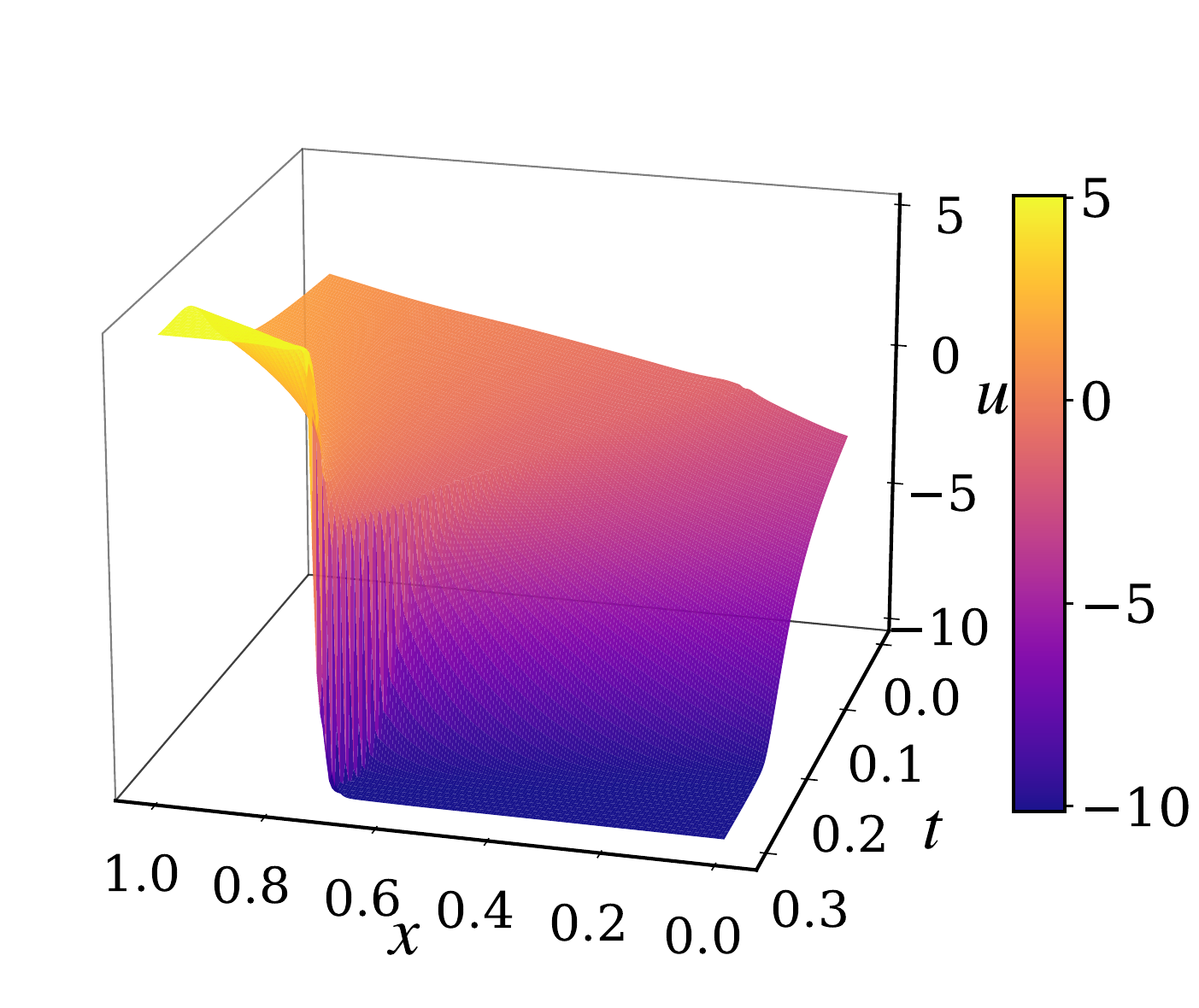} \label{fig:1d_PINN_solution_mu4} }
\vspace{-1ex}
\subfigure[]{
\includegraphics[width=0.25\linewidth]{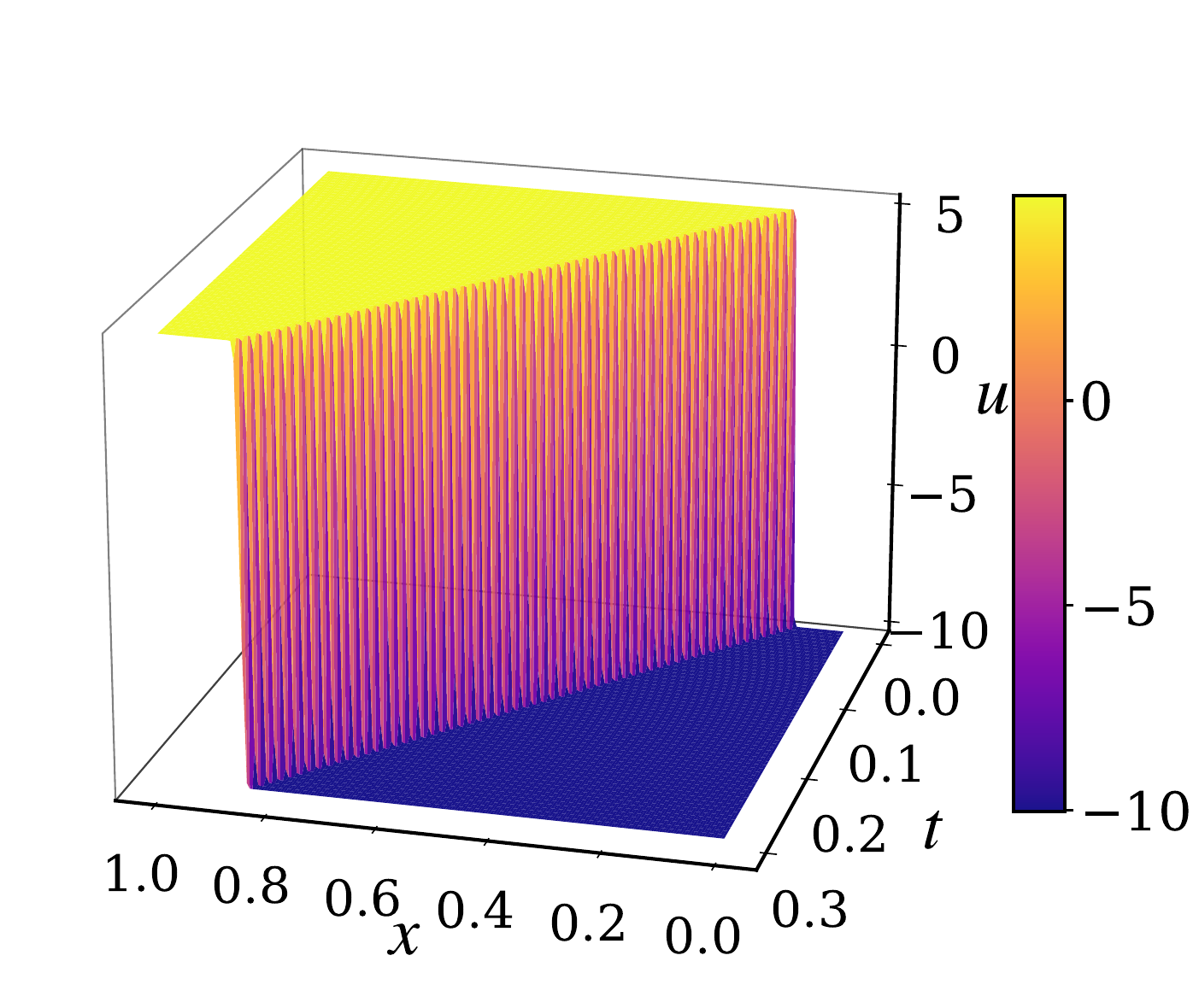} \label{fig:1d_DAE_solution_mu2} }
\subfigure[]{
\includegraphics[width=0.25\linewidth]{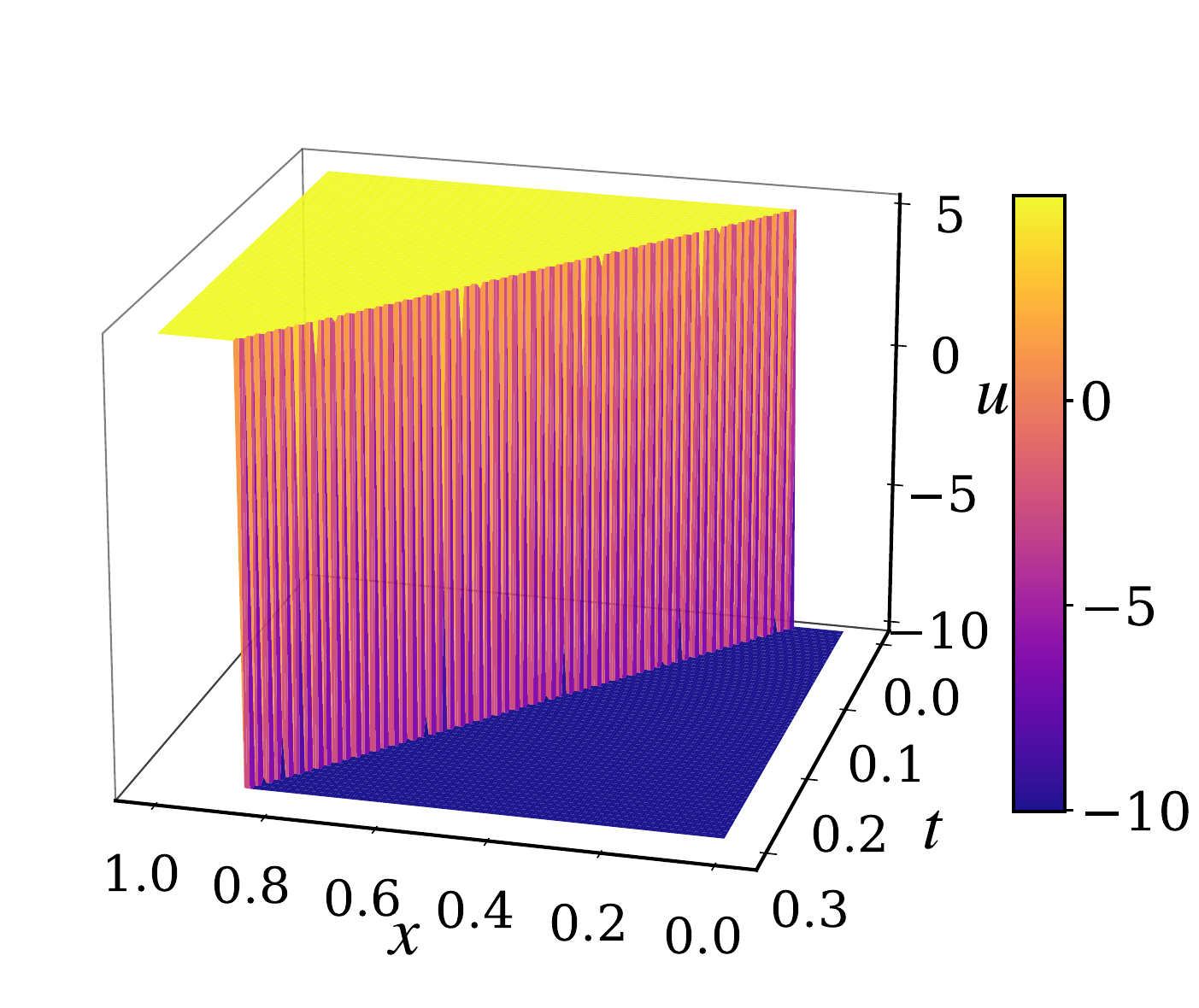} \label{fig:1d_DAE_solution_mu3} }
\subfigure[]{
\includegraphics[width=0.25\linewidth]{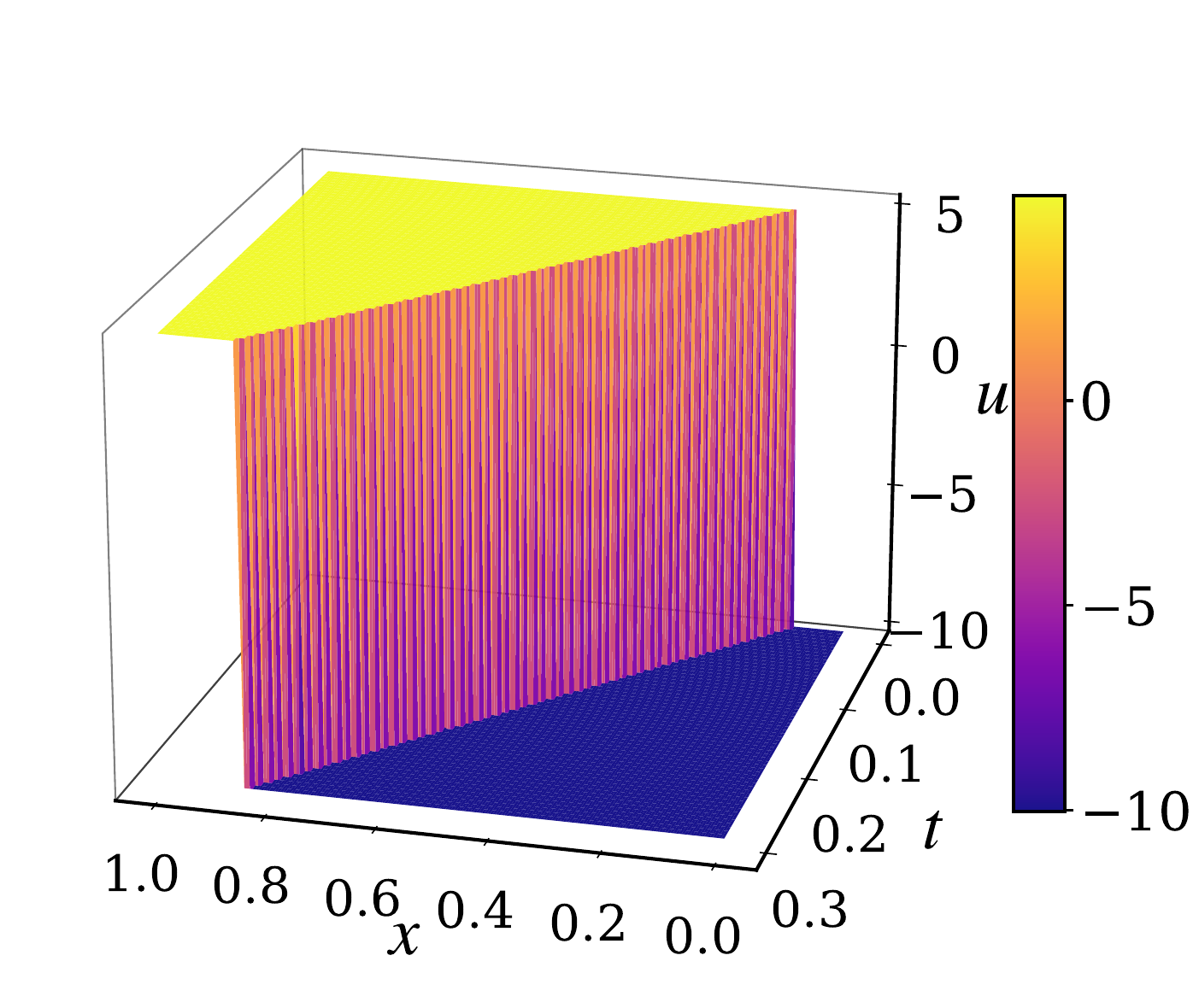} \label{fig:1d_DAE_solution_mu4} }
\caption{1D problem: Reference solutions \subref{fig:1d_true_solution_mu2} \subref{fig:1d_true_solution_mu3} \subref{fig:1d_true_solution_mu4} and predicted solutions  obtained by PINN \subref{fig:1d_PINN_solution_mu2} \subref{fig:1d_PINN_solution_mu3} \subref{fig:1d_PINN_solution_mu4} and DAE \subref{fig:1d_DAE_solution_mu2}  \subref{fig:1d_DAE_solution_mu3}  \subref{fig:1d_DAE_solution_mu4} for $\mu_1 =10^{-2}, 10^{-3}, 10^{-4}$, respectively.}
\label{fig: 1D forward problem solution}
\end{figure}

\begin{figure}[htbp]
\vspace{0ex}
\centering
\subfigure[]{
\includegraphics[width=0.35\linewidth]{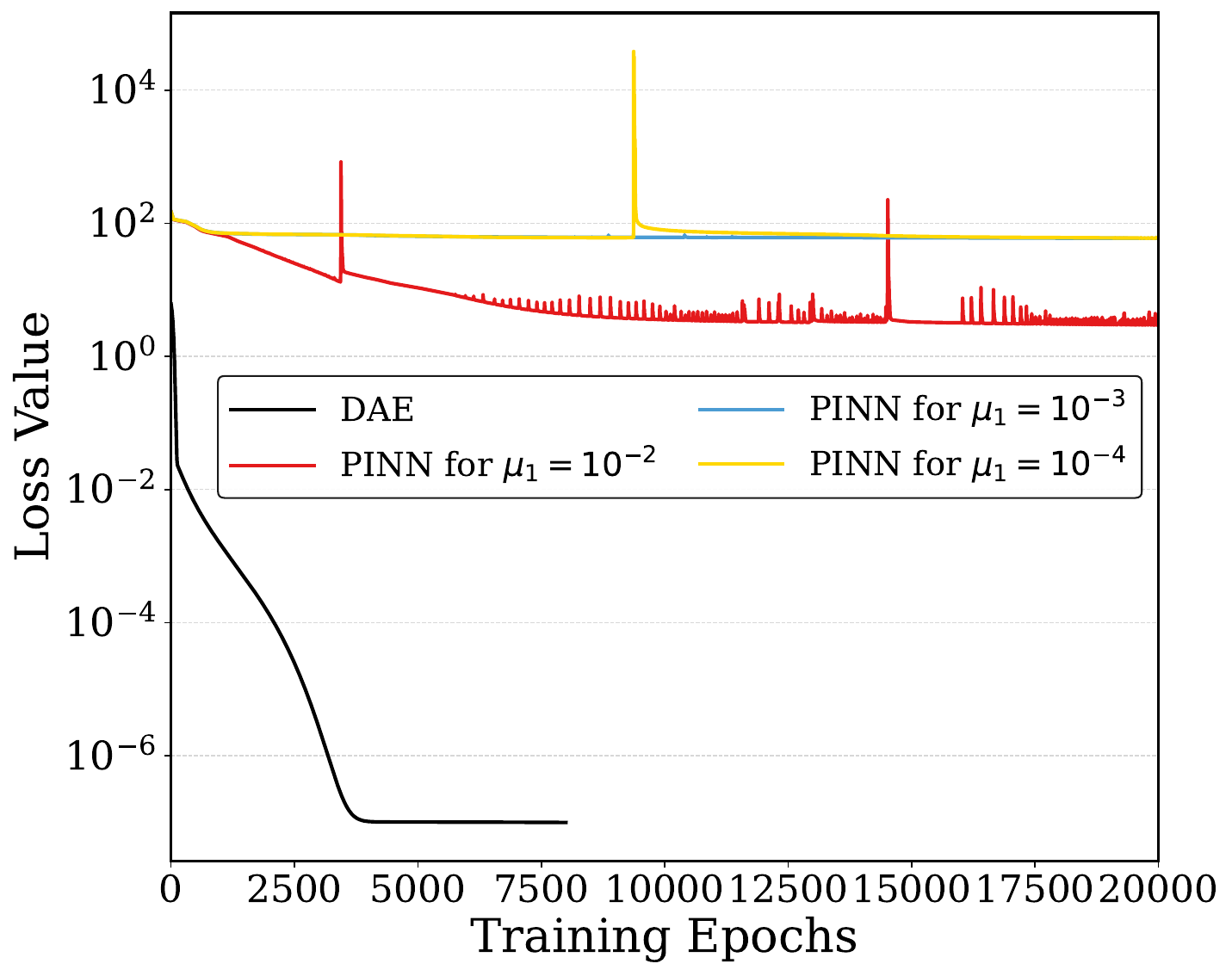} \label{fig:1d_DAE_loss} }
\subfigure[]{
\includegraphics[width=0.35\linewidth]{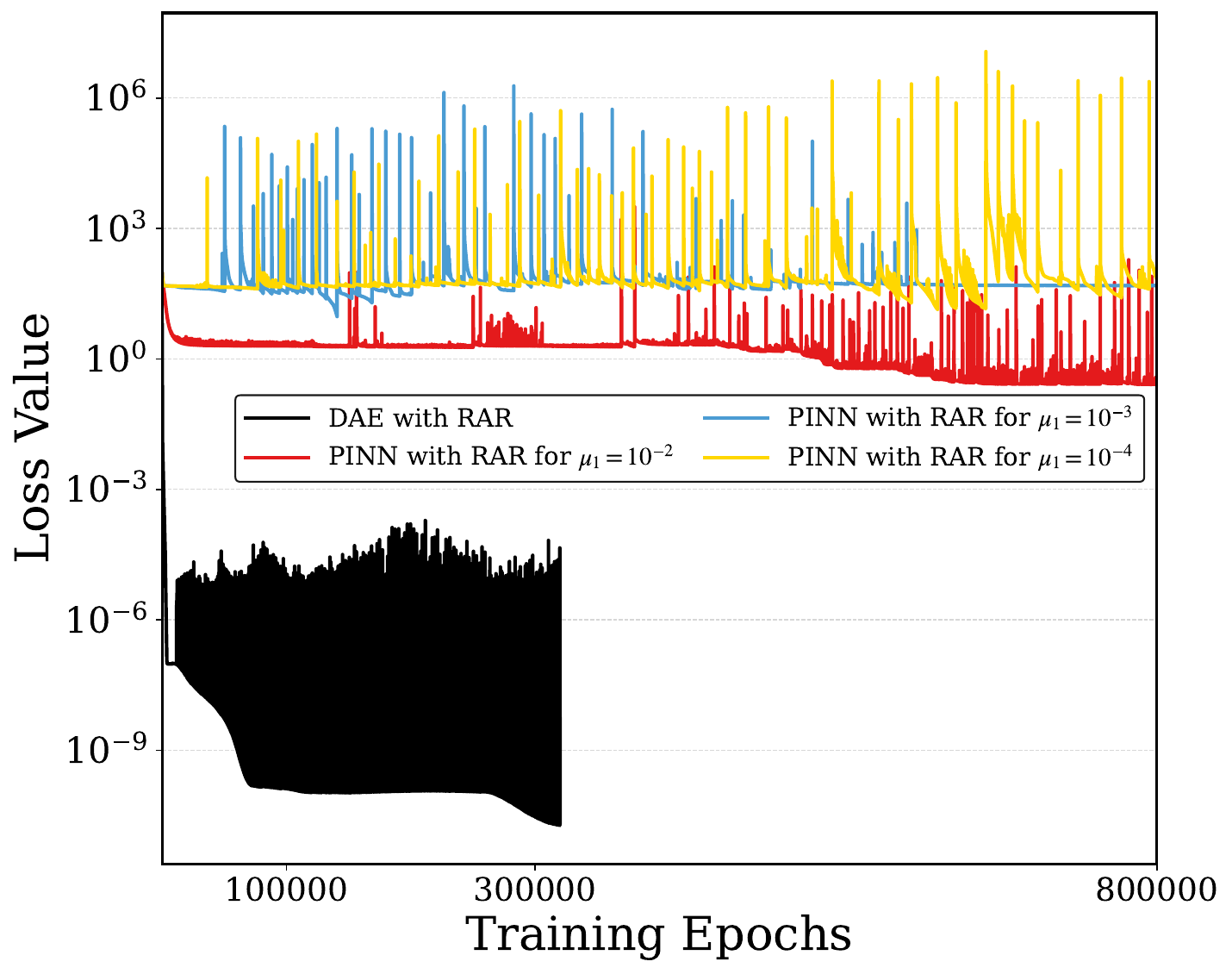} \label{fig:1d_DAE_RAR_loss} }
\caption{1D problem: Loss curves of PINN and DAE \subref{fig:1d_DAE_loss}, PINN with RAR and DAE with RAR \subref{fig:1d_DAE_RAR_loss}. The loss function of DAE is independent of \(\mu_1\).} \label{1D forward problem loss}
\end{figure}

For DAE/PINN with RAR, 800/1800 initial residual points are randomly selected, and 200/200 residual points are adaptively added with $m=5$, $|\mathcal{S}|=5000/10000$, and $\mathcal{E}_0=10^{-6}$. Table~\ref{tab:1D forward problem} and Fig.~\ref{fig: 1D forward problem solution} show that the standard PINN fails to converge and the approximation has large errors, whereas the DAE can achieve good accuracy. Table \ref{tab:1D forward problem} indicates that RAR does not improve the accuracy of PINN, while the DAE with RAR can significantly improve the accuracy.
Fig. \ref{1D forward problem loss} shows that the DAE converges faster and achieves lower loss values. The DAE consistently outperforms PINN and gPINN: the relative \(L^2\) errors is 1000 times smaller, and the DAE enjoys remarkable robustness, cf. Tables~\ref{tab:training data_1}, \ref{tab:neurons_1} and \ref{tab:random seeds_1}.

The DAE method for the first-order approximate solution, denoted as $\text{DAE}^1$, is given in Appendix \ref{appendix:DAE1}. Fig.~\ref{fig:1Dx1} shows the approximation of 
${h}_1(t)$ by $\text{DAE}^1$. Table~\ref{tab:1D forward problem} shows that the first-order asymptotic does not improve much over the zeroth-order DAE, but the training time is more than doubled. Thus we focus only on the zeroth-order asymptotic solution.

\begin{table}[htbp]
    \centering
    \fontsize{6}{7}\selectfont
    \caption{1D  problem: The loss values, $L^2$,  $L^{\infty}$ errors and training times of different methods for various $\mu_1$ values. Dashes (--) denote values identical to those for $\mu_1 = 10^{-2}$.}
    \begin{tabular}{c|ccccc}
        \toprule
        $\mu_1 $ & Method & $e_{loss}$ & $e_2$ & $e_{\infty}$ & Time (s) \\
        \midrule
        \multicolumn{6}{c}{without RAR} \\
        \midrule
        \multirow{2}{*}{$10^{-2}$}
        & PINN   & 3.04e+00 & 7.79e-01 & 1.56e+01 & 254 \\
        & \textbf{DAE}    & \textbf{9.92e-08} & \textbf{2.16e-04} & \textbf{4.53e-02} & \textbf{83} \\
        & \textbf{$\text{DAE}^1$}   & \textbf{1.42e-05} & \textbf{2.16e-04} & \textbf{4.51e-02} & \textbf{206} \\
        \cmidrule(r){1-6}
        \multirow{2}{*}{$10^{-3}$}
        & PINN   & 5.94e+01 & 7.36e-01 & 1.50e+01 & 266 \\
        & \textbf{DAE}    & -- & \textbf{9.19e-04} & \textbf{4.47e-01} & -- \\
        & \textbf{$\text{DAE}^1$} & -- & \textbf{9.16e-04} & \textbf{4.45e-01} & -- \\
        \cmidrule(r){1-6}
        \multirow{2}{*}{$10^{-4}$}
        & PINN   & 5.97e+01 & 6.40e-01 & 1.50e+01 & 272 \\
        & \textbf{DAE}    & -- & \textbf{2.99e-03} & \textbf{1.61e+00} & -- \\
        & \textbf{$\text{DAE}^1$} & -- & \textbf{2.98e-03} & \textbf{1.60e+00} & --\\
        \midrule
        \multicolumn{6}{c}{with RAR} \\
        \midrule
        \multirow{2}{*}{$10^{-2}$}
        & PINN+RAR & 2.71e-01 & 1.63e-01 & 1.47e+01 & 11515 \\
        & \textbf{DAE+RAR}  & \textbf{2.08e-11} & \textbf{2.92e-06} & \textbf{4.90e-04} & \textbf{3218} \\
        \cmidrule(r){1-6}
        \multirow{2}{*}{$10^{-3}$}
        & PINN+RAR & 4.83e+01 & 8.00e-01 & 1.51e+01 & 12797 \\
        & \textbf{DAE+RAR}  & -- & \textbf{6.72e-06} & \textbf{3.26e-03} & -- \\
        \cmidrule(r){1-6}
        \multirow{2}{*}{$10^{-4}$}
        & PINN+RAR & 6.82e+01 & 1.09e+00 & 1.47e+01 & 13184 \\
        & \textbf{DAE+RAR}  & -- & \textbf{3.29e-05} & \textbf{1.78e-02} & -- \\
        \bottomrule
    \end{tabular}
    \label{tab:1D forward problem}
\end{table}

\begin{table}[htbp]
\centering
\fontsize{6}{7}\selectfont
\caption{1D problem: Relative $L^2$ errors for $\mu_1 = 10^{-2}$ with varying training data. A 4-layer, 10-neuron network is used with 10000 iterations. Dashes (--) mark values equal to those at \(N_b/N_i = 500\).}
\vspace{0cm} 
\begin{tabular}{c@{\hskip 2pt}c@{\hskip 3pt}c@{\hskip 3pt}c@{\hskip 4pt}c@{\hskip 3pt}c@{\hskip 3pt}c@{\hskip 4pt}c@{\hskip 3pt}c@{\hskip 3pt}c@{\hskip 4pt}c}
\toprule
\multirow{2}{*}{$N_b/N_i$} & \multicolumn{9}{c}{$N_f/N_h$} \\ 
\cmidrule(lr){2-11}
& & \multicolumn{3}{c}{800} & \multicolumn{3}{c}{1600} & \multicolumn{3}{c}{2400} \\
\cmidrule(lr){3-5} \cmidrule(lr){6-8} \cmidrule(lr){9-11} 
&  & PINN &gPINN& DAE & PINN&gPINN & DAE & PINN&gPINN & DAE \\
\midrule
500  &   & 6.29e-01& 7.74e-01  & \textbf{2.16e-04} & 7.23e-01& 7.66e-01 & \textbf{2.10e-04} & 7.36e-01& 7.61e-01 & \textbf{2.09e-04} \\
1000 &   & 5.56e-01& 7.75e-01  & -- & 5.34e-01 & 7.67e-01& -- & 7.46e-01 & 7.68e-01& -- \\
1500 &   & 6.33e-01 & 7.78e-01 & -- & 7.57e-01 & 7.68e-01& -- & 7.40e-01 & 7.68e-01& -- \\
\bottomrule
\end{tabular}
\vspace{0cm} 
\label{tab:training data_1}
\end{table}

\begin{table}[htbp]
\centering
\fontsize{6}{7}\selectfont
\caption{1D forward problem: Relative $L^2$ errors for $\mu_1 = 10^{-2}$ with different network sizes. Training data fixed at $N_f/N_h = 1000$, $N_b/N_i = 1000$, with 10000 iterations.}
\vspace{0cm} 
\begin{tabular}{c@{\hskip 2pt}c@{\hskip 3pt}c@{\hskip 3pt}c@{\hskip 4pt}c@{\hskip 3pt}c@{\hskip 3pt}c@{\hskip 4pt}c@{\hskip 3pt}c@{\hskip 3pt}c@{\hskip 4pt}c}
\toprule
\multirow{2}{*}{Layers} & \multicolumn{9}{c}{Neurons} \\ 
\cmidrule(lr){2-11}
& & \multicolumn{3}{c}{8} & \multicolumn{3}{c}{12} & \multicolumn{3}{c}{15} \\
\cmidrule(lr){3-5} \cmidrule(lr){6-8} \cmidrule(lr){9-11} 
&  & PINN&gPINN & DAE & PINN&gPINN & DAE & PINN&gPINN & DAE \\
\midrule
3  &   & 6.89e-01 & 7.86e-01 & \textbf{3.37e-04} & 6.65e-01& 7.71e-01 & \textbf{4.53e-04} & 7.78e-01& 7.65e-01 & \textbf{4.44e-04} \\
4  &   & 7.53e-01 & 7.64e-01 & \textbf{8.46e-04} & 6.07e-01& 7.41e-01 & \textbf{3.03e-04} & 9.21e-01 & 7.60e-01& \textbf{3.07e-04} \\
5  &   & 6.42e-01 & 7.56e-01 & \textbf{4.68e-04} & 7.83e-01& 7.59e-01 & \textbf{3.96e-04} & 9.64e-01& 7.51e-01 & \textbf{4.88e-04} \\
\bottomrule
\end{tabular}
\vspace{0cm} 
\label{tab:neurons_1}
\end{table}

\begin{table}[htbp]
	\centering
	\fontsize{6}{7}\selectfont    
	\caption{1D forward problem: Relative $L^2$ errors for $\mu_1 = 10^{-2}$ with different random seeds.}
	\begin{tabular}{lccccc}
		\toprule
		\diagbox [width=8em,trim=l] {Methods}{{Seed}} & 33 & 99 & 202 & 5678 & 9999 \\ \midrule
		PINN & 7.60e-01 & 6.59e-01 & 7.53e-01 & 8.82e-01 & 7.27e-01\\
		gPINN & 7.63e-01 & 7.61e-01 & 7.71e-01 & 7.65e-01 & 7.62e-01\\
        DAE & \textbf{3.41e-04} & \textbf{7.30e-04} & \textbf{3.48e-04} & \textbf{6.19e-04} & \textbf{2.00e-04}\\ \bottomrule
	\end{tabular}
	\vspace{0cm}
	\label{tab:random seeds_1}
\end{table}

\subsection{2D problem}\label{2D}
Consider the 2D problem subject to periodic boundary conditions along the $y$ axis, with the source $f(x,y)=\cos{ ( \pi x /4)} \cos{(\pi y /4)}$:
\begin{align} \label{forwardexample2}
\begin{cases}
\displaystyle \mu_2 \Delta u - \partial_t u =  -u \left(\partial_x u + \partial_y u\right) +f(x,y), \\
u(x, y, t)=u(x, y+4, t),  \ u(-2,y,t)=-4, \ u(2,y,t)=2, \\
u(x,y,0)={u}_0(x,y,\mu_2 ),\\
 x \in [-2, 2], \  y \in [-2, 2], \  t \in [0,1].
\end{cases}
\end{align}
The left and right outer functions satisfy
\begin{equation*}
\begin{aligned}
&\begin{cases}
\displaystyle \varphi^{(\pm)} \left(\partial_x \varphi^{(\pm)}+\partial_y \varphi^{(\pm)}\right)=f(x, y),\\
\varphi^{(\pm)} (\pm2, y)=-1\pm3, \quad \varphi^{(\pm)} (x, y)=\varphi^{(\pm)} (x, y+4);
\end{cases}
\end{aligned}
\end{equation*}
The zeroth-order outer asymptotics are given by
\begin{align*}
\varphi^{(\pm)}(x,y) &= \pm\frac{1}{\sqrt{\pi}} \left[
    (10 \mp 6)\pi + \mp 2\pi \cos\left( \frac{\pi}{4}(-x + y) \right)
    + \pi x \cos\left( \frac{\pi}{4}(-x + y) \right) \right. \\
&\qquad\left. - 2 \sin\left( \frac{\pi}{4}(\pm4 - x + y) \right)
    + 2 \sin\left( \frac{\pi}{4}(x + y) \right)
\right]^{1/2}.
\end{align*}
In the 2D problem, we have
\begin{align} \label{h0mainequation}
\begin{cases}
\displaystyle \partial_t{h_0} =\left( \partial_y{h_0}-\tfrac{1}{2} \right) \left( \varphi^{(-)} (h_0,y)+\varphi^{(+)} (h_0,y)\right),\\
{h_0}(y, t)={h_0}(y+4, t), \quad h_0(y,0)=h_0^{*}=0.
\end{cases}
\end{align}

From the approximation of \eqref{h0mainequation} by the DAE, we determine that the transition layer is located within the region $-2 \leq h_0 (y,t) \leq 2 $ for $(y,t) \in [-2,2] \times [0,1]$. Thus Assumption~\ref{assumption}(c) is satisfied. The initial condition ${u}_0(x,y,\mu_2 ) =\frac{R-L}{2} \tanh(\frac{x-h_0^{*}}{\mu_2}+y)+ \frac{R+L}{2}= 3\tanh(\frac{x}{\mu_2}+y)-1$ satisfies Assumption \ref{assumption}(d). Thus, the solution to problem~\eqref{forwardexample2} is in the form of an autowave with a transitional moving layer localized near $h_0(y,t)$. The zeroth-order asymptotic is given by
{\small
\begin{align*} \label{asymptoticsolEXAMPLE2}
U_{0}=\begin{cases}
\displaystyle \varphi^{(-)} (x,y) +\frac{\varphi^{(+)}(h_0,y)-\varphi^{(-)}(h_0,y)}{\exp \left(\frac{\left(h_0-x \right) \left(\varphi^{(+)}(h_0,y)-\varphi^{(-)}(h_0,y)\right) \left(1-\partial_y{h_0} \right) }{2\mu_2 }  \right)+1}, \  x \in [-2,h_0],\\
\displaystyle \varphi^{(+)} (x,y)+\frac{\varphi^{(-)}(h_0,y)-\varphi^{(+)}(h_0,y)}{\exp \left(\frac{\left(h_0-x \right) \left(\varphi^{(-)}(h_0,y)-\varphi^{(+)}(h_0,y)\right) \left(1-\partial_y{h_0} \right) }{2\mu_2 }  \right)+1}, \  x \in [h_0,2].
\end{cases}
\end{align*}
}

\begin{figure}[htbp]
\centering
\subfigure[]{
\includegraphics[width=0.25\linewidth]{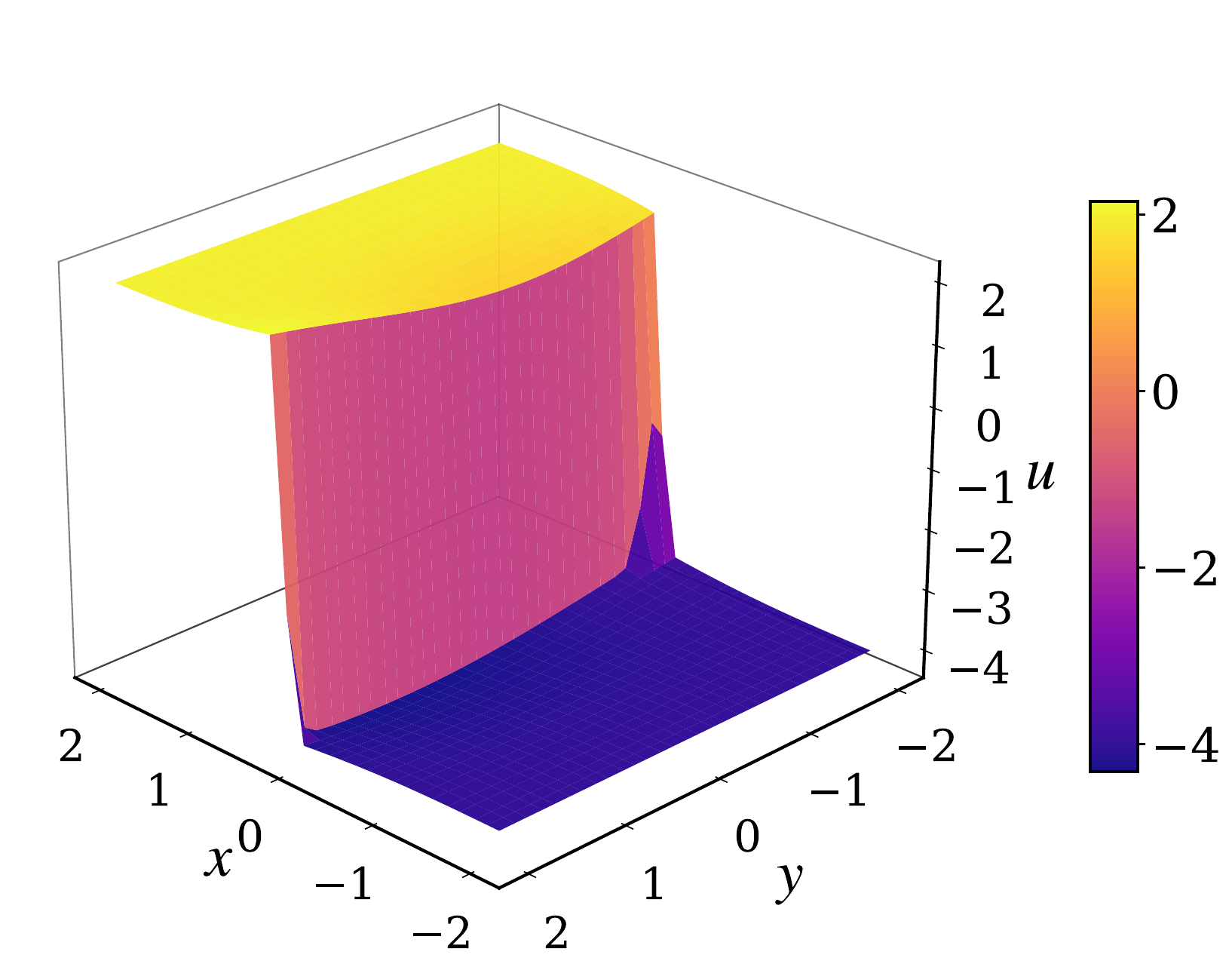} \label{fig:2d_true_solution_mu2} }
\subfigure[]{
\includegraphics[width=0.25\linewidth]{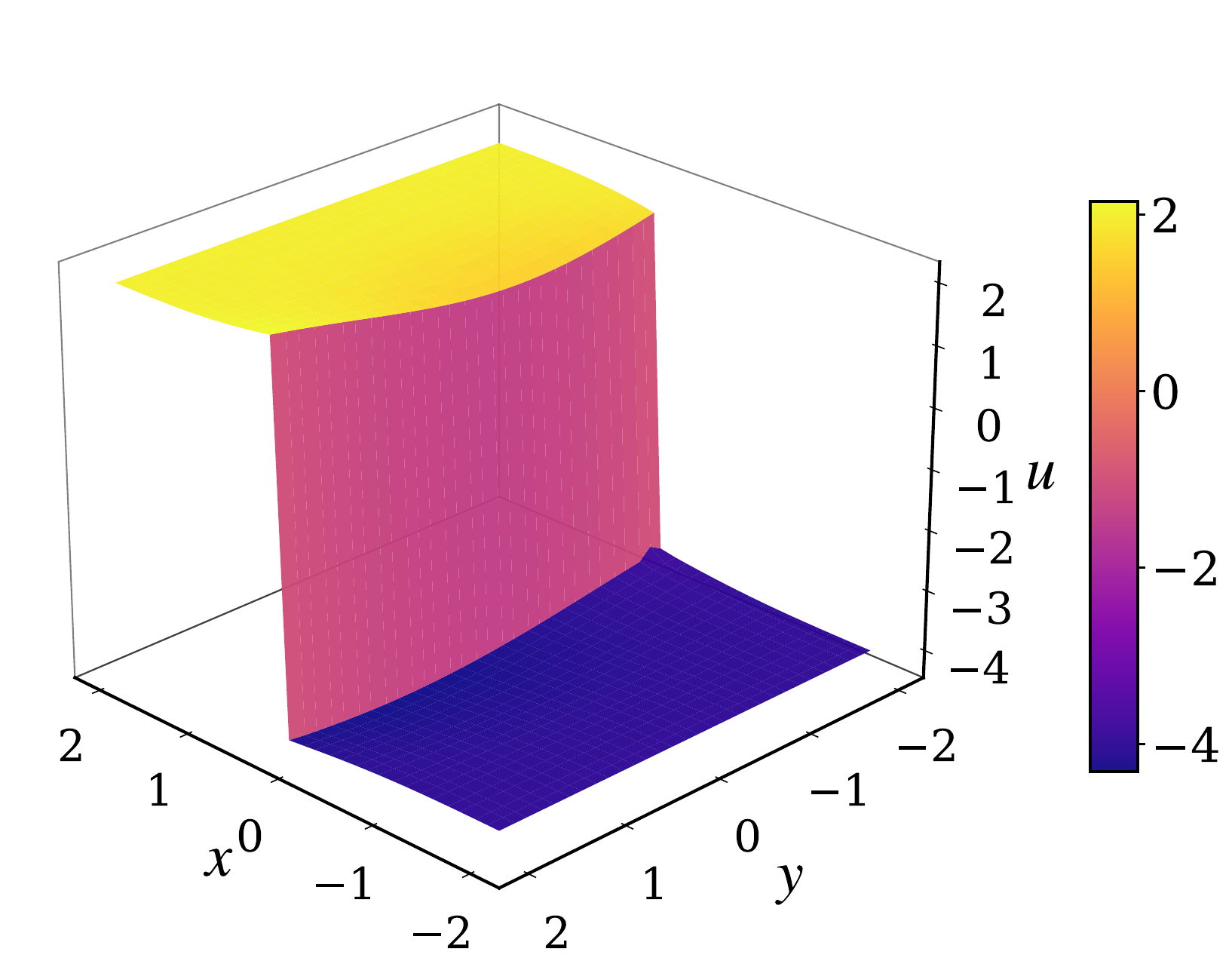} \label{fig:2d_true_solution_mu3} }
\subfigure[]{
\includegraphics[width=0.25\linewidth]{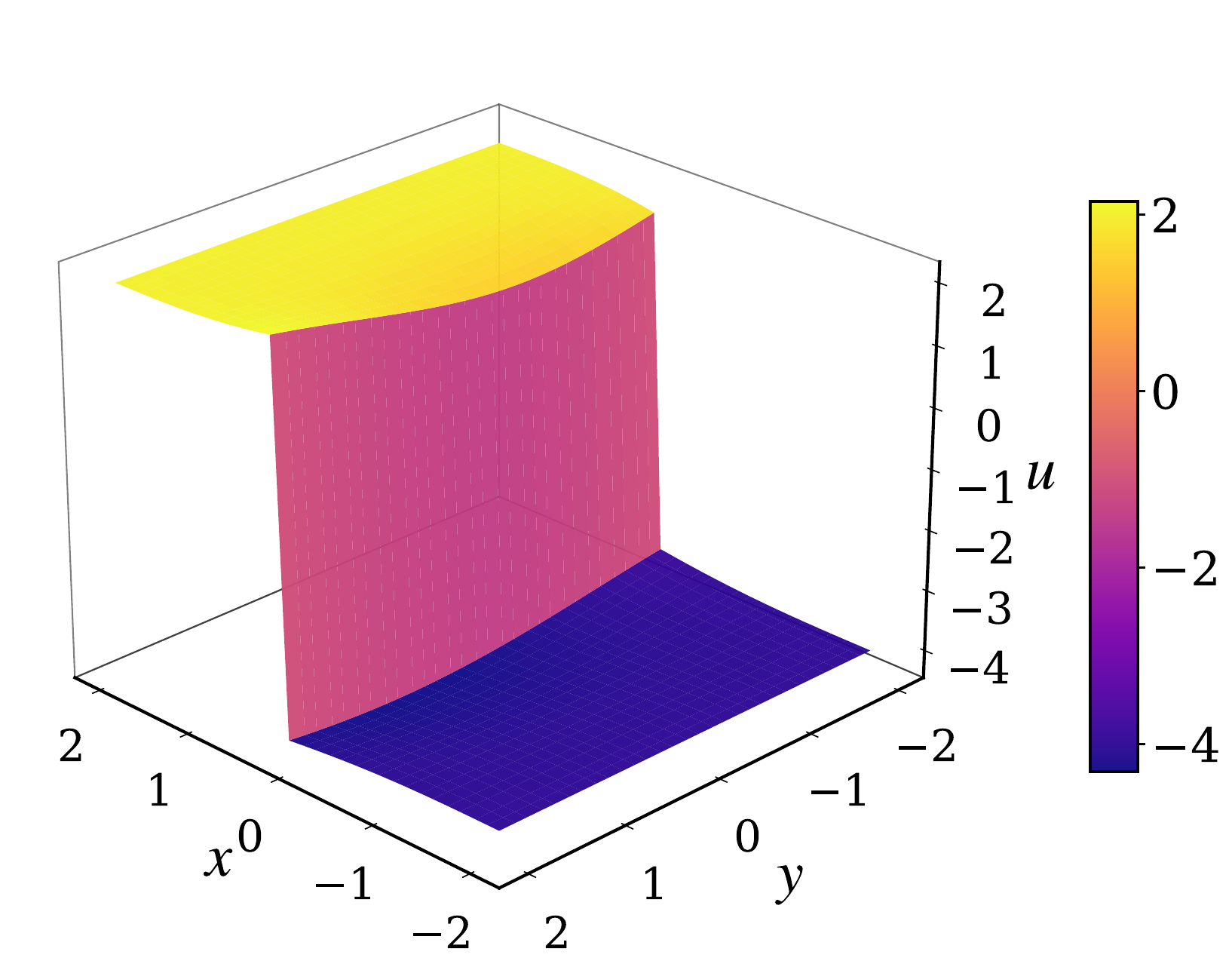} \label{fig:2d_true_solution_mu4} }
\vspace{-1ex}
\subfigure[]{
\includegraphics[width=0.25\linewidth]{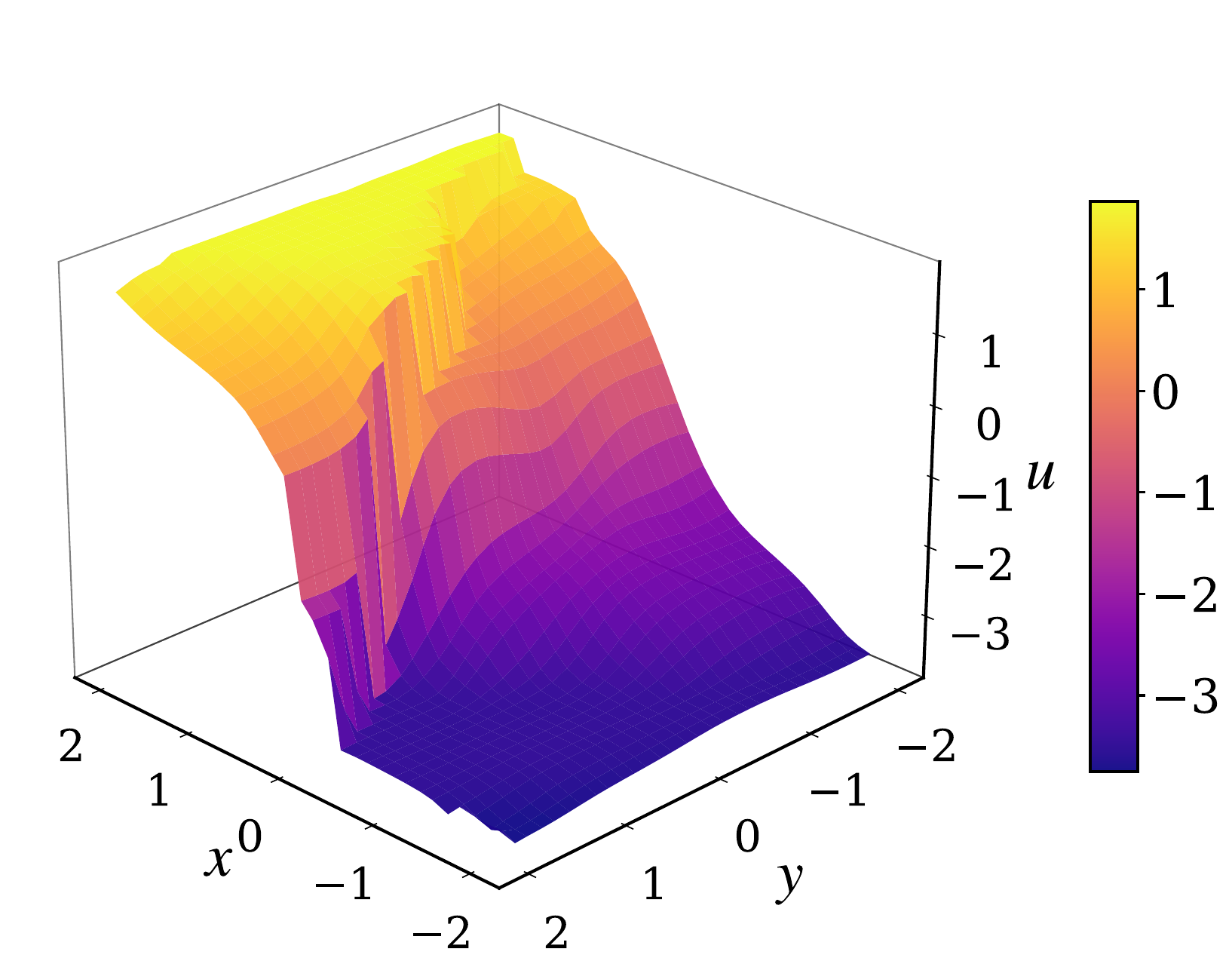} \label{fig:2d_PINN_solution_mu2} }
\subfigure[]{
\includegraphics[width=0.25\linewidth]{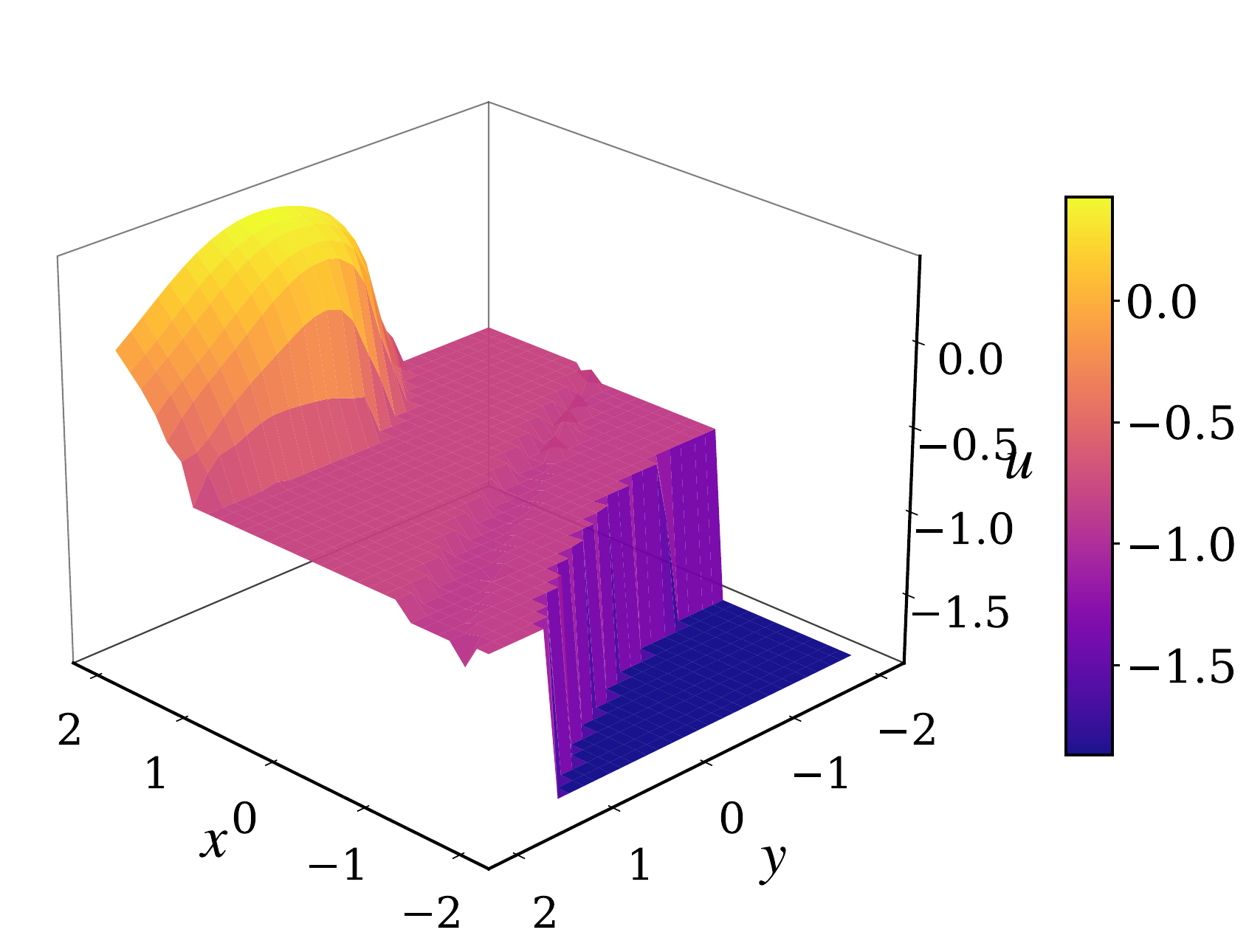} \label{fig:2d_PINN_solution_mu3} }
\subfigure[]{
\includegraphics[width=0.25\linewidth]{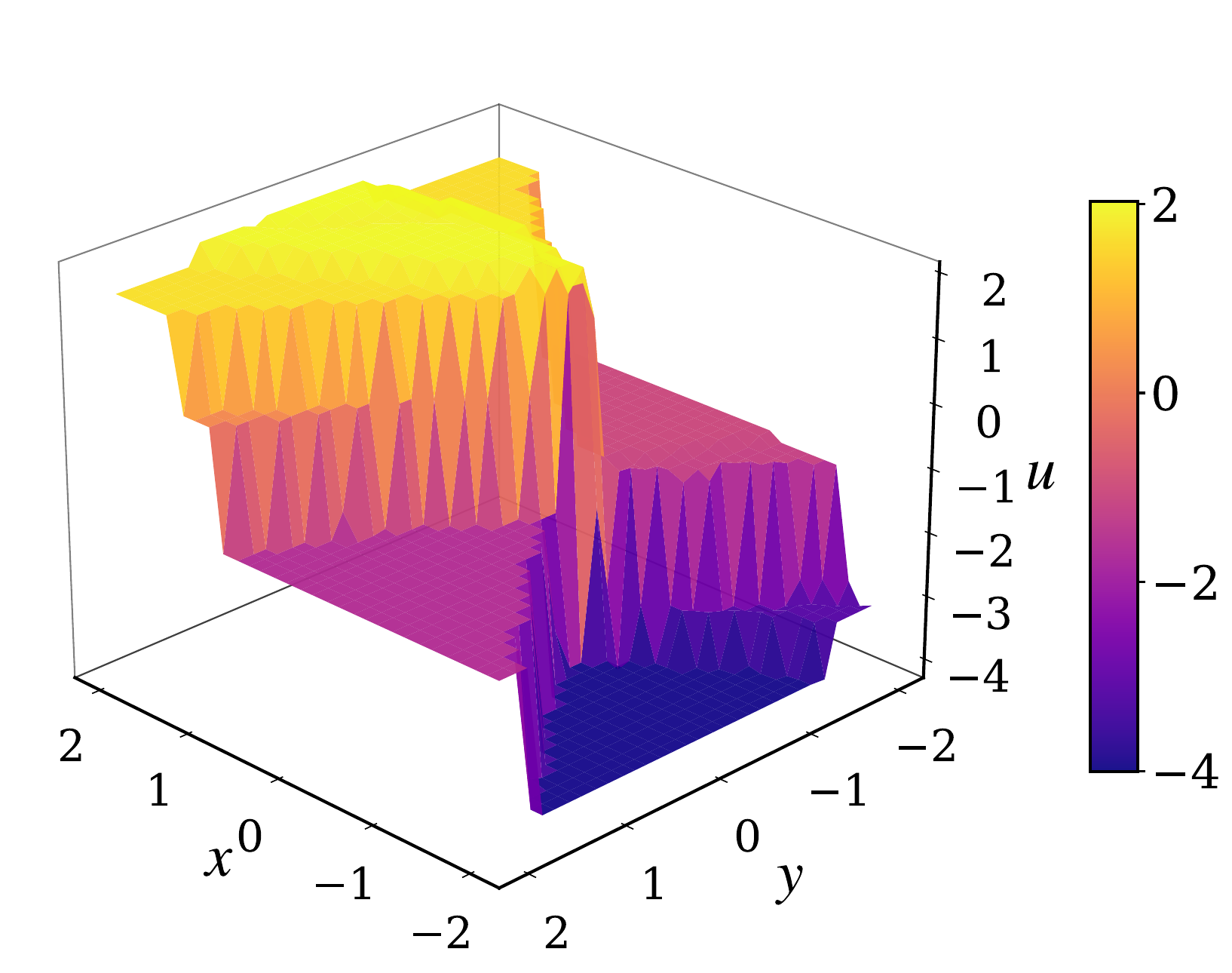} \label{fig:2d_PINN_solution_mu4} }
\vspace{-1ex}
\subfigure[]{
\includegraphics[width=0.25\linewidth]{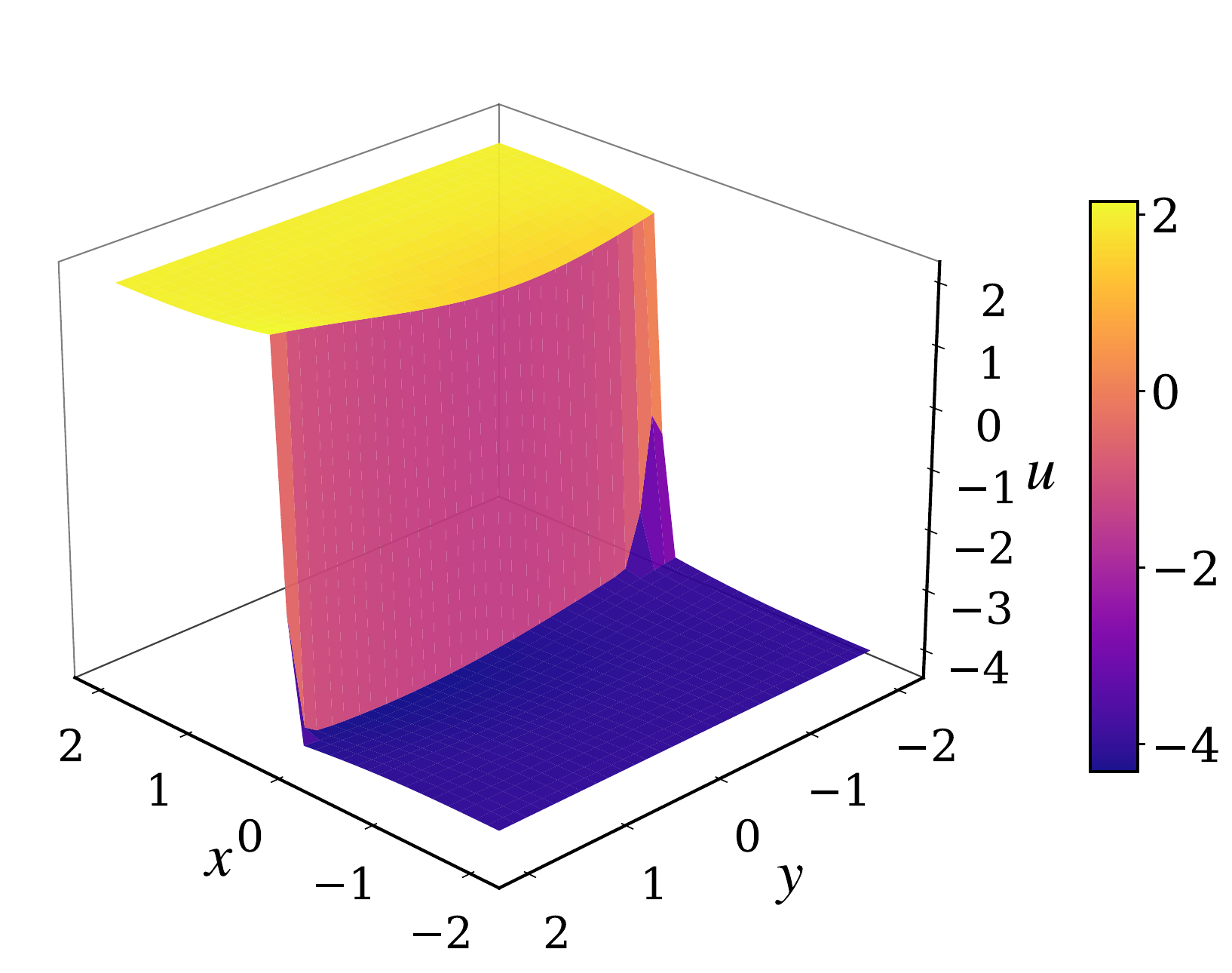} \label{fig:2d_DAE_solution_mu2} }
\subfigure[]{
\includegraphics[width=0.25\linewidth]{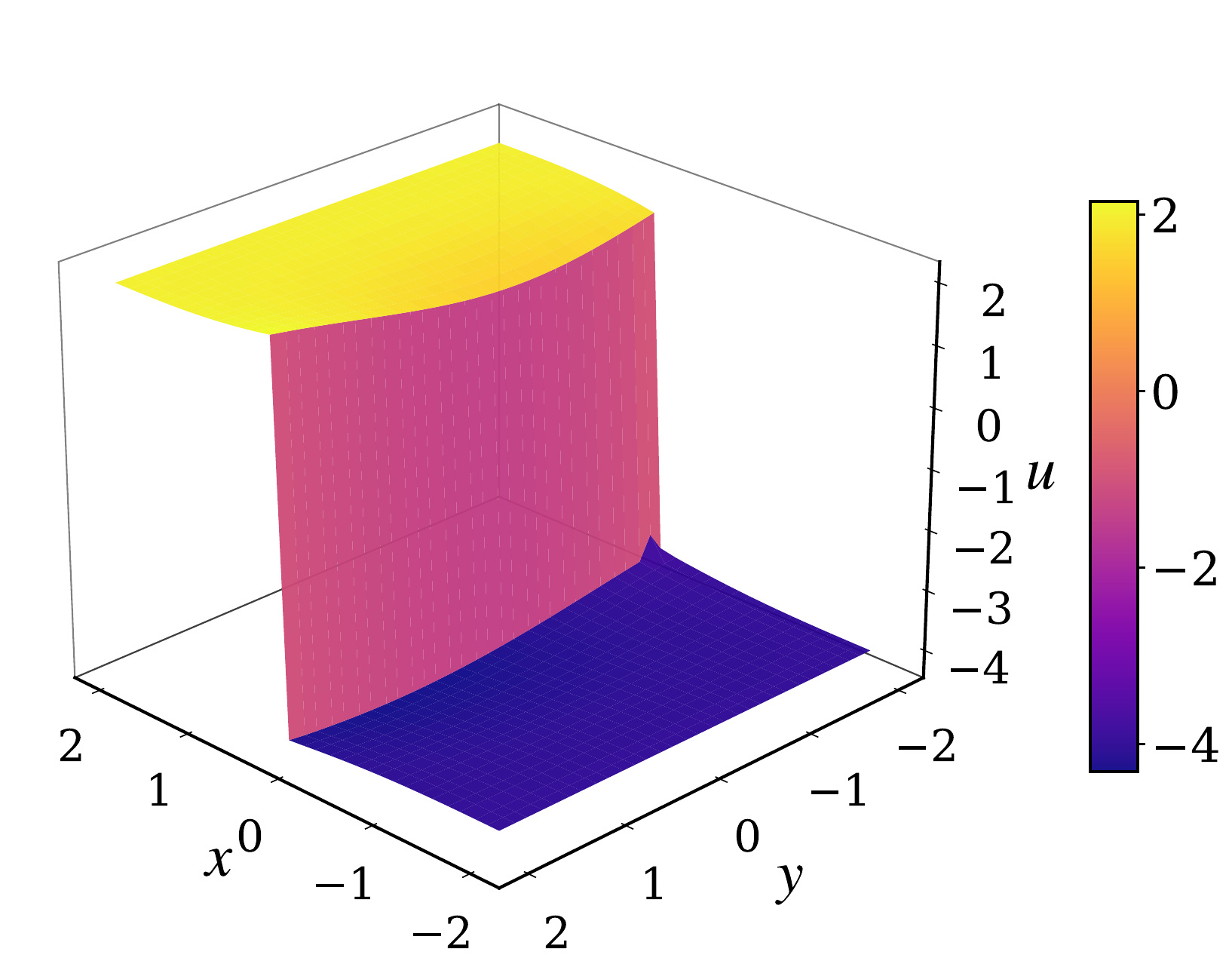} \label{fig:2d_DAE_solution_mu3} }
\subfigure[]{
\includegraphics[width=0.25\linewidth ]{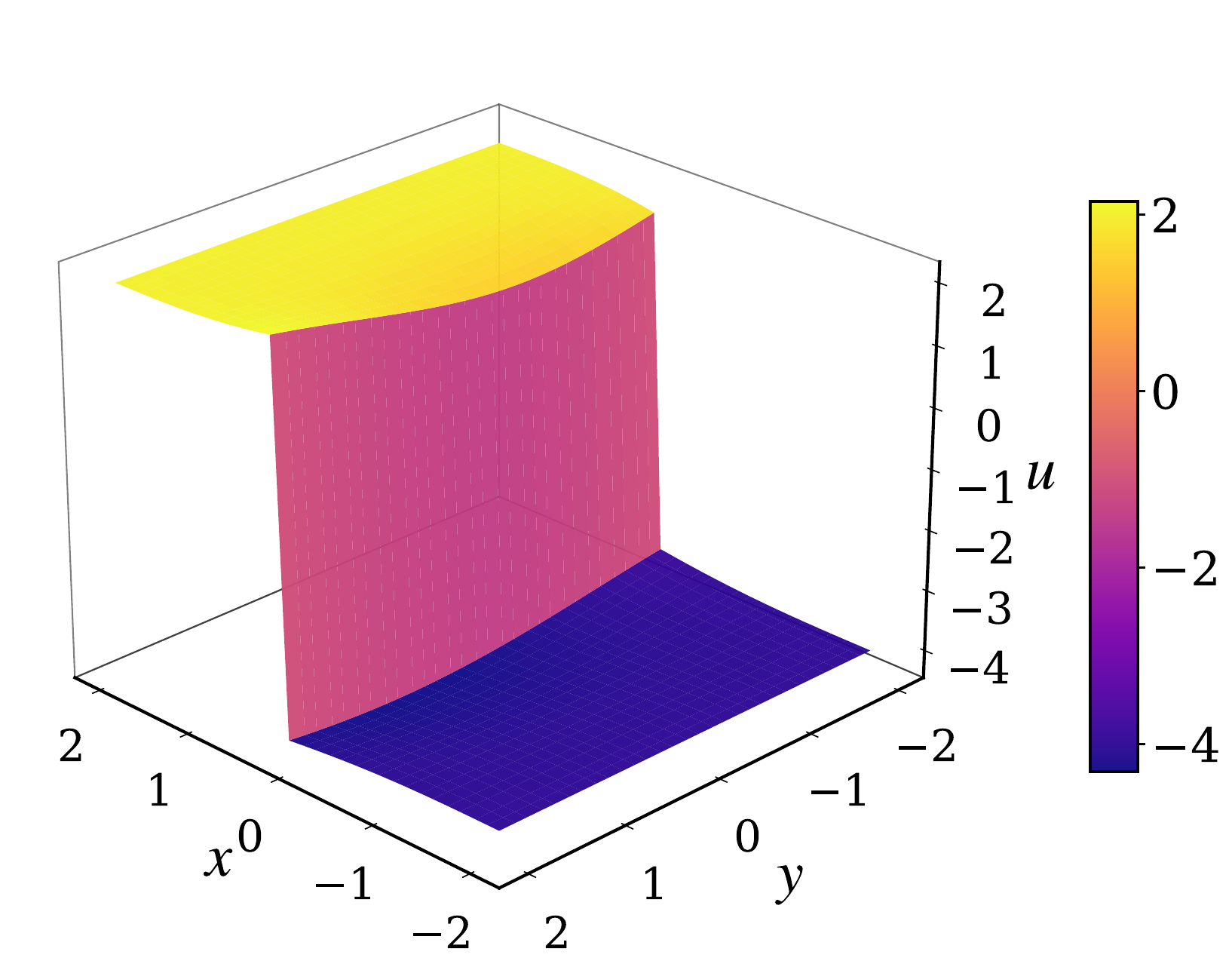} \label{fig:2d_DAE_solution_mu4} }
\caption{2D forward problem: Reference solutions \subref{fig:2d_true_solution_mu2} \subref{fig:2d_true_solution_mu3} \subref{fig:2d_true_solution_mu4} and predicted solutions  obtained by PINN with RAR \subref{fig:2d_PINN_solution_mu2} \subref{fig:2d_PINN_solution_mu3} \subref{fig:2d_PINN_solution_mu4} and DAE with RAR \subref{fig:2d_DAE_solution_mu2}  \subref{fig:2d_DAE_solution_mu3}  \subref{fig:2d_DAE_solution_mu4} for $\mu_2=10^{-2}, 10^{-3}, 10^{-4}$, respectively.}
\label{fig: 2D forward problem solution}
\end{figure}

For DAE/PINN with RAR, initially, 1800/2800 points are randomly selected as residual points. Then 400/400 residual points are adaptively added with $m=10$, $|\mathcal{S}|=30000/35000$, and $\mathcal{E}_0=10^{-6}$. Fig.~\ref{fig: 2D forward problem solution} shows that DAE with RAR better captures transition layers than PINN with RAR. Table~\ref{tab:2D forward problem} shows that DAE yields higher accuracy with less training time, while the PINN approximation remains unreliable even with more data and iterations. Although the DAE has a larger maximum absolute error for small \(\mu_2 = 10^{-3}/10^{-4}\), the proportion of errors exceeding 1 is much smaller (0.03\%/0.02\%) compared to PINN (15\%/30\%). Adaptive sampling (e.g., RAR) can further improve the DAE, while PINN with RAR shows little improvement. Tables \ref{tab:training data_2}, \ref{tab:neurons_2} and \ref{tab:random seeds_2} indicate that the DAE yields more accurate solutions under various parameter settings.

\begin{table}[htbp]
    \centering
    \fontsize{6}{7}\selectfont
    \caption{2D problem: The loss values, $L^2$,  $L^{\infty}$ errors and training times of different methods for various $\mu_2$ values.  Dashes (--) denote values identical to those for $\mu_2 = 10^{-2}$.}
    \begin{tabular}{c|ccccc}
        \toprule
        $\mu_2 $ & Method & $e_{loss}$ & $e_2$ & $e_{\infty}$ & Time (s) \\
        \midrule
        \multicolumn{6}{c}{Standard Methods (Without RAR)} \\
        \midrule
        \multirow{2}{*}{$10^{-2}$}
        & PINN   & 4.14e-01 & 5.21e-01 & 6.77e+00 & 638 \\
        & \textbf{DAE}    & \textbf{1.09e-05} & \textbf{6.41e-03} & \textbf{1.57e+00} & \textbf{213} \\
        \cmidrule(r){1-6}
        \multirow{2}{*}{$10^{-3}$}
        & PINN   & 1.84e+00 & 4.21e-01 & 6.42e+00 & 637 \\
        & \textbf{DAE}    & -- & \textbf{1.61e-02} & \textbf{5.17e+00} & -- \\
        \cmidrule(r){1-6}
        \multirow{2}{*}{$10^{-4}$}
        & PINN   & 4.36e+00 & 3.60e-01 & 5.72e+00 & 641 \\
        & \textbf{DAE}    & -- & \textbf{2.17e-02} & \textbf{5.52e+00} & -- \\
        \midrule
        \multicolumn{6}{c}{Enhanced Methods (With RAR)} \\
        \midrule
        \multirow{2}{*}{$10^{-2}$}
        & PINN+RAR & 1.88e+00 & 3.59e-01 & 5.77e+00 & 33670 \\
        & \textbf{DAE+RAR}  & \textbf{4.56e-08} & \textbf{1.19e-03} & \textbf{1.74e-01} & \textbf{8060} \\
        \cmidrule(r){1-6}
        \multirow{2}{*}{$10^{-3}$}
        & PINN+RAR & 3.07e+00 & 8.45e-01 & 3.68e+00 & 31881 \\
        & \textbf{DAE+RAR}  & -- & \textbf{1.52e-03} & \textbf{4.83e-01} & -- \\
        \cmidrule(r){1-6}
        \multirow{2}{*}{$10^{-4}$}
        & PINN+RAR & 2.34e+00 & 6.78e-01 & 6.25e+00 & 33436 \\
        & \textbf{DAE+RAR}  & -- & \textbf{5.89e-04} & \textbf{2.05e-01} & -- \\
        \bottomrule
    \end{tabular}
    \label{tab:2D forward problem}
\end{table}
\begin{table}[htbp]
\centering
\fontsize{6}{7}\selectfont
\caption{2D problem: Relative $L^2$ errors for $\mu_2 = 10^{-3}$ with varying training data. A 5-layer, 10-neuron network is used with 15000 training iterations.}%
\vspace{0cm} 
\begin{tabular}{c@{\hskip 2pt}c@{\hskip 3pt}c@{\hskip 3pt}c@{\hskip 4pt}c@{\hskip 3pt}c@{\hskip 3pt}c@{\hskip 4pt}c@{\hskip 3pt}c@{\hskip 3pt}c@{\hskip 4pt}c}

\toprule
\multirow{2}{*}{$N_{b}/N_{i}/N_{hp}$} & \multicolumn{9}{c}{$N_f/N_h$} \\ 
\cmidrule(lr){2-11}
& & \multicolumn{3}{c}{1000} & \multicolumn{3}{c}{1500} & \multicolumn{3}{c}{2000} \\
\cmidrule(lr){3-5} \cmidrule(lr){6-8} \cmidrule(lr){9-11} 
&  & PINN &gPINN& DAE & PINN&gPINN & DAE & PINN&gPINN & DAE \\
\midrule
500  &   & 2.94e-01& 5.23e-01  & \textbf{1.54e-02} & 3.70e-01& 4.86e-01 & \textbf{1.84e-02} & 4.27e-01& 4.72e-01 & \textbf{1.61e-02} \\
1000 &   & 4.34e-01& 4.29e-01  & \textbf{1.70e-02} & 3.76e-01 & 5.00e-01& \textbf{1.55e-02} & 4.00e-01 & 4.72e-01& \textbf{1.61e-02} \\
1500 &   & 2.46e-01 & 4.26e-01 & \textbf{1.63e-02} & 3.46e-01 & 4.71e-01& \textbf{1.58e-02} & 4.51e-01 & 4.73e-01& \textbf{1.53e-02} \\
\bottomrule
\end{tabular}
\vspace{0cm} 
\label{tab:training data_2}
\end{table}
\begin{table}[htbp]
\centering
\fontsize{6}{7}\selectfont
\caption{2D problem: Relative $L^2$ errors for $\mu_2 = 10^{-3}$ with different network sizes. Training data fixed at $N_f/N_h = 2000$, $N_b/N_i/N_{hp} = 2000$, with 15000 iterations.}
\vspace{0cm} 
\begin{tabular}{c@{\hskip 2pt}c@{\hskip 3pt}c@{\hskip 3pt}c@{\hskip 4pt}c@{\hskip 3pt}c@{\hskip 3pt}c@{\hskip 4pt}c@{\hskip 3pt}c@{\hskip 3pt}c@{\hskip 4pt}c}
\toprule
\multirow{2}{*}{Layers} & \multicolumn{9}{c}{Neurons} \\ 
\cmidrule(lr){2-11}
& & \multicolumn{3}{c}{8} & \multicolumn{3}{c}{12} & \multicolumn{3}{c}{16} \\
\cmidrule(lr){3-5} \cmidrule(lr){6-8} \cmidrule(lr){9-11} 
&  & PINN&gPINN & DAE & PINN&gPINN & DAE & PINN&gPINN & DAE \\
\midrule
3  &   & 3.22e-01 & 4.49e-01 & \textbf{3.19e-02} & 3.44e-01& 4.74e-01 & \textbf{3.02e-02} & 3.62e-01& 4.46e-01 & \textbf{5.35e-02} \\
4  &   & 3.34e-01 & 4.55e-01 & \textbf{2.43e-02} & 3.18e-01& 4.48e-01 & \textbf{2.98e-02} & 2.93e-01 & 4.65e-01& \textbf{4.35e-02} \\
5  &   & 4.79e-01 & 4.43e-01 & \textbf{2.66e-02} & 4.17e-01& 4.72e-01 & \textbf{4.64e-02} & 4.38e-01& 4.49e-01 & \textbf{2.84e-02} \\
\bottomrule
\end{tabular}
\vspace{0cm} 
\label{tab:neurons_2}
\end{table}
\begin{table}[htbp]
	\centering
	\fontsize{6}{7}\selectfont    
	\caption{2D forward problem: Relative $L^2$ errors for $\mu_2 = 10^{-3}$ with different random seeds.}
	\begin{tabular}{lccccc}
		\toprule
		\diagbox [width=8em,trim=l] {Methods}{{Seed}} & 33 & 99 & 202 & 1000 & 9999 \\ \midrule
		PINN & 3.09e-01 & 3.82e-01 & 3.91e-01 & 4.41e-01 & 3.80e-01\\
		gPINN & 4.42e-01 & 4.73e-01 & 4.36e-01 & 4.39e-01 & 4.61e-01\\
        DAE & \textbf{2.81e-02} & \textbf{3.37e-02} & \textbf{6.82e-02} & \textbf{7.45e-02} & \textbf{1.82e-02}\\ \bottomrule
	\end{tabular}
	\vspace{0cm}
	\label{tab:random seeds_2}
\end{table}

\subsection{3D problem}\label{3D}
Consider the following 3D problem with $\mu_3=0.01$ and the source $f(x,y,z)=\cos{ ( \pi x )}\cos{( \pi y )} \cos{( \pi z )}$:
\begin{align} \label{forwardexample3}
\begin{dcases}
\displaystyle \mu_3 \Delta u - \partial_t u = -u \left( \partial_x u + \partial_y u+ \partial_z u \right) +f(x,y,z), \\
\displaystyle u(-1,y,z,t)=-4, ~ u(1,y,z,t)=2, \\ u(x,y,z,t)=u(x,y+2,z,t)=u(x,y,z+2,t), \\
\displaystyle u(x,y,z,0)= u_0(x,y,z,\mu_3 ), \quad \\
x \in [-1, 1], y \in [-1, 1], z \in [-1, 1], t \in [0,0.5].
\end{dcases}
\end{align}
Using the asymptotic method, we reduce the problem to two problems that determine the outer functions:
\begin{align*}
\begin{cases}
\displaystyle \varphi^{(\pm)} \left( \partial_x \varphi^{(\pm)} + \partial_y \varphi^{(\pm)} + \partial_z \varphi^{(\pm)} \right) = \cos{ ( \pi x )} \cos{( \pi y )} \cos{( \pi z )}, \\
\varphi^{(\pm)}(\pm1,y,z)=-1\pm3, \ \varphi^{(\pm)}(x,y,z)=\varphi^{(\pm)}(x,y+2,z)=\varphi^{(\pm)}(x,y,z+2).
\end{cases}
\end{align*}
Using the method of characteristics, we find the zeroth-order outer asymptotic functions:
\begin{align*}
\varphi^{(-)} &= -\sqrt{16+2\int_{-1}^x \cos{ ( \pi s )}\cos{( \pi (s+y-x) )} \cos{( \pi (s+z-x) )} d s } ,\\
\varphi^{(+)} &= \sqrt{4-2\int_x^1 \cos{ ( \pi s )}\cos{( \pi (s+y-x) )} \cos{( \pi (s+z-x) )} d s}.
\end{align*}
The leading term of the asymptotics of the front $h_0(y,z,t)$ satisfies
\begin{align}\label{eq44}
\begin{cases}
\displaystyle \partial_t{h_0} =\tfrac{1}{2} \left( \partial_y{h_0}+\partial_z{h_0}-1 \right) \left( \varphi^{(-)} (h_0,y,z)+\varphi^{(+)} (h_0,y,z)\right),\\
h_0(y,z,0)=h_0^{*}=0, \quad h_0(y,z,t)=h_0(y+2,z,t)=h_0(y,z+2,t).
\end{cases}
\end{align}

\begin{figure}[htbp]
\centering
\subfigure[Predicted $h_0$]{
  \includegraphics[width=0.25\linewidth]{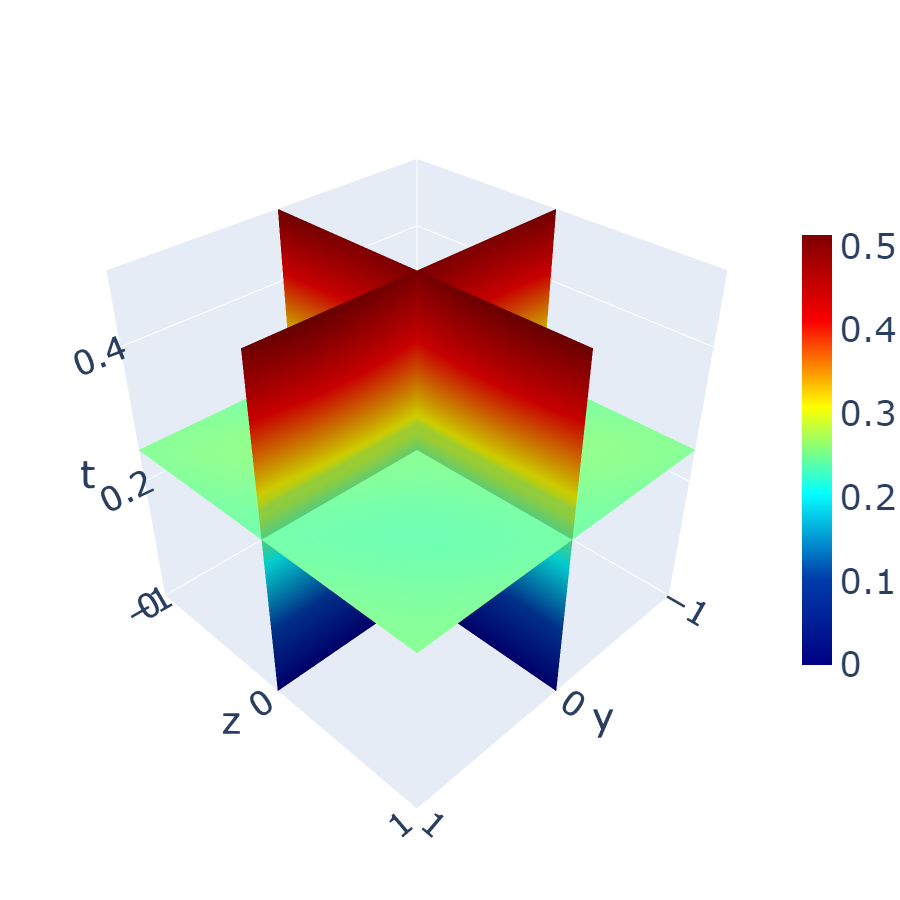}\label{fig:x0example3}
}
\subfigure[Reference solution]{
  \includegraphics[width=0.25\linewidth]{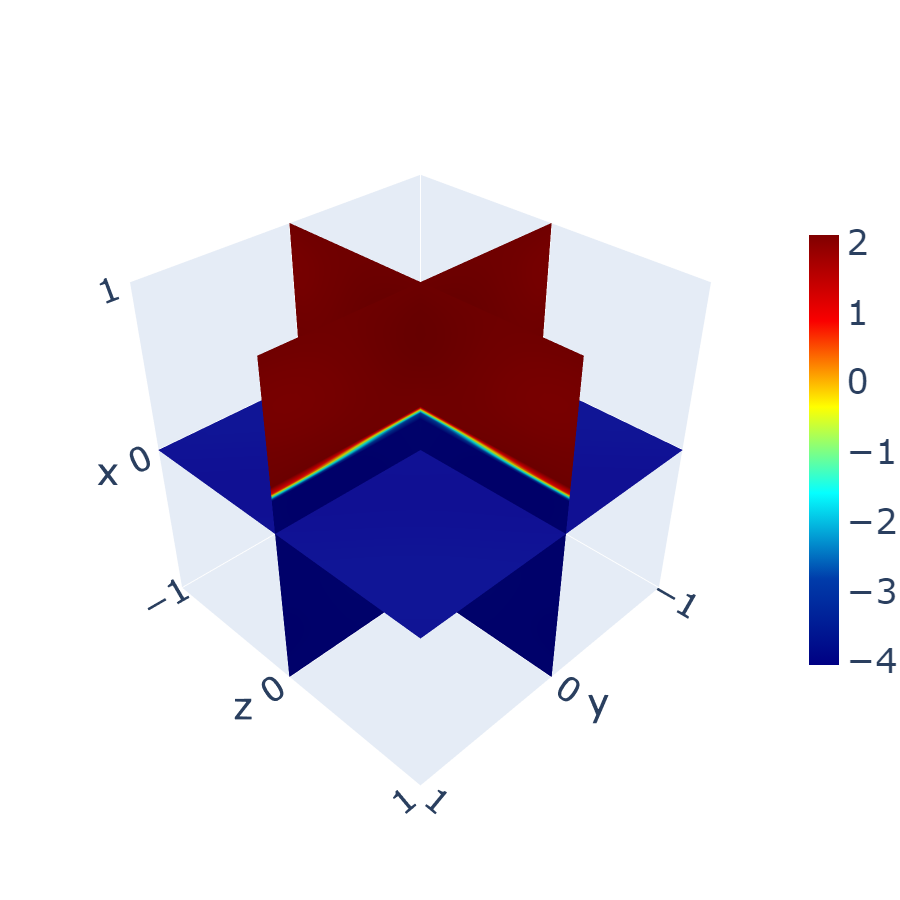}\label{fig:3d_true_solution_mu2}
}
\vspace{1ex} 

\subfigure[Prediction of PINN]{
  \includegraphics[width=0.22\linewidth]{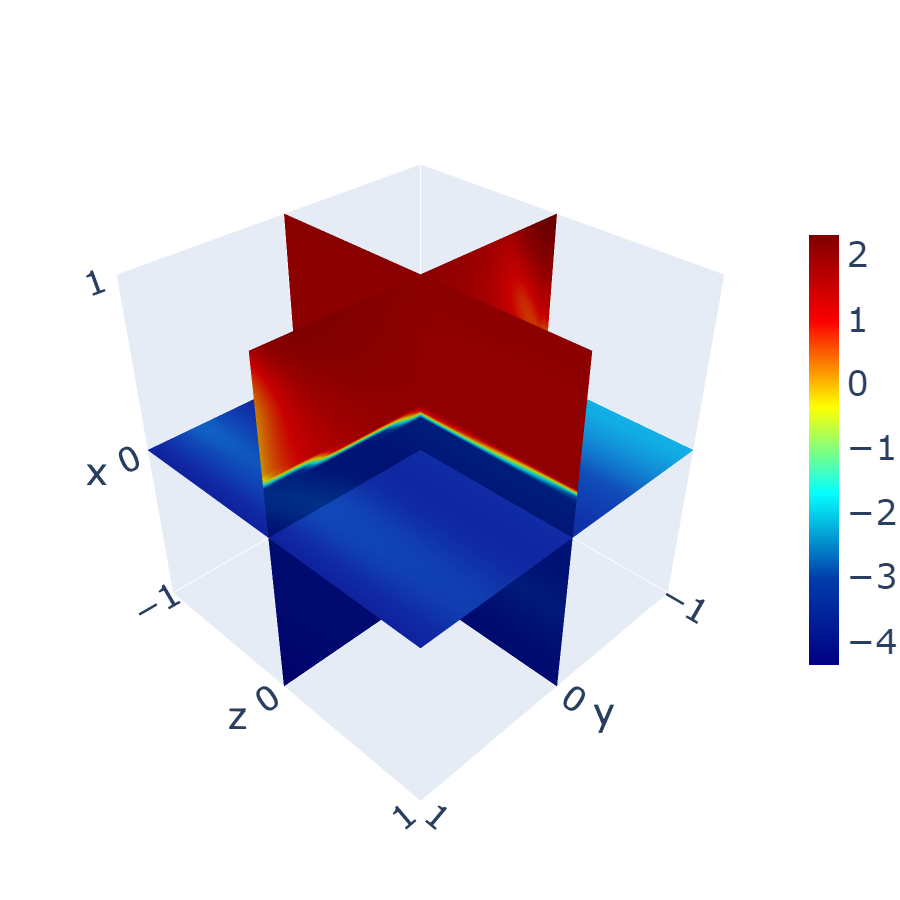}\label{fig:3d_PINN_solution_mu2}
}
\subfigure[Prediction of PINN with RAR]{
  \includegraphics[width=0.22\linewidth]{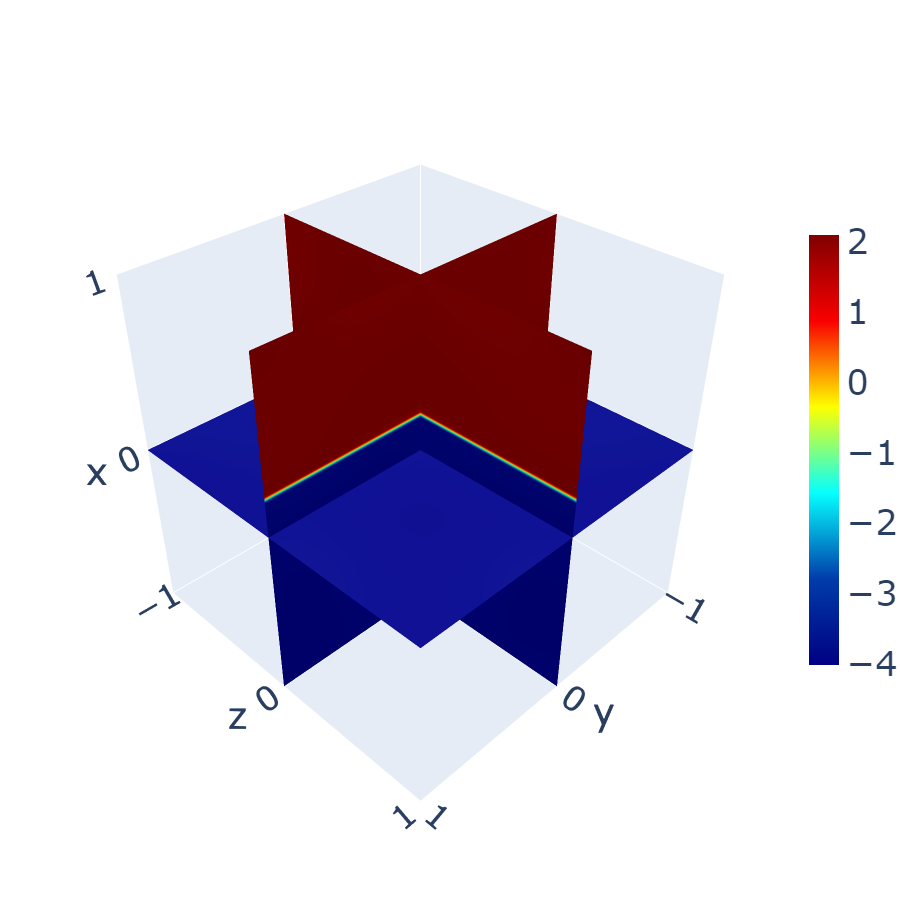}\label{fig:3d_PINN_RAR_solution_mu2}
}
\subfigure[Prediction of gPINN]{
  \includegraphics[width=0.22\linewidth]{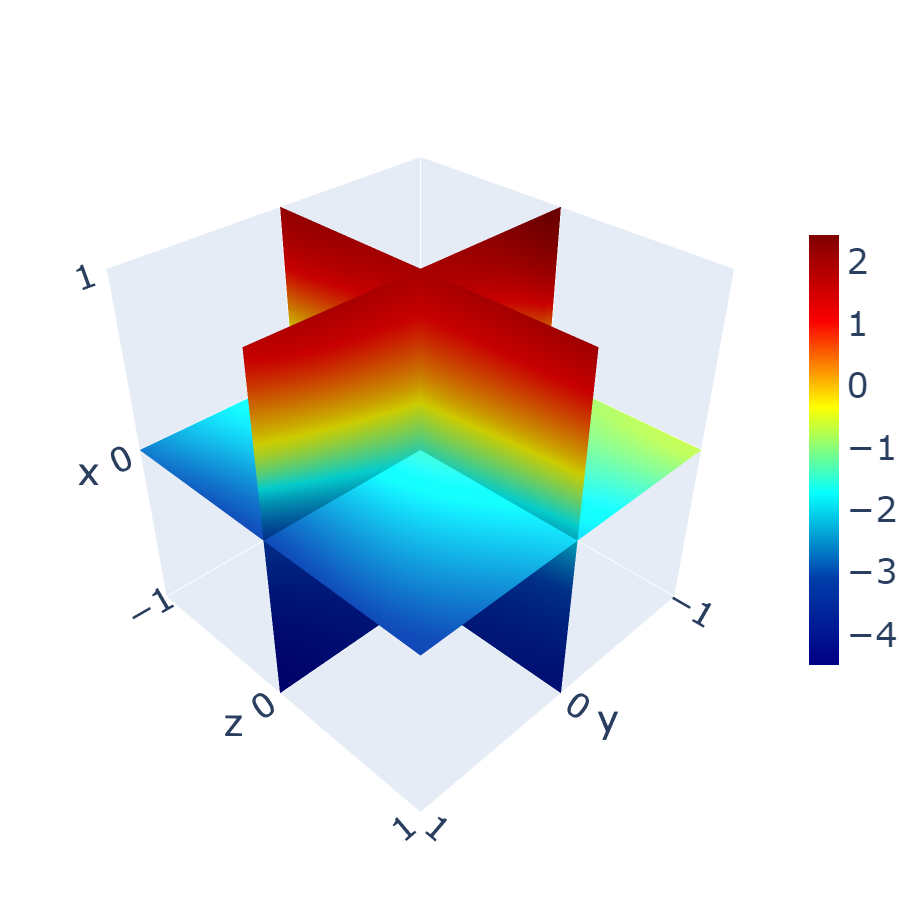}\label{fig:3d_gPINN_solution_mu2}
}
\subfigure[Prediction of DAE with RAR]{
  \includegraphics[width=0.22\linewidth]{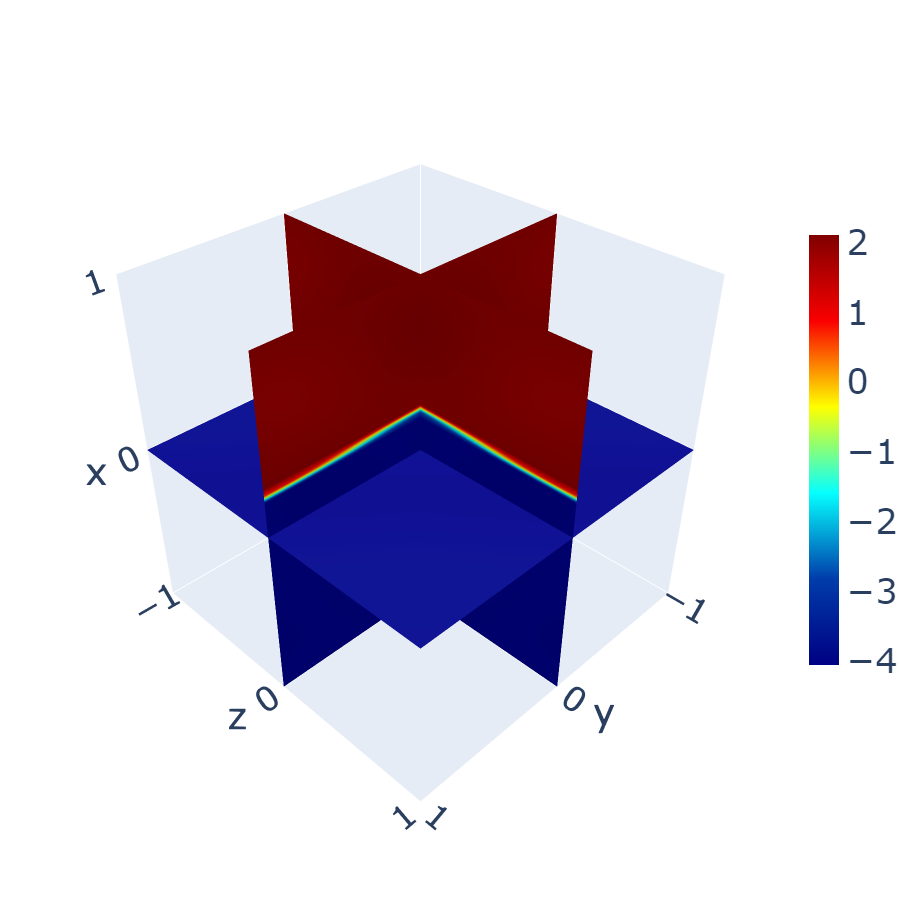}\label{fig:3d_DAE_solution_mu2}
}
\vspace{1ex} 

\subfigure[Error of PINN]{
  \includegraphics[width=0.22\linewidth]{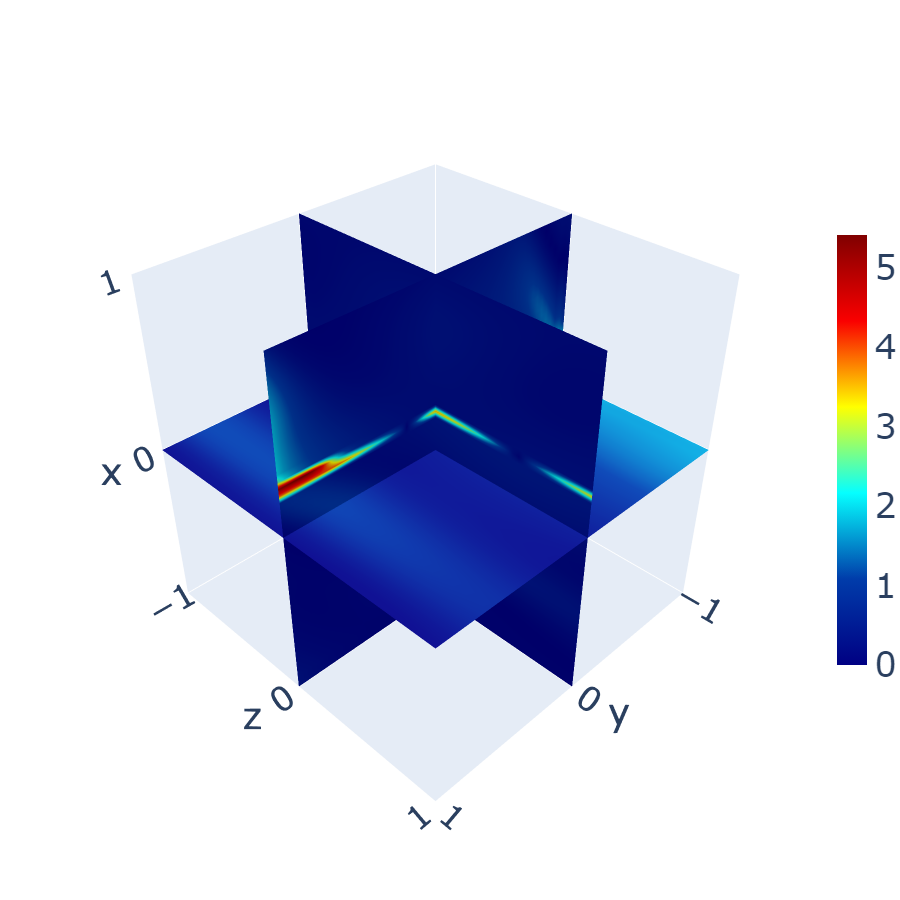}\label{fig:3d_PINN_error_mu2}
}
\subfigure[Error of PINN with RAR]{
  \includegraphics[width=0.22\linewidth]{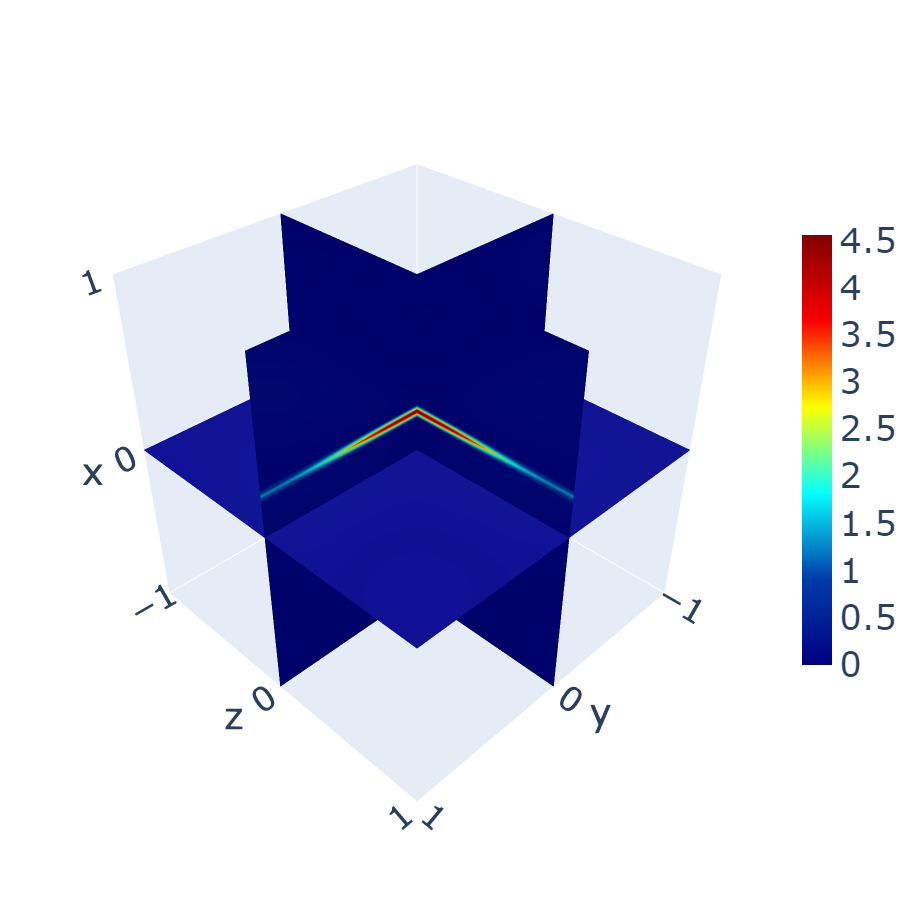}\label{fig:3d_PINN_RAR_error_mu2}
}
\subfigure[Error of gPINN]{
  \includegraphics[width=0.22\linewidth]{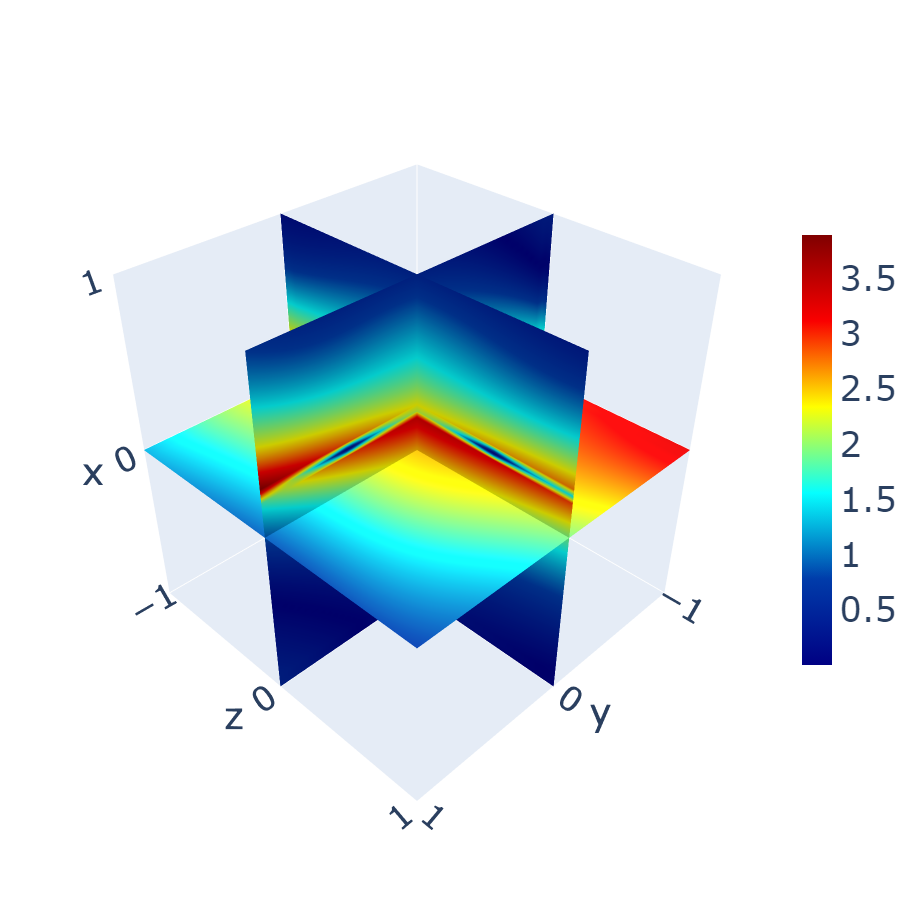}\label{fig:3d_gPINN_error_mu2}
}
\subfigure[Error of DAE with RAR]{
  \includegraphics[width=0.22\linewidth]{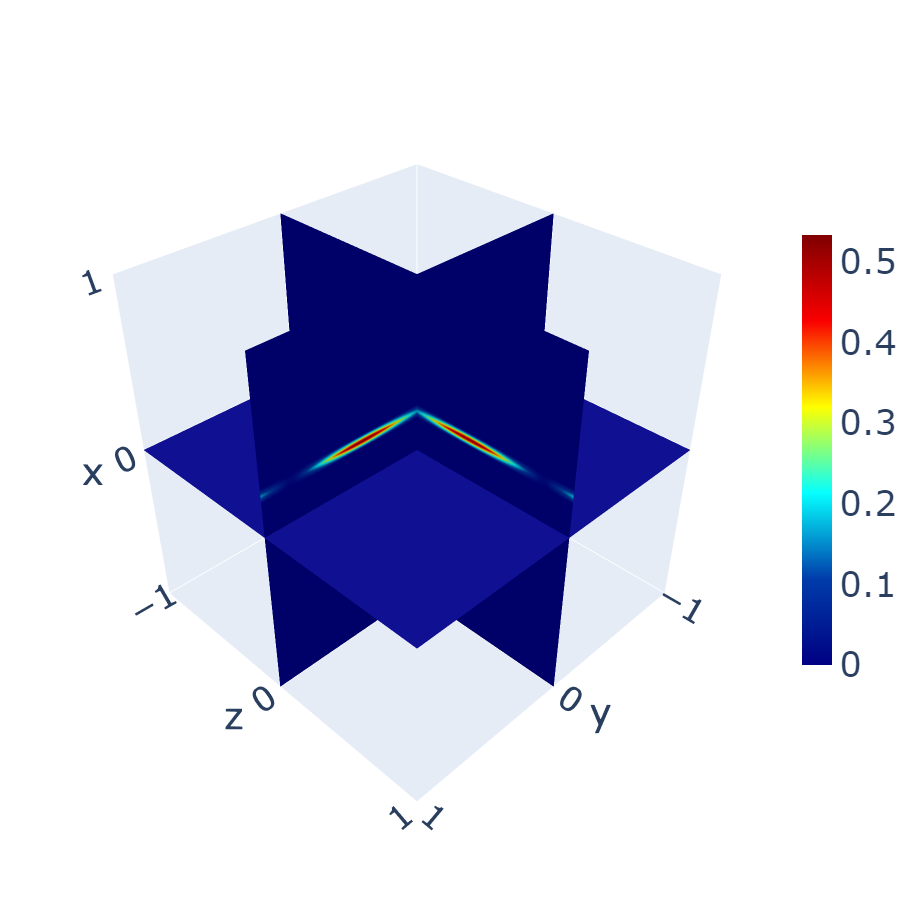}\label{fig:3d_DAE_error_mu2}
}

\caption{3D forward problem: Approximate solution of $h_0$ \subref{fig:x0example3}. Reference solution \subref{fig:3d_true_solution_mu2} and predicted solutions at $t=0.2$ of PINN \subref{fig:3d_PINN_solution_mu2}, PINN with RAR \subref{fig:3d_PINN_RAR_solution_mu2}, gPINN \subref{fig:3d_gPINN_solution_mu2} and DAE with RAR \subref{fig:3d_DAE_solution_mu2}. Absolute errors at $t=0.2$ of PINN \subref{fig:3d_PINN_error_mu2}, PINN with RAR \subref{fig:3d_PINN_RAR_error_mu2}, gPINN \subref{fig:3d_gPINN_error_mu2} and DAE with RAR \subref{fig:3d_DAE_error_mu2}.}
\label{fig: 3D forward problem solution}
\end{figure}

In view of the approximate solution of \eqref{eq44} obtained by the DAE, the transition layer is located in the region $-1 \leq h_0 (y,z,t) \leq 1$ for any $y\in [-1,1], z \in [-1,1]$, and $t \in [0,0.5]$ (Fig.~\ref{fig:x0example3}). Thus Assumption~\ref{assumption}(c) holds. The initial condition ${u}_0(x,y,z,\mu_3 )=3\tanh(\frac{x}{\mu_3}+y+z)-1$ satisfies Assumption \ref{assumption}(d). Thus, the problem has a solution in the form of an autowave with a transitional moving layer localized near $h_0(y,z,t)$, whose zeroth-order approximation is given by
\begin{align*} \label{asymptoticsolEXAMPLE1}
U_{0}={\footnotesize\begin{cases}
 \displaystyle \varphi^{(-)}(x,y,z) +\frac{\varphi^{(+)}(h_0,y,z)-\varphi^{(-)}(h_0,y,z)}{\exp \Big( \frac{ \left(h_0-x \right) (\varphi^{(+)}(h_0,y,z)-\varphi^{(-)}(h_0,y,z) ) (1-\partial_y{h_0}-\partial_z{h_0} )}{2 \mu_3} \Big)+1} , ~x \in [-1,h_0],\\
\displaystyle \varphi^{(+)}(x,y,z)+\frac{\varphi^{(-)}(h_0,y,z)-\varphi^{(+)}(h_0,y,z)}{\exp \Big( \frac{ (h_0-x ) (\varphi^{(-)}(h_0,y,z)-\varphi^{(+)}(h_0,y,z) ) (1-\partial_y{h_0}-\partial_z{h_0} )}{2 \mu_3} \Big)+1} , ~x \in [h_0,1].
\end{cases}}
\end{align*}

For DAE/PINN with RAR, initially 3600/5600  points are randomly selected as residual points, and 400/400 residual points are adaptively added with $m=20$, $|\mathcal{S}|=50000/50000$, and $\mathcal{E}_0=10^{-6}$. A 32-point Gauss-Legendre quadrature is used. Fig.~\ref{fig: 3D forward problem solution} shows that PINN yields large errors, and gPINN offers no improvement. While RAR improves PINN prediction to better match the reference solution, the errors remain relatively large. In contrast, the DAE with RAR yields lower errors. Table~\ref{tab:3D forward problem} shows that the DAE outperforms the other methods with lower maximum error and reduced training time, and the DAE with RAR can further improve the performance. Table \ref{tab:random seeds_3} shows that the accuracy of the DAE is unaffected by the choice of random seeds, indicating its robustness.

\begin{table}[htbp]
    \centering
    \fontsize{6}{7}\selectfont
    \caption{3D  problem: The loss values, $L^2$,  $L^{\infty}$ errors and training times of different methods.}
    \begin{tabular}{ccccc}
        \toprule
         Method & $e_{loss}$ & $e_2$ & $e_{\infty}$ & Time (s) \\
        \midrule

         PINN   & 1.02e+00 & 3.15e-01 & 5.95e+00 & 1842 \\
         PINN+RAR & 2.60e-02 & 1.55e-01 & 5.99e+00 & 30904 \\
         gPINN & 4.61e+00 & 4.80e-01 & 5.70e+00 & 4707 \\      
        \textbf{DAE}  & \textbf{4.04e-06} & \textbf{1.22e-02} & \textbf{1.93e+00} & \textbf{587} \\
        \textbf{DAE+RAR} & \textbf{6.77e-08} & \textbf{9.52e-03} & \textbf{8.92e-01} & \textbf{6914} \\
        \bottomrule
    \end{tabular}
    \label{tab:3D forward problem}
\end{table}

\begin{table}[htbp]
	\centering
	\fontsize{6}{7}\selectfont    
	\caption{3D forward problem: Relative $L^2$ errors for different random seeds.}
	\begin{tabular}{lccccc}
		\toprule
		\diagbox [width=8em,trim=l] {Methods}{{Seed}} & 33 & 99 & 202 & 5678 & 9999 \\ \midrule
		PINN & 4.43e-01 & 5.51e-01 & 3.37e-01 & 3.34e-01 & 4.06e-01\\
		gPINN & 4.96e-01 & 4.93e-01 & 4.90e-01 & 4.89e-01 & 4.81e-01\\
        DAE & \textbf{1.16e-02} & \textbf{1.13e-02} & \textbf{1.59e-02} & \textbf{1.04e-02} & \textbf{1.52e-02}\\ \bottomrule
	\end{tabular}
	\vspace{0cm}
	\label{tab:random seeds_3}
\end{table}

\section{Conclusions}
\label{sec5}

In this paper, we have developed the DAE method which integrates asymptotic analysis with deep learning techniques to efficiently solve singularly perturbed time-dependent reaction-advection-diffusion equations exhibiting internal transition layers. We first derive the governing equations for transition layers in a general $d$-dimensional setting, and then solve the resulting equations using the standard PINN. The effectiveness and robustness of the approach were shown on several numerical experiments. The DAE consistently outperforms traditional PINN methods, achieving significantly higher accuracy and faster convergence across a variety of test cases. Furthermore, by integrating DAE with sampling strategy RAR, we observed notable performance improvements.


\appendix
\section{First-order deep asymptotic expansion method}\label{appendix:DAE1} Here, we provide the derivation of the first-order asymptotic solution for problem \eqref{forwardexample1}. 

The general problems defining the first-order outer approximations $\bar{u}_{1}^{(\mp)}(\boldsymbol{x})$ reads
\begin{equation} \label{firstorderregularequation}
\begin{dcases}
 \sum_{i=1}^{d} \frac{\partial\bar{u}_{1}^{(\mp)}}{\partial x_i} +\bar{u}_{1}^{(\mp)} K^{(\mp)}(\boldsymbol{x})=W^{(\mp)}(\boldsymbol{x}),\\  
\bar{u}_{1}^{(-)}(-a,\boldsymbol{x}^* )=0, ~\bar{u}_{1}^{(+)}(a,\boldsymbol{x}^*)=0,~\bar{u}_{0}^{(\mp)}(\boldsymbol{x}+\boldsymbol{P}) = \bar{u}_{0}^{(\mp)}(\boldsymbol{x}),
\end{dcases}
\end{equation}
with 
$K^{(\mp)}(\boldsymbol{x})= \frac{1}{A(\varphi^{(\mp)},\boldsymbol{x})} \frac{\partial A(\varphi^{(\mp)},\boldsymbol{x})}{\partial \varphi^{(\mp)}} \sum_{i=1}^{d}  \frac{\partial\varphi^{(\mp)}}{\partial x_i}$, and $W^{(\mp)}(\boldsymbol{x})= \frac{\Delta \varphi^{(\mp)}}{A(\varphi^{(\mp)},\boldsymbol{x})}.$ The solutions $\bar{u}_{1}^{(\mp)}(\boldsymbol{x})$ can be explicitly obtained for any dimension $d$ from \eqref{firstorderregularequation}.

For problem (\ref{forwardexample1}), we derive the following equations from \eqref{firstorderregularequation}:
\begin{align*}
\begin{split}
&\displaystyle  \varphi^{(\mp)}(x)(\bar{u}_{1}^{(\mp)}(x))'+\bar{u}_{1}^{(\mp)}(x)(\varphi^{(\mp)}(x))'=-  (\varphi^{(\mp)}(x))'', \\
&\bar{u}_{1}^{(-)}(0)=0,\quad \bar{u}_{1}^{(+)}(1)=0.
\end{split}
\end{align*}

The solutions to these problems are given by
\begin{align}\label{baru}
\begin{split}
\bar{u}_{1}^{(-)}(x) = \exp \left( \int_{0}^{x} -W(s) ds \right) \int_{0}^{x}- \exp \left( \int_{0}^{s'} W(s) ds \right) Y(s') ds', \\
\bar{u}_{1}^{(+)}(x) = \exp \left( \int_{x}^{1} W(s) ds \right) \int_{x}^{1} \exp \left( \int_{s'}^{1} -W(s) ds \right) Y(s') ds',
\end{split}
\end{align}
with $W(x)=\displaystyle \frac{1}{\varphi^{(\mp)}(x)}\frac{d \varphi^{(\mp)}(x) }{d x}$ and $Y(x)= \displaystyle \frac{1}{\varphi^{(\mp)}(x)}\frac{d^2 \varphi^{(\mp)}(x) }{d x^2} $.

We also obtain the equations for the first-order transition-layer functions:
\begin{align*}
\begin{split}
 &\frac{\partial^{2} Q_{1}^{(\mp)} }{\partial \xi^2}+ \left((\varphi^{(\mp)}(h_0(t))+Q_{0}^{(\mp)} ) +\frac{dh_0(t)}{dt} \right)\frac{\partial Q_{1}^{(\mp)} }{\partial \xi}  +    Q_{1}^{(\mp)}  \Upsilon^{(\mp)}(\xi,t) \\
 &=   \left(  -(h_1(t) \frac{d \varphi^{(\mp)} }{d x}(h_0(t))  + \bar{u}_{1}^{(\mp)}(h_0(t))) - \frac{dh_1(t)}{dt} \right) \Upsilon^{(\mp)} + r_{1}^{(\mp)}(\xi,t) :=H_1^{(\mp)} (\xi,t),
\end{split}
\end{align*}
with
$ \Upsilon^{(-)}(\xi,t) = \frac{\partial \tilde{u}}{\partial \xi} (\xi,t), \ \xi \leq 0,~ \Upsilon^{(+)}(\xi,t) = \frac{\partial \tilde{u}}{\partial \xi} (\xi,t), \ \xi \geq 0,
$ and
\begin{align*}
r_{1}^{(\mp)}(\xi,t) = -\left( \left(   \xi  \frac{d \varphi^{(\mp)} }{d x}(h_0(t))  \right)\Upsilon^{(\mp)}+Q_{0}^{(\mp)} \frac{d   \varphi^{(\mp)} }{ d x}(h_0(t)) \right)+\frac{\partial Q_{0}^{(\mp)}}{\partial t}.
\end{align*}

We derive additional conditions for $Q_{1}^{(\mp)} (\xi,t)$:
\begin{align*}
Q_{1}^{(\mp)} (0,t)=-\bar{u}_{1}^{(\mp)}(h_0(t))- h_1(t) \frac{d \varphi^{(\mp)} }{d x}(h_0(t))\equiv p_{1}^{(\mp)} (t), \\
Q_{1}^{(-)} ( \xi ,t) \rightarrow 0 \ \text{for} \  \xi \rightarrow -\infty, \quad Q_{1}^{(+)} (\xi,t) \rightarrow 0 \ \text{for} \ \xi \rightarrow +\infty.
\end{align*}
The functions $Q_{1}^{(\mp)}(\xi,t)$ are given by
\begin{equation} \label{Q1function}
Q_{1}^{(\mp)} (\xi,t)=z^{(\mp)}(\xi,t) \left( p_{1}^{(\mp)} (t)   - \int_{0}^{\xi} \frac{1}{z^{(\mp)}(s,t)} \int_{s}^{\mp \infty} H_{1}^{(\mp)} (\eta,t) d\eta ds \right),
\end{equation}
with $ \displaystyle z^{(\mp)}(\xi,t)= \left( \Upsilon^{(\mp)}(0,t) \right)^{-1} \Upsilon^{(\mp)}(\xi,t)$, 
$p_{1}^{(\mp)} (t)=-\bar{u}_{1}^{(\mp)}(h_0(t))- h_1(t) \frac{d \varphi^{(\mp)} }{d x}(h_0(t))$, and  $
    H_1^{(\mp)} (\xi,t)=\left(  -\Big(h_1(t) \frac{d \varphi^{(\mp)} }{d x}(h_0(t))  + \bar{u}_{1}^{(\mp)}(h_0(t))\Big) - \frac{d h_1(t)}{dt} \right) \Upsilon^{(\mp)} + r_{1}^{(\mp)}(\xi,t)$.

From the first-order $C^1$-matching condition
\begin{equation*} \label{matchingfirstord}
\frac{\partial Q_{1}^{(-)}}{\partial \xi }(0,t) + \frac{d \varphi^{(-)}}{d x} (h_0 (t))=\frac{\partial Q_{1}^{(+)}}{\partial \xi }(0,t) + \frac{d \varphi^{(+)}}{d x} (h_0 (t)), \quad t \in \bar{\mathcal{T}},
\end{equation*}
we obtain the equation that determines $h_1(t)$:
\begin{equation} \label{eqforx1}
(\varphi^{(-)} (h_0(t))-\varphi^{(+)} (h_0(t)))\frac{d h_1 (t)}{d t}+\Phi_1(t) h_1(t)=\Phi_2(t),
\end{equation}
with initial condition $h_1(0)=0$, and
\begin{align*}
\Phi_1(t)&=  \frac{k}{2}\left(\varphi^{(+)}(h_0(t))-\varphi^{(-)}(h_0(t))  \right) \left( \frac{d \varphi^{(+)} }{d x}(h_0(t))+ \frac{d \varphi^{(-)} }{d x}(x_0(t)) \right),\\
 \Phi_2(t)&=  \frac{k}{2} \left(  \bar{u}_{1}^{(-)}(h_0(t))   + \bar{u}_{1}^{(+)}(h_0(t))  \right) \left(\varphi^{(-)}(h_0(t)) - \varphi^{(+)}(h_0(t)) \right)  -\frac{d \varphi^{(-)} }{d x} (h_0(t)) \\
 & \quad +\frac{d \varphi^{(+)} }{d x} (h_0(t))
  + \int_{0}^{- \infty} r_{1}^{(-)} (\eta,t) d\eta -\int_{0}^{+ \infty} r_{1}^{(+)} (\eta,t) d\eta.
\end{align*}

The first-order asymptotic solution have the following form:
{\small
\begin{align} \label{asymptoticsolU1}
U_{1}(x,t)=\begin{cases}
\displaystyle \varphi^{(-)} (x) +\frac{ \left(\varphi^{(+)}(h_{0}(t))-\varphi^{(-)}(h_{0}(t)) \right)}{  \exp \left( \left(x-h_0(t) \right) \left( \frac{\varphi^{(-)}(h_{0}(t))-\varphi^{(+)}(h_{0}(t)) }{2\mu_1} \right) \right)+ 1} \\ \qquad \quad +\mu_1 \left(\bar{u}_1^{(-)}(x)+Q_1^{(-)}(\xi,t) \right), \quad x \in [0,h_0(t)+\mu_1 h_1(t)],\\
\displaystyle \varphi^{(+)} (x)+\frac{ -\left(\varphi^{(+)}(h_{0}(t))-\varphi^{(-)}(h_{0}(t))\right)}{  \exp \left(\left(x-h_0(t) \right)  \left( \frac{\varphi^{(+)}(h_{0}(t))-\varphi^{(-)}(h_{0}(t)) }{2\mu_1}\right) \right)+ 1}  \\  \qquad \quad +\mu_1 \left(\bar{u}_1^{(+)}(x)+Q_1^{(+)}(\xi,t)\right) , \quad x \in [h_0(t)+\mu_1 h_1(t),1].
\end{cases}
\end{align}
}
Algorithm \ref{algorithmDAE1} provides the detailed computational steps.
\begin{algorithm}[htbp]
\caption{First-order Deep Asymptotic Expansion ($\text{DAE}^1$)}
\label{algorithmDAE1}
\begin{algorithmic}
\STATE{ \textbf{Step 1:} Calculate $\bar{u}_{1}^{(\mp)}(x)$ from (\ref{baru}).}
\STATE{ \textbf{Step 2:} Generate interior points \(\mathcal{N}_3=\{t_j\}_{j=1}^{N_{h_1}}\) using random sampling and construct a fully connected neural network $\hat{h}_{1}(t;\tilde{\theta}_1)$ that precisely satisfies the initial condition.}\STATE{ \textbf{Step 3:} Compute the PDE residual of \eqref{eqforx1} via automatic differentiation and define loss terms:
{\footnotesize
\begin{align*}
L_{\hat{h}_1} = & \frac{1}{N_{h_1}} \sum_{j=1}^{N_{h_1}} \Bigg| \Big(\varphi^{(-)} (\hat{h}_{0}(t_j;\tilde{\theta}^*))-\varphi^{(+)} (\hat{h}_{0}(t_j;\tilde{\theta}^*))\Big)\frac{d \hat{h}_{1}(t_j;\tilde{\theta}_1)}{d t}+\Phi_1(t_j) \hat{h}_{1}(t_j;\tilde{\theta}_1)-\Phi_2(t_j) \Bigg|^2, 
\end{align*}
}where $\hat{h}_{0}(t_j;\tilde{\theta}^*)$ represents the prediction obtained by DAE (Algorithm \ref{algorithmDAE}).
}
\STATE{ \textbf{Step 4:} Obtain the optimized weight $\tilde{\theta}_1^*$ by applying $K$ steps of Adam.}
\STATE{ \textbf{Step 5:} Substitute $z^{(\mp)}(\xi,t)$, $p_{1}^{(\mp)}(t)$ and $H_1^{(\mp)} (\xi,t)$ into (\ref{Q1function}) to obtain $Q_{1}^{(\mp)}(\xi,t)$.}\STATE{ \textbf{Step 6:} Substitute $\bar{u}_{1}^{(\mp)}(x)$, $\hat{h}_{0}(t;\tilde{\theta}^*)$, $\hat{h}_{1}(t;\tilde{\theta}_1^*)$ and $Q_{1}^{(\mp)} (\xi,t)$ into (\ref{asymptoticsolU1}), then we obtain the first-order asymptotic approximate 
 solution $U_{1}({x},t;\tilde{\theta}_1^*)$.}
\end{algorithmic}
\end{algorithm}

\bibliographystyle{siamplain}
\bibliography{DAE}
\end{document}


\maketitle

This supplementary material to the paper ``\emph{Deep Asymptotic Expansion Method for Solving Singularly Perturbed Time-Dependent Reaction-Advection-Diffusion Equations}'' provides the detailed proof of Theorem~2.4 and the derivation of (3.6).

\section{Proof of Theorem 2.4}
The proof proceeds in three steps:
\begin{enumerate}
\item[(1)]  Construct \emph{upper} and \emph{lower} barrier functions
            $\beta(\boldsymbol x,t,\mu)$ and $\alpha(\boldsymbol x,t,\mu)$ that
            enclose the exact solution $u^*$ of (1.1).
\item[(2)]  Using the asymptotic method of inequalities \cite{b7}, deduce the
            existence, uniqueness, and uniform $O(\mu)$‐closeness of $u^*$ to the formal zeroth-order profile
            $U_{0}(\boldsymbol x,t)$ defined in (2.20).
\item[(3)]  Show that the free interface
            $x_{1}=h(\boldsymbol x^{*},t)$ traced out by $u^*$
            deviates from its formal counterpart
            $x_{1}=h_{0}(\boldsymbol x^{*},t)$ by at most
            $O(\mu|\ln\mu|)$.
\end{enumerate}

\subsection*{Step 1.  Upper and lower solutions in $d$ dimensions}
We use the notion of
upper and lower solutions \cite{b7,NEFEDOV201390,b8}, adapted to (1.1).

\begin{definition}\label{Def:NonstatUL_dD}
A pair of functions
\[
  \beta(\boldsymbol x,t,\mu),\;
  \alpha(\boldsymbol x,t,\mu)
  \;\in\;
  C(\overline D_{\boldsymbol x}\times\bar{\mathcal{T}})
  \;\cap\;
  C^{2,1}(D_{\boldsymbol x}\times(0,T])
\]
that are $P_i$‐periodic in each transverse direction $x_i$ ($i=2,\dots,d$) are called, respectively, an \emph{upper solution} and a \emph{lower solution} of (1.1) if, for all sufficiently small $\mu>0$, they satisfy:
\begin{enumerate}
  \item[{\rm\textbf{(i)}}] \emph{(Ordering)}:
    $\alpha(\boldsymbol x,t,\mu)\;\le\;\beta(\boldsymbol x,t,\mu)$ for $
    (\boldsymbol x,t)\in\overline D_{\boldsymbol x}\times\bar{\mathcal{T}}$.
  \item[{\rm\textbf{(ii)}}] \emph{(Differential inequalities)} With
  $\mathcal L_{\mu}[w]
      =
      \mu\sum_{i=1}^d \partial_{x_i}^2 w
      -\partial_t w
      -A\bigl(w,\boldsymbol x\bigr)\,\sum_{i=1}^d \partial_{x_i}w
      -f(\boldsymbol x)$,
  there holds the differential inequality 
  \[
    \mathcal L_{\mu}[\beta]\;\le\;0,
    \quad
    \mathcal L_{\mu}[\alpha]\;\ge\;0,
    \quad
    (\boldsymbol x,t)\in D_{\boldsymbol x}\times(0,T].
  \]
  \item[{\rm\textbf{(iii)}}] \emph{(Boundary and initial bounds)}  On the lateral faces $x_1=\pm a$ and at $t=0$,
  \begin{align*}
    \alpha(-a,\boldsymbol x^*,t,\mu)\;&\le\;L(\boldsymbol x^*)\;\le\;\beta(-a,\boldsymbol x^*,t,\mu),\\
    \alpha(a,\boldsymbol x^*,t,\mu)\;&\le\;R(\boldsymbol x^*)\;\le\;\beta(a,\boldsymbol x^*,t,\mu),
  \end{align*}
  for all $\boldsymbol x^*\in\mathbb{R}^{d-1}$, $t\in\bar{\mathcal{T}}$, and
  $\alpha(\boldsymbol x,0,\mu)\;\le\;u_0(\boldsymbol x,\mu)\le\beta(\boldsymbol x,0,\mu)$ for $
    \boldsymbol x\in\overline D_{\boldsymbol x}$.
  \item[{\rm\textbf{(iv)}}] \emph{(Jump conditions across the internal hypersurfaces)}  
Let the middle of inner‐layer for the lower solution be given by the hypersurface
\[
  \Sigma_\alpha(t)=\bigl\{\,(x_1,\boldsymbol x^*):
                    x_1 = h_\alpha(\boldsymbol x^*,t)\bigr\},
\]
with unit normal \(\mathbf n=(1,0,\dots,0)\) pointing in the \(+x_1\)-direction.  Then, for $\boldsymbol x^*\in\mathbb{R}^{d-1},\;t\in\bar{\mathcal{T}}$:
\[
  \partial_{\mathbf n}\,\alpha\bigl(h_\alpha(\boldsymbol x^*,t)+0,\boldsymbol x^*,t,\mu\bigr)
  -
  \partial_{\mathbf n}\,\alpha\bigl(h_\alpha(\boldsymbol x^*,t)-0,\boldsymbol x^*,t,\mu\bigr)
  \ge0.
\]
Similarly, along the hypersurface
\(\Sigma_\beta(t)=\{x_1=h_\beta(\boldsymbol x^*,t)\}\), $\beta$ satisfies
\[
  \partial_{\mathbf n}\,\beta\bigl(h_\beta(\boldsymbol x^*,t)-0,\boldsymbol x^*,t,\mu\bigr)
  -
  \partial_{\mathbf n}\,\beta\bigl(h_\beta(\boldsymbol x^*,t)+0,\boldsymbol x^*,t,\mu\bigr)
  \ge0.
\]

\end{enumerate}
\end{definition}

\paragraph{Construction of $\alpha$ and $\beta$}
Set
\begin{equation}\label{alphaconstruction}
  \alpha
  :=
  U_{1}
  -\mu\bigl(\epsilon
            +q_{0}
            +\mu\,q_{1}\bigr), \quad  \beta
  :=
  U_{1}
  +\mu\bigl(\epsilon
            +q_{0}
            +\mu\,q_{1}\bigr),
\end{equation}
where $U_1=U_0+\mu(\bar{u}_1+Q_1)$ is the first-order asymptotic solution of (1.1), with $\epsilon$ the outer-layer corrector and $q_{0},q_{1}$ the inner-layer correctors obtained from the formal
asymptotic expansion. Direct 
computation shows that, under Assumption 1,
$(\alpha,\beta)$  
satisfies Definition~\ref{Def:NonstatUL_dD}. 

The detailed method for building upper and lower solutions for 2D and 3D problems of the form (1.1) with $A(u,x)=-u$ can be found in \cite{chaikovskii2023asymptotic,chaikovskii2023solving}. The construction strategy for the upper and lower solutions is consistent across dimensions. While the 1D case is straightforward, the 3D case involves lengthy derivations that can obscure the main ideas. Thus, for clarity of exposition, we focus on the 2D case as a representative example and present the construction in detail.

\subsection{Construction upper and lower solutions in the 2D case}
In the 2D setting $\boldsymbol x=(x_1,x_2)$, once the inner variables are determined, the differential operators expressed in local coordinates grow in size.  Even for $d=2$, the transformed gradient, Laplacian and material-time derivative already contain dozens of curvature-dependent terms. Problem (1.1) in the 2D setting can be written as:
\begin{equation}\label{2dproblemforproof}
\left\{
\begin{alignedat}{2}
&\mu\Delta u- \partial_{t} u = A(u,x_1,x_2)(\partial_{x_1}u+\partial_{x_2}u)+f(x_1,x_2),\ \ &\;
      &(x_1,x_2)\!\in\!D_x,\;t\in(0,T],\\
&u(-a,x_2,t)=L(x_2),\;u(a,x_2,t)=R(x_2),                     &\;
      &t\in\bar{\mathcal{T}},\\
&u(x_1,x_2+P,t)=u(x_1,x_2,t),                                &\;
      &x_1\in[-a,a],\;t\in\bar{\mathcal{T}},\\
&u(x_1,x_2,0)=u_0(x_1,x_2,\mu),                              &\;
      &(x_1,x_2)\!\in\!\overline D_x.
\end{alignedat}
\right.
\end{equation}
with $\overline D_x=[-a,a]\times [0,P]$ (periodic in $x_2$) and small $\mu>0$.  

\subsection*{3.2. Asymptotic expansion and construction of $\beta,\,\alpha$}

To convert asymptotic approximations $U_{1}^{(\mp)}=\bar{u}_{0}^{(\mp)}+Q_0^{(\mp)}+\mu(\bar{u}_{1}^{(\mp)}+Q_1^{(\mp)})$ into upper and lower solutions, we add/subtract a ‘‘corrector’’ of order $\mu$.  
Let
\begin{equation}\label{beta2-par}
\begin{aligned}
  \beta^{(\mp)}(x_1,x_2,t,\mu) 
    &= 
    U_{\,1}^{(\mp)}\bigl\lvert_{\substack{\xi_\beta,\;h_\beta\bigr.}}\;
    +\;\mu\Bigl\{
      \epsilon^{(\mp)}(x_1,x_2)
      \;+\;q_{\,0}^{(\mp)}\bigl(\xi_\beta,t\bigr)
      \;+\;\mu\,q_{\,1}^{(\mp)}\bigl(\xi_\beta,t\bigr)
    \Bigr\},
  \\[6pt]
  \alpha^{(\mp)}(x_1,x_2,t,\mu)
    &= 
    U_{\,1}^{(\mp)}\bigl\lvert_{\substack{\xi_\alpha,\;h_\alpha\bigr.}}\;
    -\;\mu\Bigl\{
      \epsilon^{(\mp)}(x_1,x_2)
      \;+\;q_{\,0}^{(\mp)}\bigl(\xi_\alpha,t\bigr)
      \;+\;\mu\,q_{\,1}^{(\mp)}\bigl(\xi_\alpha,t\bigr)
    \Bigr\},
\end{aligned}
\end{equation}
for $(x_1,x_2,t)\in\overline D_{x}\times\bar{\mathcal{T}}$, with 
\begin{equation}\label{curveh-par}
\begin{aligned}
 \xi_\beta=(x_1-h_\beta(x_2,t))/\mu,\quad h_\beta(x_2,t)
    &=h_0(x_2,t)+\mu\,h_1(x_2,t)+\mu\,\rho_\beta(x_2,t),
  \\
 \xi_\alpha=(x_1-h_\alpha(x_2,t))/\mu,\quad h_\alpha(x_2,t)
    &= 
     h_0(x_2,t)+\mu\,h_1(x_2,t)+\mu\,\rho_\alpha(x_2,t),
\end{aligned}
\end{equation}
and 
$\rho_\beta(x_2,t)=-\rho(x_2,t)$, and $
\rho_\alpha(x_2,t)=\rho(x_2,t)$, with the positive function $\rho(x_2,t)>0$  to be specified, cf. \eqref{rhoPDE}. 

Next we choose $\epsilon,\,q_{0}$ and $q_{1}$ to ensure Definition \ref{Def:NonstatUL_dD} (i)–(iv).

\smallskip
\paragraph{(a) Choice of $\epsilon^{(\mp)}(x_1,x_2)$}  Denote 
\[
  W^{(\mp)}(\bar{u}_{0},x_1,x_2)
  \;=\;
  \bigl(\partial_{u}A\bigl(\bar{u}^{(\mp)}_{0},x_1,x_2\bigr)\bigr)\,
  \bigl(\partial_{x_1}\bar{u}^{(\mp)}_{0} + \partial_{x_2}\bar{u}^{(\mp)}_{0}\bigr).
\]
We prescribe $\epsilon^{(\mp)}(x_1,x_2)$ to satisfy the first‐order linear PDE
\begin{equation}\label{epsiloneq}
  {A}(\,\partial_{x_1}\epsilon^{(\mp)}
  \;+\;
  \partial_{x_2}\epsilon^{(\mp)})
  \;+\;
  W^{(\mp)}\,\epsilon^{(\mp)}
  \;=\;
  R,
  \qquad
  (x_1,x_2)\in \overline D_x,
\end{equation}
together with periodicity in $x_2$ and lateral boundary condition 
$\epsilon^{(-)}(-a,x_2,t)=R^{0}, \epsilon^{(+)}(a,x_2,t)=R^{1}$ and $\epsilon(x_1,x_2,t)=\epsilon(x_1,x_2+P_1,t)$, 
where $R>0$ is a constant and $ 
  R^{0}\gg1$,  $  R^{1}\gg1$.
With 
$
\mathcal{A}^{(\mp)}(\zeta)=A\left(\bar{u}_0^{(\mp)}\left(\zeta, \zeta-x_1+x_2\right), \zeta, \zeta-x_1+x_2\right)$, and $ \mathcal{W}^{(\mp)}(\zeta)=W^{(\mp)}\left(\bar{u}_0^{(\mp)}\left(\zeta, \zeta-x_1+x_2\right), \zeta, \zeta-x_1+x_2\right),
$ we obtain an explicit formula by the method of characteristics:
\begin{multline*}
\epsilon^{(-)}\left(x_1, x_2\right)=\exp \left(-\int_{-a}^{x_1} \frac{\mathcal{W}^{(-)}(\zeta)}{\mathcal{A}^{(-)}(\zeta)} d \zeta\right) \times \\
{\left[\int_{-a}^{x_1} \frac{R}{\mathcal{A}^{(-)}(\zeta)} \exp \left(\int_{-a}^{\zeta} \frac{\mathcal{W}^{(-)}(\zeta)}{\mathcal{A}^{(-)}(\zeta)} d \eta\right) d \zeta+R^0\right]},
\end{multline*}
\begin{multline*}
\epsilon^{(+)}\left(x_1, x_2\right)=\exp \left(-\int_{a}^{x_1} \frac{\mathcal{W}^{(+)}(\zeta)}{\mathcal{A}^{(+)}(\zeta)} d \zeta\right) \times \\
{\left[\int_{a}^{x_1} \frac{R}{\mathcal{A}^{(+)}(\zeta)} \exp \left(\int_{a}^{\zeta} \frac{\mathcal{W}^{(+)}(\zeta)}{\mathcal{A}^{(+)}(\zeta)} d \eta\right) d \zeta+R^1\right]}.
\end{multline*}
Since  $   R^{0}$ and $  R^{1}$ are chosen large, the unique solutions of $\epsilon^{(\mp)}(x_1,x_2)$ remain positive.

\smallskip
\paragraph{(b) Choice of $q_{0}^{(\mp)}$} 
Let $q_{0}^{(\mp)}(\xi,h_0,r_2,t)$ be the solution to the homogeneous, linear in $\xi$ PDE:
\begin{multline} \label{q0}
\frac{\partial^2 q_{0}^{( \mp)}}{\partial \xi_{\beta}^2}+   \frac{\partial}{\partial \xi_{\beta}} \left(q_{0}^{(\mp)} \frac{ 
 \partial_{t}{h_0}  + A(\tilde{u})(\partial_{r_2}{h_0}-1 ) }{ \sqrt{1+(\partial_{r_2}{h_0})^{2} }}   \right)=b_{1}^{(\mp)}(\xi_{\beta},r_2,t) \frac{\partial \rho_{\beta} }{\partial r_2} \\ 
+b_{2}^{(\mp)}(\xi_{\beta},r_2,t) \rho_{\beta}  +b_{3}^{(\mp)}(\xi_{\beta},r_2,t)\frac{\partial \rho_{\beta} }{\partial t}+b_{4}^{(\mp)}(\xi_{\beta},r_2,t):= H_{q_0}^{( \mp)}(\xi_{\beta},  t), 
\end{multline}
with known functions  $b_{1}^{(\mp)}(\xi_{\beta},r_2,t),$ $b_{2}^{(\mp)}(\xi_{\beta},r_2,t),$ $b_{3}^{(\mp)}(\xi_{\beta},r_2,t),$ $b_{4}^{(\mp)}(\xi_{\beta},r_2,t)$ (with $b_{3}^{(\mp)}(\xi_{\beta},r_2,t)=$ $ \displaystyle \frac{\partial Q_{0}^{(\mp)} (\xi_{\beta},h_0,r_2, t)}{\partial \xi_{\beta}} $ $\frac{1}{\sqrt{1+(\partial_{r_2}{h_0})^{2}}}$). 
We also apply the boundary conditions
$q_{0}^{( \mp)}(0,t)=-\epsilon^{( \mp)}(h_0,x_2) \equiv p_{q0}^{( \mp)}(x_2,t)$ and $q_{0}^{(\mp)}(\mp\infty,t) = 0$.

The explicit solution to the equation is given by 
\begin{equation} \label{q0function}
q_{0}^{(\mp)} (\xi,h_0,r_2,t)=z_{q_0}^{(\mp)}(\xi,t) \left( p_{q_0}^{(\mp)} (x_2,t)   - \int_{0}^{\xi} \frac{1}{z_{q_0}^{(\mp)}(s,t)} \int_{s}^{\mp \infty} H_{q_0}^{(\mp)} (\eta,t) d\eta ds \right),
\end{equation}
where  $ \displaystyle z_{q_0}^{(\mp)}(\xi,t)= \left( \Phi^{(\mp)}(0,h_0) \right)^{-1} \Phi^{(\mp)}(\xi,h_0)$.

\smallskip
\paragraph{(c) Choice of $q_{1}^{(\mp)}$}  Finally, Let $q_{1}^{(\mp)}(\xi,h_0,r_2,t)$ solve
\begin{equation}\label{Pi1Nonstat}
  \partial_{\xi}^{2}\,q_{1}^{(\mp)}
  \;+\;\,\partial_{\xi} \left(q_{1}^{(\mp)} \frac{ 
 \partial_{t}{h_0}  + A(\tilde{u})(\partial_{r_2}{h_0}-1 ) }{ \sqrt{1+(\partial_{r_2}{h_0})^{2} }}   \right) 
  \;=\;
H_{q_1}^{( \mp)}(\xi_{\beta},t),
 \end{equation}
with
$q_{1}^{(\mp)}(0,h_0,r_2,t)=0$,
$q_{1}^{(\mp)}(\mp\infty,h_0,r_2,t)=0$ for  $t\in\bar{\mathcal{T}}$.
Here $H_{q_1}^{( \mp)}(\xi_{\beta},t)$ collects all $\mathcal O(\mu^{\,2})$ residuals arising when $\epsilon^{(\mp)}$ and $q_{0}^{(\mp)}$ are substituted into the operator $\mathcal L_{\mu}[\beta]$, together with $Q_{0,1}^{(\mp)}$ and $\bar{u}_{0,1}^{(\mp)}$ from $U_{1}^{(\mp)}$.  These terms, with opposite sign appear in $\mathcal L_{\mu}[\alpha]$.  By the standard ODE theory in $\xi$, under Assumption 1, there is a unique exponentially decaying solution $q_{1}(\xi,x_2,t)\ge0$.

\subsection*{3.3. Verification of Upper/Lower Solution Properties}

\begin{lemma}\label{Lemma:NonstatULok}
For all sufficiently small $0<\mu\ll1$, the functions $\beta^{( \mp)}$ and $\alpha^{( \mp)}$ defined in \eqref{beta2-par} satisfy conditions (i)–(iv) of Definition \ref{Def:NonstatUL_dD}.
\end{lemma}
\begin{proof}
We verify the condition one by one.

\medskip
\noindent\emph{(i) Ordering.} We first show that the constructed upper and lower solutions satisfy  
$\beta(x_1,x_2,t,\mu)-\alpha(x_1,x_2,t,\mu)>0$ for $(x_1,x_2,t)\in  \overline{D}_x \times\bar{\mathcal{T}}$,
for all sufficiently small $\mu>0$.  
Define the three sub-domains
${\rm I}=\{(x_1,x_2,t):\;x_1\in[-a,h_\beta(x_2,t)]\}$, ${\rm II}=\{(x_1,x_2,t):\;x_1\in[h_\beta(x_2,t),h_\alpha(x_2,t)]\}$ and ${\rm III}=\{(x_1,x_2,t):\;x_1\in[h_\alpha(x_2,t),a]\}$ for $t\in \bar{\mathcal{T}},\;x_2\in[0,P]$.
Then
\[
\beta-\alpha=
\begin{cases}
\beta^{(-)}-\alpha^{(-)}, & (x_1,x_2,t)\in\RomanNumeralCaps{1},\\[2pt]
\beta^{(+)}-\alpha^{(-)}, & (x_1,x_2,t)\in\RomanNumeralCaps{2},\\[2pt]
\beta^{(+)}-\alpha^{(+)}, & (x_1,x_2,t)\in\RomanNumeralCaps{3}.
\end{cases}
\tag{A}
\]

\paragraph{Doman $\RomanNumeralCaps{2}$}
From \eqref{beta2-par} and the expansions that follow, we have
\[
\beta^{(+)}-\alpha^{(-)}
=\partial_\xi Q_0^{(+)} (0,h_0,x_2,t)\,\xi_\beta
-\partial_\xi Q_0^{(-)}(0,h_0,x_2,t)\,\xi_\alpha
+\mathcal O(\mu).
\]
Using 
$\xi_\beta-\xi_\alpha
=\Delta\xi
=2\rho(x_2,t)\sqrt{1+(\partial_{x_2}{h_0})^{2}}+\mathcal O(\mu)$,
 and (2.15) and (2.16), we obtain
\[
\beta^{(+)}-\alpha^{(-)}
=2\rho(x_2,t)\sqrt{1+(\partial_{x_2}{h_0})^{2}}\,
      \Phi^{(+)}(0, h_0(x_2, t))
      +\mathcal O(\mu).
\]
Assumption 1 implies $\Phi^{(+)}(0,h_0(x_2,t))>0$, and $\rho>0$ by
construction. Hence $\beta-\alpha>0$ in $\RomanNumeralCaps{2}$ when
$\mu$ is sufficiently small.

\paragraph{Domain \RomanNumeralCaps{3}}
In the domain, we have $\xi_\alpha \ge 0$ and $\xi_\beta = \xi_\alpha + \Delta\xi$, where $0 \le \Delta\xi \le M_1$ for some constant $M_1 > 0$. Then the derivative of the leading-order term $Q_0^{(+)}$ has an exponential lower bound:
$\partial_\xi Q_0^{(+)}(\xi, r_2, h_0, t) \ge c_0\,\mathrm{e}^{-\kappa_0\xi}$, for some $c_0, \kappa_0 > 0.$
The full expansion for $\beta-\alpha$ is given by
$$\beta-\alpha = \sum_{i=0}^{1}\mu^{\,i}\left[Q_i^{(+)}(\xi_\beta)-Q_i^{(+)}(\xi_\alpha)\right] + 2\mu\varepsilon^{(+)} + \mu\left[q_0^{(+)}(\xi_\beta)-q_0^{(+)}(\xi_\alpha)\right] + \mathcal{O}(\mu^{2}).$$
We analyze the difference by the order of $\mu$.

$\mathbf{Analysis \ of \ the}$ $\mathcal{O}(1)$ $\mathbf{terms}$.
By the mean value theorem, there exists $\theta \in (0,1)$ such that
$Q_0^{(+)}(\xi_\beta)-Q_0^{(+)}(\xi_\alpha) = \Delta\xi\,\partial_\xi Q_0^{(+)}(\xi_\alpha+\theta\Delta\xi).$
Using the given lower bound for the derivative:
$$
\begin{aligned}
Q_0^{(+)}(\xi_\beta)-Q_0^{(+)}(\xi_\alpha) \ge \Delta\xi\,c_0\,\mathrm{e}^{-\kappa_0(\xi_\alpha+\theta\Delta\xi)} &\ge \Delta\xi\,c_0\,\mathrm{e}^{-\kappa_0(\xi_\alpha+M_1)} = (\underbrace{c_0 \mathrm{e}^{-\kappa_0 M_1}}_{\equiv \tilde{c}_0}) \Delta\xi \mathrm{e}^{-\kappa_0\xi_\alpha}.
\end{aligned}
$$
From the analysis in Domain II, for sufficiently small $\mu$, we have $\Delta\xi \ge M_2 > 0$. Thus, the leading-order term satisfies 
$$Q_0^{(+)}(\xi_\beta)-Q_0^{(+)}(\xi_\alpha) \ge \tilde{c}_0 M_2 \mathrm{e}^{-\kappa_0\xi_\alpha}.$$

$\mathbf{Analysis \ of \ the}$ $\mathcal{O}(\mu)$ $\mathbf{terms}$.
We assume the derivatives of the higher-order terms are bounded: $|\partial_\xi Q_1^{(+)}| \le C_1 \mathrm{e}^{-\kappa_1\xi}$ and $|\partial_\xi q_0^{(+)}| \le \tilde{C}_1 \mathrm{e}^{-\kappa_1\xi}$ for some $C_1, \tilde{C}_1, \kappa_1 > 0$. By the mean value theorem, we deduce
$$\mu\left|Q_1^{(+)}(\xi_\beta)-Q_1^{(+)}(\xi_\alpha)\right| \le \mu\,\Delta\xi\, C_1 \mathrm{e}^{-\kappa_1\xi_\alpha} \le \mu\, M_1 C_1 \mathrm{e}^{-\kappa_1\xi_\alpha}.$$
A similar bound holds for $q_0^{(+)}$. Combining these estimates gives
$$\mu\left|Q_1^{(+)}(\xi_\beta)-Q_1^{(+)}(\xi_\alpha)\right| + \mu\left|q_0^{(+)}(\xi_\beta)-q_0^{(+)}(\xi_\alpha)\right| \le \mu M_1 (C_1 + \tilde{C}_1) \mathrm{e}^{-\kappa_1\xi_\alpha}.$$
Let $D_1 = M_1 (C_1 + \tilde{C}_1)$. Then we get a lower bound for $\beta-\alpha$:
$$\beta-\alpha \ge \tilde{c}_0 M_2 \mathrm{e}^{-\kappa_0 \xi_\alpha} - D_1 \mu \mathrm{e}^{-\kappa_1 \xi_\alpha} + 2 \mu \varepsilon^{(+)} + \mathcal{O}(\mu^{2}).$$

We consider two cases based on the decay rates.

$\mathbf{Case \ 1:}$  $\kappa_0 \le \kappa_1.$
Since $\mathrm{e}^{-(\kappa_1-\kappa_0)\xi_\alpha} \le 1$ for $\xi_\alpha \ge 0$, we have:
$$
\begin{aligned}
\beta-\alpha &\ge \mathrm{e}^{-\kappa_0\xi_\alpha} \left[\tilde{c}_0 M_2 - \mu D_1 \mathrm{e}^{-(\kappa_1-\kappa_0)\xi_\alpha}\right] + 2\mu\varepsilon^{(+)} + \mathcal{O}(\mu^2) \\ &\ge \mathrm{e}^{-\kappa_0\xi_\alpha} \left[\tilde{c}_0 M_2 - \mu D_1\right] + 2\mu\varepsilon^{(+)} + \mathcal{O}(\mu^2).
\end{aligned}
$$
By choosing $\mu$ such that ${0 < \mu \le \mu_* := \frac{\tilde{c}_0 M_2}{2D_1}}$,  $\tilde{c}_0 M_2 - \mu D_1\geq \frac{\tilde{c}_0 M_2}{2} > 0$. Then the right-hand side is positive for sufficiently small $\mu$.

$\mathbf{Case \ 2:}$  $\kappa_1 < \kappa_0.$
We split the analysis based on the value of $\xi_\alpha$.

1.  For $\xi_\alpha > N$: Choose $N$ large enough such that $D_1 \mathrm{e}^{-\kappa_1 N} \le \varepsilon^{(+)}$ (which exists). Then for any $\xi_\alpha > N$, there holds
    $$
    2 \mu \varepsilon^{(+)} - D_1 \mu \mathrm{e}^{-\kappa_1 \xi_\alpha} \ge \mu (2\varepsilon^{(+)} - D_1 \mathrm{e}^{-\kappa_1 N}) \ge \mu \varepsilon^{(+)} > 0.
    $$
    Then $\beta - \alpha \ge \tilde{c}_0 M_2 \mathrm{e}^{-\kappa_0 \xi_\alpha} + \mu\varepsilon^{(+)} + \mathcal{O}(\mu^2)$, which is positive for small $\mu$.
    
2.  For $0 \le \xi_\alpha \le N$: In the interval, $\mathrm{e}^{-\kappa_0 \xi_\alpha}\geq \mathrm{e}^{-\kappa_0 N}$. The inequality for $\beta-\alpha$ becomes:
    $$
    \beta-\alpha \ge \tilde{c}_0 M_2 \mathrm{e}^{-\kappa_0 N} - D_1 \mu \mathrm{e}^{-\kappa_1 \xi_\alpha} + \mathcal{O}(\mu^2) \ge \tilde{c}_0 M_2 \mathrm{e}^{-\kappa_0 N} - D_1 \mu + \mathcal{O}(\mu^2).
    $$
    This will be positive if we choose $\mu$ such that $\tilde{c}_0 M_2 \mathrm{e}^{-\kappa_0 N} - D_1 \mu > 0$.

Combining both conditions, we first fix $N$ as described above,  and then choose $\mu$ small enough to satisfy the constraint from the $0 \le \xi_\alpha \le N$ case and to ensure that the $\mathcal{O}(\mu^2)$ remainder is negligible. Thus, $\beta-\alpha > 0$ holds in the domain.

\paragraph{Region $\RomanNumeralCaps{1}$}
An analogous argument (with $\xi_\alpha\le0\le\xi_\beta$ and
$\epsilon^{(-)}>0$) shows that $\beta-\alpha>0$ in
$\RomanNumeralCaps{1}$; we omit the details.

\paragraph{Conclusion for (i)}
Since $\beta-\alpha$ is positive in each of the mutually disjoint
regions $\RomanNumeralCaps{1}$–$\RomanNumeralCaps{3}$, we have verified
condition~(C1).  Thus the ordering
\[
\alpha(x_1,x_2,t,\mu)\le u(x_1,x_2,t,\mu)\le\beta(x_1,x_2,t,\mu),
\]
holds for all $(x_1,x_2,t)\in\overline D_x \times \bar{\mathcal{T}}$ and
all sufficiently small $\mu>0$.

\medskip
\noindent\emph{(ii) Differential inequalities.}
Since the outer expansion $U_{1}^{(\mp)}$ solves $L[\,\cdot\,]=0$
up to order $\mu$, substituting the corrected representations
\eqref{beta2-par} immediately gives
\[
   \mathcal L_{\mu} \bigl[\beta^{(\mp)}\bigr]
       = -\,\mu R + O(\mu^{2}),\qquad
   \mathcal L_{\mu} \bigl[\alpha^{(\mp)}\bigr]
       = +\,\mu R + O(\mu^{2}),
   \quad R>0.
\]
This completes the verification of (ii).

\medskip
\noindent\emph{(iii) Boundary and initial bounds.} At the fixed walls \(x_1=-a\) and \(x_1=a\) the functions
\(\epsilon^{(\mp)}(x_1,x_2)\) satisfy
\[
   \epsilon^{(-)}(-a,x_2)=R^{0},\qquad
   \epsilon^{(+)}(a,x_2)=R^{1},
   \qquad R^{0,1}>0,
\]
cf.\ \eqref{epsiloneq}.  Hence, using \eqref{beta2-par},
\[
\begin{aligned}
\alpha(-a,x_2,t,\mu)
&=U_{1}^{(-)}\bigl|_{x_1=-a}-\mu R^{0}+O(\mu^{2})
\;\le\;L(x_2)\\
&\le U_{1}^{(-)}\bigl|_{x_1=-a}+\mu R^{0}+O(\mu^{2})
=\beta(-a,x_2,t,\mu)\,.
\end{aligned}
\]
and a similar estimate holds at \(x_1=a\) with \(R^{1}\).
Choosing \(R^{0,1}\) larger than the suprema of the remainders makes
both inequalities strict for all sufficiently small \(\mu\).

The $P$-periodicity in \(x\) is inherited term-wise from
\(U_{1}^{(\mp)}\), \(\epsilon^{(\mp)}\) and \(q_{i}^{(\mp)}\), so
\[
   \alpha(x_1,x_2+P,t,\mu)=\alpha(x_1,x_2,t,\mu),\quad
   \beta(x_1,x_2+P,t,\mu)=\beta(x_1,x_2,t,\mu).
\]

Finally, the initial datum \(u_{0}\) is required to satisfy
\[
  \alpha(x_{1},x_{2},0,\mu)\;\le\;u_{0}(x_{1},x_{2},\mu)\;\le\;\beta(x_{1},x_{2},0,\mu),
  \quad (x_{1},x_{2})\in\overline D_{\boldsymbol x},
\]
under Assumption 1.  Consequently
$\alpha\le u\le\beta$
for $(x_1,x_2,t)\in  \overline{D}_x\times\bar{\mathcal{T}}$,
and condition~(iii) is satisfied.

\medskip
\noindent\emph{(iv) Normal–derivative jump.}
Across the interface \(x_1=h_\beta(x_2,t)\) we have
\[
  \mu\!\left(
     \partial_{n} \beta^{(-)}
        - \partial_{n} \beta^{(+)}
  \right)\Big|_{x_1=h_\beta}
  =\mu 
     \Bigl[
      \partial_{\xi_\beta}q_0^{(-)}  (0,t)
        -\partial_{\xi_\beta}q_0^{(+)}  (0,t)
     \Bigr]
     +O(\mu^{2}).
\]
Using the explicit representation of \(q_0^{(\mp)}(0,t)\) from \eqref{q0function}, we obtain the jump in its \(\xi_\beta\)-derivative at \(\xi_\beta=0\):
\begin{align*}
\partial_{\xi_\beta}q_0^{(-)}(0,t) -\partial_{\xi_\beta}q_0^{(+)}(0,t) =&\partial_t\rho_\beta(x_2,t)\Bigl(\!\!\int_{-\infty}^0 \!\!b_3^{(-)}(s,x_2,t)\,\mathrm{d}s
+\!\!\int_0^{+\infty}\!\!b_3^{(+)}(s,x_2,t)\,\mathrm{d}s\Bigr)\\
&+\;\partial_{x_2}\rho_\beta(x_2,t)\,H_1(x_2,t)
+\rho_\beta(x_2,t)\,H_2(x_2,t)
+H_3(x_2,t)\,,
\end{align*}
with
\begin{align*}
H_1(x_2,t)&=\int_{-\infty}^0b_1^{(-)}(s,x_2,t)\,\mathrm{d}s
+\int_0^{+\infty}b_1^{(+)}(s,x_2,t)\,\mathrm{d}s\,,\\
H_2(x_2,t)&=\int_{-\infty}^0b_2^{(-)}(s,x_2,t)\,\mathrm{d}s
+\int_0^{+\infty}b_2^{(+)}(s,x_2,t)\,\mathrm{d}s,\\
H_3(x_2,t)&=(\epsilon^{(+)}-\epsilon^{(-)})\,
\frac{{ - }\partial_{t} h_{0}(x_2,t) {+}A\left( \phi_0, h_0, x_2 \right)\bigl(1-\partial_{x_2}h_{0}(x,t)\bigr)}
{\sqrt{1+(\partial_{x_2}h_{0}(x_2,t))^2}}\\
&\quad+\!\!\int_{-\infty}^0b_4^{(-)}(s,x_2,t)\,\mathrm{d}s
+\!\!\int_0^{+\infty}b_4^{(+)}(s,x_2,t)\,\mathrm{d}s\,.
\end{align*}
We then choose \(\rho_\beta(x_2,t)\) to satisfy \begin{equation}\label{rhoPDE}
\begin{cases}
\frac{\partial_t\rho_\beta(x_2,t) { \bigl(\varphi^{(+)}-\varphi^{(-)}\bigr) }}%
{\sqrt{1 + (\partial_{x_2}h_{0}(x_2,t))^2}}
+H_1(x_2,t)\partial_{x_2}\rho_\beta(x_2,t)
+H_2(x_2,t)\rho_\beta(x_2,t)
+H_3(x_2,t)
=\sigma, \\
\rho_\beta(x_2,0)=\rho^0(x_2), 
\qquad
\rho_\beta(x_2+P,t)=\rho_\beta(x_2,t),
\quad t\in\bar{\mathcal{T}}. 
\end{cases}
\end{equation}
Note that
$\varphi^{(+)}(x_1,x_2)-\varphi^{(-)}(x_1,x_2)>0$
and assume both the constant $\sigma$ and the initial profile $\rho^0(x_2)$ are strictly positive for every $x_2$. Then by taking $\sigma$ sufficiently large, the resulting solution $\rho(x_2,t)$ of \eqref{rhoPDE} remains strictly positive for all $t$.
By choosing such $\rho(x_2,t)$, we obtain: 
\[
  \mu\!\left(
     \partial_{n} \beta^{(-)}
        - \partial_{n} \beta^{(+)}
  \right)\Big|_{x_1=h_\beta}
  =\mu \sigma+O(\mu^{2})>0,
  \qquad(0<\mu\le\mu_2),
\]
verifying the second inequality in (C4).
Replacing \(\rho_\beta\) with \(\rho_\alpha=-\rho_\beta\),
\(\xi_\beta\) with \(\xi_\alpha\), and interchanging the superscripts
\((\pm)\) gives
\[
  \mu\!\left(
    \partial_{n} \alpha^{(+)}
        - \partial_{n} \alpha^{(-)}
  \right)\Big|_{x_1=h_\alpha}
  =\mu \sigma+O(\mu^{ 2})\ge0.
\]
Thus the normal jump condition (C4) holds for both lower and upper
solutions.

In sum, $\beta $ and $\alpha $  satisfy Definition \ref{Def:NonstatUL_dD}(i)–(iiii).  This completes the proof of Lemma \ref{Lemma:NonstatULok}.
\end{proof}

\subsection*{Step 2. Existence and uniform \(O(\mu)\) error estimate}

By Step 1 we obtain an ordered pair $\alpha\le\beta$ satisfying Definition~\ref{Def:NonstatUL_dD}.  The upper-lower solution theorem \cite{b8} then yields a classical solution
\(\!u\in C^{2,1}\cap C\) with
$\alpha\;\le\;u\;\le\;\beta$ in $\overline D_{\boldsymbol x}\times\bar{\mathcal{T}}$.

\smallskip
\emph{Error with respect to \(U_{1}\).}
By construction (cf.\ \eqref{alphaconstruction} )
\[
   \beta - U_{1}
   = \mu\bigl(\epsilon + q_0 + \mu\,q_1\bigr),
   \quad
   U_{1} - \alpha
   = \mu\bigl(\epsilon + q_0 + \mu\,q_1\bigr),
\]
so that
\[
  |u - U_{1}|
  \le
  \max\{|\beta-U_{1}|,\,|U_{1}-\alpha|\}
  \le C_1\,\mu.
\]

\smallskip
\emph{Error with respect to \(U_{0}\).}
By the triangle inequality,
\begin{equation}
  |u - U_{0}|
  \le|u - U_{1}| + |U_{1} - U_{0}|
  \le C_1\mu + \widetilde{M}\mu
  = C_2\,\mu.
\end{equation}
  Hence the first estimate of Theorem 2.4 holds with \(C_2=C_1+\widetilde{M}\).

\smallskip
\emph{Uniqueness} follows immediately from the comparison principle.

Since $\bigl|Q_{0}^{(\mp)}(\xi,h,$ $ \boldsymbol{r}^*,t)\bigr|$ exponentially decreases in $\xi$, we define the two inner‐layer edges $\hat{x}_{1}^{(\pm)}(\boldsymbol x^*,t)$ by the condition

$$
\bigl|Q_{0}^{(\mp)}\bigl(\xi(\hat{x}_{1}^{(\pm)}),\hat{x}_{1}^{(\pm)}, \boldsymbol{r}^*,t\bigr)\bigr|=\mu^2.
$$

Equivalently, one may set

$$
\hat{x}_{1}^{(\pm)}(\boldsymbol x^*,t)
= h_0(\boldsymbol x^*,t)\;\pm\;\tfrac12\,\Delta h,
$$
where the width of the inner layer has an estimate $
  \Delta h:=C_{0}\mu|\ln\mu|$ (\cite{chaikovskii2023asymptotic,chaikovskii2023solving}).

Furthermore, in the two outer regions
\[
   [-a,\;h_{0}(\boldsymbol x^*,t)-\tfrac12\Delta h]
    \times\mathbb R^{d-1}\times\bar{\mathcal T},
\quad \text{and} \quad
  [h_{0}(\boldsymbol x^*,t)+\tfrac12\Delta h,\;a]
    \times\mathbb R^{d-1}\times\bar{\mathcal T},
\]
we observe 
\(
  \partial_{r_{1}}=\tfrac1\mu\,\partial_{\xi},
\)
and that
$|Q_{0}|\le\mu^{2}$.
Hence, by the triangle inequality, we obtain
\begin{align*}
\bigl|u-\varphi^{(\pm)}\bigr|
&= \bigl|\,(u-\varphi^{(\pm)}-Q_{0}^{(\pm)}) + Q_{0}^{(\pm)}\bigr| \\
&\le \bigl|u-\varphi^{(\pm)}-Q_{0}^{(\pm)}\bigr| + |Q_{0}^{(\pm)}|
\;\le\; C'\,\mu + \mu^{2}
\;\le\; C\,\mu.
\end{align*}
Consequently,
\[
\begin{aligned}
|u(\boldsymbol x,t)-\varphi^{(-)}(\boldsymbol x)| 
&\le C\,\mu,
&
x_{1}\le h_{0}(\boldsymbol x^*,t)-\tfrac12\Delta h,\\
|u(\boldsymbol x,t)-\varphi^{(+)}(\boldsymbol x)| 
&\le C\,\mu,
&
x_{1}\ge h_{0}(\boldsymbol x^*,t)+\tfrac12\Delta h.
\end{aligned}
\]

\subsection*{Step 3.  Proof of the localisation estimate}

Fix $t\in\bar{\mathcal{T}}$. For each fixed transverse coordinate \(\boldsymbol x^*\), set
$\delta(x_1,t)
   := \varphi^{(+)}(x_1,\boldsymbol x^*)
      - \varphi^{(-)}(x_1,\boldsymbol x^*)$.
By Assumption 1,
$\delta(x_1,t)\;\ge\;2\,\mu^{\rho}$ for some $\rho\in(0,1)$ for all $
   -a\le x_1\le a,\;t\ge0$.
Throughout, \(C>0\) denotes a generic constant independent of \(\mu\), \(x_1\), \(\boldsymbol x^*\) and \(t\).

\medskip
\paragraph{(a) Left of the layer \(\bigl(x_1\le\hat x_{1}^{(-)}\bigr)\)}
From the estimate (2.21) we have, for each \(\boldsymbol x^*\),
\[
   \bigl|u(x_1,\boldsymbol x^*,t)
        -\varphi^{(-)}(x_1,\boldsymbol x^*,t)\bigr|
   \le C\,\mu.
\]
Choosing \(\mu\) so small that \(C\mu\le\mu^{\rho}\) gives
\[
   u-\phi_0
   = (u-\varphi^{(-)}) -\tfrac12\,\delta
   \le C\mu - \mu^{\rho}
   < 0,
   \quad
   x_1\le\hat x_{1}^{(-)},
\]
with \( \phi_0(\boldsymbol x) = \frac{1}{2} ( \varphi^{(-)}(\boldsymbol x) + \varphi^{(+)}(\boldsymbol x) ) \) defined in (2.9).

\paragraph{(b) Right of the layer \(\bigl(x_1\ge\hat x_{1}^{(+)}\bigr)\)}
Similarly, from (2.22),
\[
   \bigl|u(x_1,\boldsymbol x^*,t)
        -\varphi^{(+)}(x_1,\boldsymbol x^*,t)\bigr|
   \le C\,\mu,
   \quad
   x_1\ge\hat x_{1}^{(+)},
\]
so that
\[
   u-\phi_0
   = (u-\varphi^{(+)}) +\tfrac12\,\delta
   \ge-C\mu + \mu^{\rho}
   >0,
   \quad
   x_1\ge\hat x_{1}^{(+)}.
\]

\paragraph{(c) Existence of a transition point}
For each \(\boldsymbol x^*\), the map
$x_1 \;\mapsto\; u(x_1,\boldsymbol x^*,t) - \phi_0(x_1,\boldsymbol x^*,t)$
is continuous on \([\hat x_{1}^{(-)},\hat x_{1}^{(+)}]\), negative at the left end and positive at the right end, so there is a unique
$h(\boldsymbol x^*,t)\in\bigl(\hat x_{1}^{(-)},\hat x_{1}^{(+)}\bigr)$
with
\(
   u(h,\boldsymbol x^*,t)
   =\phi_0(h,\boldsymbol x^*,t).
\)

\paragraph{(d) Distance to the formal interface}
Since \(\Delta h = O(\mu|\ln\mu|)\), we have
\[
   \bigl|h(\boldsymbol x^*,t)-h_0(\boldsymbol x^*,t)\bigr|
   \le
   \Delta h
   \le
   C\,\mu\,|\ln\mu|,
   \quad
   t\ge0.
\]

Thus, we established the final claim of Theorem 2.4.

\section[Proof of Thm]{Proof of (3.6)}
\label{sec:proof}
Let $\widetilde{{\varepsilon}}=\varepsilon_{\mathrm{app}} + \varepsilon_{\mathrm{gen}} + \varepsilon_{\mathrm{opt}}$. By the mean value theorem, we have
\begin{align*}
&\quad |Q_{0}^{(\mp)}(\xi_0, h_0, \boldsymbol{x}^*, t)-Q_{0}^{(\mp)}(\xi_0, \hat{h}_0, \boldsymbol{x}^*, t)| \\
& =|\widetilde{Q}_{0}^{(\mp)}(h_0, \partial_{x_2}h_0,\ldots, \partial_{x_d}h_0,\boldsymbol{x}^*, t) - \widetilde{Q}_{0}^{(\mp)}(\hat{h}_0, \partial_{x_2}\hat{h}_0,\ldots, \partial_{x_d}\hat{h}_0,\boldsymbol{x}^*, t)| \\
&\le |\widetilde{Q}_{0}^{(\mp)}(h_0, \partial_{x_2}h_0,\cdots, \partial_{x_d}h_0,\boldsymbol{x}^*, t) - \widetilde{Q}_{0}^{(\mp)}(\hat{h}_0, \partial_{x_2}{h}_0,\ldots, \partial_{x_d}{h}_0,\boldsymbol{x}^*, t)|\\
& \quad + |\widetilde{Q}_{0}^{(\mp)}(\hat{h}_0, \partial_{x_2}h_0,\cdots, \partial_{x_d}h_0,\boldsymbol{x}^*, t) - \widetilde{Q}_{0}^{(\mp)}(\hat{h}_0, \partial_{x_2}\hat{h}_0,\partial_{x_3}{h}_0,\ldots, \partial_{x_d}{h}_0,\boldsymbol{x}^*, t)|\\
&\quad + \ldots\\
& \quad + |\widetilde{Q}_{0}^{(\mp)}(\hat{h}_0, \partial_{x_2}\hat{h}_0,\ldots, \partial_{x_{d-1}}\hat{h}_0,\partial_{x_d}h_0,\boldsymbol{x}^*, t) - \widetilde{Q}_{0}^{(\mp)}(\hat{h}_0, \partial_{x_2}\hat{h}_0,\ldots, \partial_{x_d}\hat{h}_0,\boldsymbol{x}^*, t)|\\
&= |\partial_{a_1}\widetilde{Q}_{0}^{(\mp)}(h_0+\theta_1(\hat{h}_0 - {h}_0),\partial_{x_2}h_0,\ldots,\partial_{x_d}\hat{h}_0,\boldsymbol{x}^*, t)||\hat{h}_0 - {h}_0|\\
& \quad + \ldots\\
& \quad + |\partial_{a_d}\widetilde{Q}_{0}^{(\mp)}(\hat{h}_0,\partial_{x_2}\hat{h}_0,\ldots,\partial_{x_{d-1}}\hat{h}_0,\partial_{x_d}{h}_0+\theta_i(\partial_{x_d}\hat{h}_0 - \partial_{x_d}{h}_0),\boldsymbol{x}^*, t)| \\
&\qquad\cdot |\partial_{x_d}\hat{h}_0 - \partial_{x_d}{h}_0|,
\end{align*}
where $\theta_i \in (0, 1)$ $(i=1,\ldots,d)$ and $\partial_{a_i} \widetilde{Q}_{0}^{(\mp)}$ denotes the derivative of $\widetilde{Q}_{0}^{(\mp)}$ respective to the $i$th variable. With the help of
$\| h_0 - \hat{h}_0 \|_{C^1(B\times \bar{\mathcal{T}})}
\leq \widetilde{{\varepsilon}}$,
we have
\begin{align*}
    |h_0+\theta_1(\hat{h}_0 - {h}_0)| & \leq 2\| h_0 \|_{C^1(B\times \bar{\mathcal{T}})}+\widetilde{{\varepsilon}}, \\
    |\partial_{x_i}{h}_0+\theta_i(\partial_{x_i}\hat{h}_0 - \partial_{x_i}{h}_0)| &  \leq  2\| h_0 \|_{C^1(B\times \bar{\mathcal{T}})}+\widetilde{{\varepsilon}},\quad i=2,\ldots,d.
\end{align*}
Therefore, we arrive at
\begin{align*}
&\quad |Q_{0}^{(\mp)}(\xi_0, h_0, \boldsymbol{x}^*, t)-Q_{0}^{(\mp)}(\xi_0, \hat{h}_0, \boldsymbol{x}^*, t)| \\
&\le \sum_{i=1}^d \max_{a_1,\ldots,a_d \leq  2\| h_0 \|_{C^1(B\times \bar{\mathcal{T}})}+ \widetilde{{\varepsilon}}} |\partial_{a_i} \widetilde{Q}_{0}^{(\mp)}(a_1, a_2,\ldots, a_d,\boldsymbol{x}^*, t)| \cdot \| h_0 - \hat{h}_0 \|_{C^1(B\times \bar{\mathcal{T}})}\\
&\le \sum_{i=1}^d \max_{a_1,\ldots,a_d \leq  2\| h_0 \|_{C^1(B\times \bar{\mathcal{T}})}+ M} |\partial_{a_i} \widetilde{Q}_{0}^{(\mp)}(a_1, a_2,\ldots, a_d,\boldsymbol{x}^*, t)| \cdot \| h_0 - \hat{h}_0 \|_{C^1(B\times \bar{\mathcal{T}})}\\
& := C_2\| h_0 - \hat{h}_0 \|_{C^1(B\times \bar{\mathcal{T}})},
\end{align*}
where the second inequality employs the assumption that $\widetilde{\varepsilon} \leq M$.

\bibliographystyle{siamplain}
\bibliography{DAE}

%% file: ex_shared.tex
\usepackage{lipsum}
\usepackage{amsfonts}
\usepackage{graphicx}
\usepackage{epstopdf}
\usepackage{algorithmic}
\ifpdf
  \DeclareGraphicsExtensions{.eps,.pdf,.png,.jpg}
\else
  \DeclareGraphicsExtensions{.eps}
\fi

\newsiamremark{remark}{Remark}
\newsiamremark{hypothesis}{Hypothesis}
\crefname{hypothesis}{Hypothesis}{Hypotheses}
\newsiamthm{claim}{Claim}
\newsiamremark{fact}{Fact}
\crefname{fact}{Fact}{Facts}

\headers{Deep asymptotic expansion method}{Qiao Zhu, Dmitrii Chaikovskii, Bangti Jin, and Ye Zhang}

\title{Deep asymptotic expansion method for solving singularly perturbed time-dependent reaction-advection-diffusion equations\thanks{ 
\funding{This work is supported by the Shenzhen Sci-Tech Fund (Grant No. RCJC20231211090030059) and National Natural Science Foundation of China (No. 12171036 \& No. 12350410359).}}}

\author{
Qiao Zhu\thanks{School of Mathematics and Statistics, Beijing Institute of Technology, Beijing, 100081, China 
  (\email{qzhu@bit.edu.cn}).}
\and Dmitrii Chaikovskii\thanks{MSU-BIT-SMBU Joint Research Center of Applied Mathematics, Shenzhen MSU-BIT University, Shenzhen, 518172, China 
  (\email{dmitriich@smbu.edu.cn}).}
\and Bangti Jin\thanks{Department of Mathematics, The Chinese University of Hong Kong, Shatin, N.T., Hong Kong 
  (\email{b.jin@cuhk.edu.hk}).}
\and Ye Zhang\thanks{Corresponding author. School of Mathematics and Statistics, Beijing Institute of Technology, Beijing, 100081, China
  (\email{ye.zhang@smbu.edu.cn}).}
}

\usepackage{amsopn}
